\documentclass[twoside,11pt]{article}

\usepackage{jmlr2e}
\usepackage{amsfonts}
\usepackage{amsmath}
\usepackage{color}
\usepackage{graphicx}
\usepackage{caption}
\usepackage{verbatim}
\usepackage{subcaption}
\usepackage{rotating}
\usepackage{multirow}
\newcommand{\qed}{{\unskip\nobreak\hfil\penalty50\hskip2em\vadjust{}
            \nobreak\hfil$\Box$\parfillskip=0pt\finalhyphendemerits=0\par}}

\newtheorem{thm}[theorem]{Theorem}
\newtheorem{cor}[theorem]{Corollary}
\newcommand{\bed}{\begin{definition}}
\newcommand{\eed}{\end{definition}}
\newcommand{\beq}{\begin{equation}}
\newcommand{\eeq}{\end{equation}}

\newcommand{\eps}{\epsilon}

\newcommand{\bitem}{\begin{itemize}}
\newcommand{\eitem}{\end{itemize}}

\newcommand{\goto}{\rightarrow}

\newcommand{\margmin}{\mathrm{argmin}}

\newcommand{\beqn}{\begin{equation}}
\newcommand{\eeqn}{\end{equation}}
\newcommand{\balign}{\begin{align}}
\newcommand{\ealign}{\end{align}}

\newcommand{\bphi}{\bar{\Phi}}

\newcommand{\rank}{\mathrm{rank}}
\newcommand{\sgn}{\mathrm{sgn}}

\newcommand{\ty}{\tilde{Y}}

\newcommand{\cg}{{\cal G}}

\newcommand{\call}{{\cal I}}
\newcommand{\hb}{\hat{\beta}}
\newcommand{\hamm}{\mathrm{Hamm}}

\newcommand{\bel}{\begin{eqnarray}\label}
\newcommand{\eel}{\end{eqnarray}}
\newcommand{\bes}{\begin{eqnarray*}}
\newcommand{\ees}{\end{eqnarray*}}

\def\sgn{\hbox{\rm sgn}}

\newcommand{\bi}{\begin{itemize}}
\newcommand{\ei}{\end{itemize}}
\newcommand{\hd}{\hat{D}}
\newcommand{\he}{\hat{F}}
\newcommand{\tq}{{\cal Q}}
\def\R{\mathbb{R}}

\jmlrheading{15}{2014}{xx-xx}{05/12; Revised 01/14}{xx/14}{Jiashun Jin, Cun-Hui Zhang, and Qi Zhang}
\ShortHeadings{Optimality of Graphlet Screening  in High Dimensional Variable Selection}{ Jin, Zhang and Zhang}
\firstpageno{1}

\begin{document}
\title{Optimality of Graphlet Screening  in High Dimensional Variable Selection}


\author{\name Jiashun Jin \email jiashun@stat.cmu.edu \\
       \addr Department of Statistics\\
       Carnegie Mellon University\\
       Pittsburgh, PA 15213, USA
       \AND
       \name Cun-Hui Zhang \email cunhui@stat.rutgers.edu \\
       \addr Department of Statistics\\
       Rutgers University  \\
       Piscataway,   NJ 08854, USA
       \AND
       \name Qi Zhang \email  qizhang@stat.wisc.edu   \\
       \addr Department of Biostatistics $\&$ Medical Informatics  \\
University of Wisconsin-Madison \\
Madison, WI 53705, USA
       }

\editor{Bin Yu}
\maketitle

 \begin{abstract}
Consider a linear  model $Y = X \beta + \sigma z$,   where $X$ has $n$ rows and $p$ columns and $z \sim N(0, I_n)$.  We assume both $p$ and $n$ are large, including the case of  $p \gg n$.
The unknown signal vector $\beta$ is assumed to be sparse in the sense that
only a small fraction of its components is nonzero.   The goal is to identify such nonzero coordinates (i.e., variable selection).

We are primarily interested in the regime where signals are both {\it rare and weak} so that
successful variable selection is challenging but is still possible.
We assume the Gram matrix $G = X'X$ is sparse in the sense that
each row  has relatively few large entries
(diagonals of $G$ are normalized to $1$).
The sparsity of $G$ naturally induces the sparsity of the so-called {\it Graph of Strong Dependence} (GOSD).     The key insight is that there is an interesting interplay between the   signal sparsity and graph sparsity: in a broad context,  the signals decompose into many  small-size components of  GOSD  that are disconnected to each other.

We propose {\it Graphlet Screening}  for variable selection.  This is a two-step Screen and Clean procedure, where in the first step, we screen subgraphs of GOSD with sequential $\chi^2$-tests, and in the second step, we clean with penalized MLE. The main methodological innovation is to use GOSD to guide both the screening and cleaning processes.

For any variable selection procedure $\hat{\beta}$, we measure its performance by the Hamming distance between the sign vectors of $\hat{\beta}$ and $\beta$, and assess the optimality by the minimax Hamming distance. Compared with more stringent criterions such as exact support recovery
or oracle property,
which demand strong signals, the Hamming distance criterion is more appropriate for weak signals since it naturally allows a small fraction of errors.

We show that in a broad class of situations,  Graphlet Screening achieves the optimal rate of convergence in terms of  the Hamming distance.  Unlike Graphlet Screening,
well-known procedures  such as the $L^0/L^1$-penalization methods do not utilize local graphic structure for variable selection,  so they generally do not achieve the optimal rate of convergence, even in very simple settings and even if the tuning parameters are ideally set.

The the presented algorithm is implemented as R-CRAN package {\it ScreenClean} and in \emph{matlab} (available at http://www.stat.cmu.edu/$\sim$jiashun/Research/software/GS-matlab/).
\end{abstract}

\begin{quote}
\end{quote}

\section{Introduction}
\label{sec:Intro}
Consider a linear regression model
\begin{equation} \label{model}
Y = X \beta + \sigma z, \qquad  X = X_{n, p}, \qquad   z  \sim N(0, I_n).
\end{equation}
We write
\begin{equation} \label{columnX}
X   =  [x_1, x_2, \ldots, x_p],   \qquad \mbox{and} \qquad X' = [X_1, X_2, \ldots, X_n],
\end{equation}
so that $x_j$ is the $j$-th design vector and $X_i$ is the $i$-th sample.
Motivated by the recent interest in `Big Data',   we assume both $p$ and $n$ are large but  $p \geq n$ (though this should not be taken as a restriction).
The  vector $\beta$ is unknown to us, but is presumably {\it sparse} in the sense that  only a small proportion of its entries  is nonzero. Calling a nonzero entry of $\beta$ a signal,  the main interest of this paper  is to  identify all signals (i.e.,  variable selection).

Variable selection is one of the most studied problem in statistics.
However,  there are important regimes
where our understanding  is very limited.

One of such regimes is the {\it rare and weak} regime, where the signals are both {\it rare} (or sparse) and  {\it individually weak}. Rare and weak signals are frequently found in research areas such as Genome-wide Association Study (GWAS) or next generation sequencing.
Unfortunately, despite urgent demand in applications,   the literature of variable selection has  been focused on the regime where the signals are rare but {\it individually strong}.  This motivates a revisit to variable selection, focusing on the rare and weak regime.

For variable selection in this regime,  we need new methods and  new theoretical frameworks. In particular, we need a loss function that is appropriate for
rare and weak signals to evaluate the optimality. In the literature,
given a variable selection procedure $\hat{\beta}$,
we usually use the probability of exact recovery $P(\sgn(\hat{\beta}) \neq  \sgn(\beta))$
as the measure of loss \citep{FanLi}; $\sgn(\hat{\beta})$ and $\sgn(\beta)$ are the sign vectors of $\hat{\beta}$ and $\beta$ respectively.
In the rare and weak regime,
the signals are so rare and weak that exact recovery is impossible,
and the Hamming distance between  $\sgn(\hat{\beta})$ and $\sgn(\beta)$
is a more appropriate measure of loss.

Our focus on the rare and weak regime and  the Hamming distance loss provides  new perspectives to variable selection, in methods and in theory.

Throughout this paper, we assume the diagonals of the Gram matrix
\begin{equation} \label{diag}
G  = X'X
\end{equation}
are normalized to $1$ (and approximately $1$ in the random design model), instead of $n$ as
often used in the literature.
The difference between two normalizations is non-essential,
but the signal vector $\beta$ are different by a factor of $n^{1/2}$. 

We also assume the Gram matrix
$G$ is `sparse'  (aka. graph sparsity)  in the sense that each of its rows has relatively few large entries.   Signal sparsity and graph sparsity can be simultaneously found in the following application areas.
\begin{itemize}
\item {\it Compressive sensing}.   We are interested in a very high dimensional sparse vector $\beta$.
The goal is  to store or transmit $n$ linear functionals of $\beta$ and then reconstruct it.
For $1 \leq i \leq n$,   we choose a $p$-dimensional coefficient vector $X_i$ and  observe
$Y_i = X_i' \beta + \sigma z_i$ with an error  $\sigma z_i$.
The so-called Gaussian design is often considered \citep{Donoho2006a, Donoho2006b,  Nowak2007}, where $X_i \stackrel{iid}{\sim} N(0, \Omega/n)$ and $\Omega$ is  sparse; the
  sparsity of $\Omega$ induces the sparsity  of $G = X' X$.
\item {\it Genetic Regulatory Network (GRN)}.   For $1 \leq i \leq n$,  $W_i =  (W_i(1), \ldots, W_i(p))'$ represents the expression level of $p$ different genes of   the $i$-th patient.   Approximately, $W_i \stackrel{iid}{\sim} N(\alpha, \Sigma)$, where the contrast mean vector $\alpha$ is sparse reflecting that only few genes are differentially expressed between a  normal patient and a diseased one \citep{Pengetal2009}.    Frequently,   the concentration matrix $\Omega =  \Sigma^{-1}$ is believed to be  sparse, and can be effectively estimated in some cases (e.g., \citet{Bickel2008a} and \citet{CaiLuo2010}), or can be assumed as known in others, with the so-called  ``data about data"   available \citep{LiLi}.
Let $\hat{\Omega}$ be a positive-definite estimate of $\Omega$,  the setting  can be re-formulated as the linear model    $(\hat{\Omega})^{1/2}  Y  \approx  \Omega^{1/2} Y  \sim N(\Omega^{1/2}\beta, I_p)$,  where $\beta = \sqrt{n} \alpha$ and
the Gram matrix $G \approx \Omega$,   and both are sparse.
\end{itemize}
Other examples can be found in Computer Security \citep{JiJin} and Factor Analysis \citep{Hanetal2011}.


The sparse Gram matrix $G$ induces a sparse graph  which we call
the {\it Graph of Strong Dependence} (GOSD), denoted by  $\cg = (V, E)$, where $V = \{1,2, \ldots, p\}$ and there is an edge between nodes
$i$ and $j$ if and only the design vectors $x_i$ and $x_j$ are {\it strongly correlated}.   Let
\begin{equation}  \label{DefineS}
S = S(\beta) = \{1 \leq j \leq p: \beta_j \neq 0\}
\end{equation}
be the support of $\beta$ and $\cg_S$ be the subgraph of $\cg$  formed by all nodes in $S$.
The key insight is that,  there is an interesting interaction between signal sparsity and graph sparsity, which yields the subgraph $\cg_S$ decomposable:  $\cg_S$ splits into many ``graphlet"; each  ``graphlet" is a small-size component and different components are not connected (in $\cg_S$).

While we can always decompose $\cg_S$ in this way, our emphasis in this paper is that,  in many cases,  the maximum size of the graphlets is small; see Lemma \ref{lemma:sparsegraph} and related discussions.

The decomposability of $\cg_S$ motivates a new approach to variable selection, which we call {\it Graphlet Screening} (GS). GS is a Screen and Clean method \citep{Wasserman}. In the screening stage,  we use  multivariate screening to identify candidates for all the graphlets.
Let $\hat{S}$ be all the nodes that survived the screening, and let $\cg_{\hat{S}}$ be the subgraph of GOSD formed by all nodes in $\hat{S}$.
 Although $\hat{S}$ is expected to be somewhat larger than $S$, the subgraph $\cg_{\hat{S}}$ is still likely to resemble $\cg_S$ in structure in the sense that it, too,
splits into many small-size disconnected components.  We then clean each component  separately to remove  false positives.

The objective of the paper is two-fold.
\begin{itemize}
\item To propose a ``fundamentally correct" solution in the rare and weak paradigm along with a computationally fast algorithm for the solution.
\item To show that GS achieves the optimal rate of convergence in terms of the Hamming distance, and achieves the optimal phase diagram for variable selection.
\end{itemize}
Phase diagram can be viewed as an optimality criterion which is especially appropriate for rare and weak signals. See  \citet{DonohoJin2004} and \citet{PNAS} for example.

In the settings we consider,  most popular approaches  are  {\it not} rate optimal; we explain this in Sections  \ref{subsec:subsetselection}-\ref{subsec:BMSC}.
In Section \ref{subsec:interplay}, we explain the basic idea of GS and why it works.
\subsection{Non-optimality of  the $L^0$-penalization method for rare and weak signals}
\label{subsec:subsetselection}
When $\sigma = 0$, Model (\ref{model}) reduces to  the ``noiseless"  model  $Y = X \beta$.  In this model, 
\citet{UP2} (see also \citet{DonohoHuo})
reveals a fundamental phenomenon on  sparse representation.
Fix $(X, Y)$ and consider the equation $Y = X \beta$.  Since $p > n$,  the equation has infinitely many solutions. However,  a {\it very sparse solution}, if exists, is {\it unique} under mild conditions on the design $X$, with all other solutions being much denser.
In fact, if the sparsest solution $\beta_0$ has $k$ elements, then all other solutions of
the equation $Y=X\beta$ must have at least $(\rank(X)-k+1)$ nonzero elements, and
$\rank(X)=n$ when $X$ is in a ``general position''.

From a practical viewpoint, we frequently believe that this unique sparse solution is the truth
(i.e., Occam's razor). Therefore, the
problem of variable selection can be solved by some {\it global methods}
designed for finding the  sparsest  solution to the equation $Y = X \beta$.

Since the $L^0$-norm is  (arguably)  the most natural way to measure the sparsity of a vector,
the above idea suggests that
the $L^0$-penalization method  is a ``fundamentally correct" (but computationally intractable)  method for variable selection, provided that some mild conditions hold (noiseless,
Signal-to-Noise Ratio (SNR) is  high, signals are sufficiently sparse \citep{UP2, DonohoHuo}).

Motivated by this,
in the past two decades, a long list of computationally tractable algorithms have been proposed
that approximate the solution of the $L^0$-penalization method,
including the lasso, SCAD, MC+,  and many more
\citep{AIC,  CandesTao, LARS, FanLi,  BIC, Tibshirani,  Zhang10, ZhangT11a,  YuB, Zou}.

With that being said, we must  note that these methodologies   were built upon a framework with
four tightly woven core components:  ``signals are rare but strong",   ``the truth is also the sparsest solution to $Y = X \beta$",   ``probability of exact recovery is an appropriate loss function", and ``$L^0$-penalization method is  a  fundamentally correct approach".
Unfortunately, when signals are rare and weak, such a  framework is no longer suitable.
\begin{itemize}
\item  When signals are ``rare and weak", the fundamental uniqueness property of the sparse solution
in the noiseless case is no longer valid in the noisy case.
Consider the model $Y = X \beta + \sigma z$ and
suppose that a sparse $\beta_0$ is the true signal vector. There are many
vectors $\beta$ that are small perturbations of $\beta_0$
 such that the two models
$Y =  X \beta + \sigma z$ and $Y = X \beta_0 + \sigma z$ are
indistinguishable (i.e., all tests are asymptotically powerless).
In the ``rare and strong'' regime, $\beta_0$ is the sparsest solution among
all such ``eligible'' solutions of $Y = X \beta + \sigma z$.
However, this claim no longer holds in the ``rare and weak'' regime and
the principle of Occam's razor may not be as relevant as before.
\item The $L^0$-penalization method is originally designed for  ``rare and strong" signals where ``exact recovery" is used to measure its performance \citep{UP2, DonohoHuo, Donoho2006a}.
 When we must consider ``rare and weak" signals and when we use the
Hamming distance as the loss function,  it is unclear whether the $L^0$-penalization method
is still ``fundamentally correct".
\end{itemize}
In fact,   in  Section \ref{subsec:lasso} (see also \citet{JiJin}),  we show that the $L^0$-penalization method is not optimal in Hamming distance when signals are rare and weak, even with very simple designs (i.e., Gram matrix is tridiagonal or block-wise)  and even when the tuning parameter is ideally set.
Since the $L^0$-penalization method is used as the benchmark in the development of
many other penalization methods, its sub-optimality is expected to imply the sub-optimality
of other methods designed to match its performance
(e.g.,  lasso, SCAD, MC+).

\subsection{Limitation of Univariate Screening and UPS}
\label{subsec:UPS}
Univariate Screening (also called marginal regression or Sure Screening \citep{FanLv, GJW2011}) is a well-known variable selection method.  For $1 \leq j \leq p$,  recall that $x_j$ is the $j$-th column of $X$.  Univariate Screening selects variables with large marginal correlations:  $|(x_j, Y)|$, where $(\cdot, \cdot)$ denotes the inner product.  The method is computationally  fast, but it can be seriously corrupted by the so-called phenomenon of ``signal cancellation"  \citep{Wasserman}.
In our model (\ref{model})-(\ref{diag}),  the  SNR
associated with $(x_j, Y)$ is
\begin{equation} \label{signalcancel1}
\frac{1}{\sigma}  \sum_{\ell = 1}^p  (x_j, x_{\ell}) \beta_{\ell}  = \frac{\beta_j}{\sigma} +  \frac{1}
{\sigma}  \sum_{\ell \neq j} (x_j, x_{\ell}) \beta_{\ell}.
\end{equation}
``Signal cancellation" happens if  SNR is significantly smaller than $\beta_j/ \sigma$.
For this reason, the success of  Univariate Screening needs relatively strong conditions (e.g., Faithfulness Condition \citep{GJW2011}),  under which signal cancellation does not have a major effect.

In \citet{JiJin}, Ji and Jin proposed  {\it Univariate Penalized Screening (UPS)} as a refinement of Univariate Screening,
where it was showed to  be optimal in the rare  and weak paradigm, for the following two scenarios. The first scenario is where the nonzero effects of variables are all positively correlated:
$(x_j\beta_j)' (x_k\beta_k)  \geq 0$  for all $\{j,k\}$.
This guarantees the faithfulness of the univariate association test.
The second scenario is a Bernoulli model where the ``signal cancellation" only has
negligible effects over the Hamming distance of UPS.

With that being said, UPS attributes its success mostly to the cleaning stage; the
screening stage of UPS uses nothing but Univariate Screening, so UPS does not adequately address the challenge of ``signal cancellation".  For this reason,  we should not expect UPS  to be optimal in much more general settings.

\subsection{Limitations of Brute-force Multivariate Screening}
\label{subsec:BMSC}
One may attempt to overcome  ``signal cancellation" by
 multivariate screening, with
Brute-force Multivariate Screening (BMS) being the most straightforward version.
Fix an integer $1 \leq m_0 \ll  p$.   BMS
consists of a series of screening phases,   indexed by $m$, $1 \leq m \leq m_0$,  that are  increasingly more ambitious.
In Phase-$m$ BMS,
we test the significance of the association between $Y$ and any set of $m$ different
design variables $\{x_{j_1}, x_{j_2}, \ldots, x_{j_m}\}$, $j_1 < j_2 < \ldots < j_m$, and retain all such design variables if the test is significant.
The problem of BMS is,  it  enrolls too many
candidates for screening, which is both unnecessary and unwise.
\begin{itemize}
\item ({\it Screening inefficiency}). In Phase-$m$ of BMS, we test about $\binom{p}{m}$ hypotheses involving different subsets of $m$ design variables. The larger the number of  hypotheses we consider, the higher the   threshold we need to set for the tests, in order to control the false positives. When we enroll too many candidates for hypothesis testing, we need
signals that are stronger than necessary in order for them to survive the screening.
\item ({\it Computational challenge}).   Testing $\binom{p}{m}$ hypotheses is
computationally infeasible when $p$ is large, even when $m$ is very small (say,  $(p, m) = (10^4, 3)$).
\end{itemize}

\subsection{Graphlet Screening: how it is different and how it works}
\label{subsec:interplay}
Graphlet Screening (GS) uses a similar screening strategy as BMS does,
except for a major difference.  When it comes to the  test of significance
between $Y$ and design variables  $\{x_{j_1}, x_{j_2}, \ldots, x_{j_m}\}$, $j_1 < j_2 < \ldots < j_m$,  GS only carries out such a test if $\{j_1, j_2, \ldots, j_m\}$ is
a connected subgraph of the GOSD. Otherwise, the test is
safely skipped!

Fixing an appropriate threshold $\delta > 0$,  we let $\Omega^{*, \delta}$ be the regularized Gram matrix:
\begin{equation} \label{DefineOmega*}
\Omega^{*, \delta}(i,j) = G(i,j) 1\{|G(i,j)| \geq \delta\},   \qquad 1 \leq i, j \leq p.
\end{equation}
The GOSD $\cg \equiv \cg^{*, \delta} = (V, E)$ is the graph where $V = \{1, 2, \ldots, p\}$ and there is an edge between nodes $i$ and $j$  if and only if $\Omega^{*, \delta}(i,j)  \neq 0$.  See Section \ref{subsec:UB} for the choice of $\delta$.

{\bf Remark}.
GOSD and $\cg$ are generic terms which vary from case to case, depending on  $G$  and $\delta$. GOSD is very different from the Bayesian conditional independence graphs \citep{Pearl}.

Fixing $m_0 \geq 1$ as in BMS,   we  define
\begin{equation} \label{DefinecalA}
{\cal A}(m_0) =  {\cal A}(m_0;  G, \delta) =  \{ \mbox{all connected subgraphs of $\cg^{*, \delta}$ with size $\leq m_0$} \}.
\end{equation}
GS is a Screen and Clean method, consisting of   a graphical screening step (GS-step) and a graphical cleaning step (GC-step).
\begin{itemize}
\item {\it $GS$-step}. We test the significance of association between $Y$ and $\{x_{j_1}, x_{j_2}, \ldots, x_{j_m}\}$ if and only if $\{j_1, j_2, \ldots, j_m \} \in {\cal A}(m_0)$ (i.e., graph guided multivariate screening). Once   $\{j_1,\ldots,j_m\}$  is retained, it remains there until the end of the GS-step.
\item  {\it $GC$-step}. The set of surviving  nodes decompose into many small-size components, which we  fit separately using an  efficient low-dimensional test for small graphs.
\end{itemize}
GS is similar to  \cite{Wasserman} for  both of them have   a screening
and a cleaning stage, but is more sophisticated.
For clarification, note that Univariate Screening or BMS introduced earlier does not contain a cleaning stage and can be viewed as a counterpart of the GS-step.

We briefly explain why GS works. We discuss the GS-step and GC-step separately.

Consider the GS-step first.  Compared with BMS,  the GS-step  recruits far fewer candidates for screening, so it is able to overcome the two major shortcomings of BMS aforementioned:   high computational cost and low statistical efficiency. In fact,
fix $K \geq 1$ and suppose $\cg^{*, \delta}$ is $K$-sparse (see Section \ref{subsec:notations} for the definition).  By a well-known result in graph theory \citep{Frieze},
\begin{equation} \label{cardcalA}
|{\cal A}(m_0)|  \leq C m_0  p (e K)^{m_0}.
\end{equation}
The right hand side is much smaller than the term $\binom{p}{m_0}$ as we encounter in BMS.

At the same time, recall that $S = S(\beta)$ is the support of $\beta$. Let $\cg_S \equiv \cg_S^{*, \delta}$ be the subgraph of $\cg^{*, \delta}$ consisting all signal nodes. We can always split $\cg_S^{*, \delta}$  into ``graphlets" (arranged lexicographically) as follows:
\begin{equation} \label{graphlets}
\cg_{S}^{*, \delta}  =  \cg_{S,1}^{*, \delta} \cup \cg_{S,2}^{*, \delta} \ldots \cup \cg_{S,M}^{*, \delta},
\end{equation}
where each $\cg_{S,i}^{*,\delta}$ is a component (i.e., a maximal connected subgraph) of $\cg_S^{*, \delta}$, and different $\cg_{S,i}^{*, \delta}$ are not connected in $\cg_S^{*, \delta}$. Let
\begin{equation} \label{Definem0}
m_0^* = m_0^*(S(\beta), G, \delta) = \max_{1 \leq i \leq M} | \cg_{S, i}^{*, \delta}|
\end{equation}
be the maximum size of such graphlets (note that $M$ also depends on ($S(\beta)$, $G$, $\delta$)).

In many cases, $m_0^*$ is small. One such case is when we have a Bernoulli signal model.
\begin{lemma} \label{lemma:sparsegraph}
Fix $K \geq 1$ and $\eps > 0$. If $\cg^{*, \delta}$ is $K$-sparse and $\sgn(|\beta_1|), \sgn(|\beta_2|), \ldots,
\sgn(|\beta_p|)$ are $iid$   from $\mathrm{Bernoulli}(\eps)$, then except for a probability  $p (e \eps  K)^{m_0+1}$,
$m_0^*(S(\beta), G, \delta) \leq m_0$.
\end{lemma}
Lemma \ref{lemma:sparsegraph} is not tied to the Bernoulli model and holds more generally. For example, it holds when $\{\sgn(|\beta_i|)\}_{i = 1}^p$ are generated according to certain Ising models \citep{Ising}.

We recognize that in order for the GS-step to be efficient both in screening and in computation,  it is sufficient that
\begin{equation} \label{m0ineq}
m_0 \geq m_0^*.
\end{equation}
In fact, first,   if (\ref{m0ineq}) holds,  then for each $1 \leq \ell \leq M$,
$\cg_{S, \ell}^{*, \delta}  \in {\cal A}(m_0)$. Therefore,
at some point of the screening process of the GS-step,  we must have considered a significance test between $Y$ and the set of design variables $\{x_j:  j \in \cg_{S, \ell}^{*, \delta} \}$. Consequently,
the GS-step is able to overcome the ``signal cancellations" (the explanation is a little bit long, and we slightly defer it).
Second,  since $m_0^*$ is small in many situations,
we could choose a relatively small $m_0$ such that
 (\ref{m0ineq}) holds.   When $m_0$ is small,  as long as $K$ is small or moderately large,  the GS-step is computationally feasible. In fact,  the right hand side of (\ref{cardcalA}) is only larger than $p$ by a   moderate factor. See Section \ref{subsec:compu} for more discussion on the computation complexity.

We now explain the first point above.
The  notations below are frequently used.
\bed
For $X$ in Models (\ref{model})-(\ref{columnX}) and any subset $\call\subset \{1,2,...,p\}$,  let $P^{\call} = P^{\call}(X)$ be the projection from $\mathbb{R}^n$ to the subspace spanned
by $\{x_j: j \in {\cal I} \}$.
\eed
\bed
For an  $n  \times p$ matrix $A$ and  sets ${\cal I} \subset \{1,  \ldots, n\}$ and ${\cal J}  \subset \{1,   \ldots, p\}$, $A^{{\cal I}, {\cal J}}$
is the $|{\cal I}|  \times |{\cal J}|$ sub-matrix  formed by restricting the rows of $A$ to ${\cal I}$ and columns  to ${\cal J}$.
\eed
When $p = 1$, $A$ is a vector, and $A^{\call}$ is the sub-vector of $A$  formed by restricting the rows of $A$ to $\call$ .
When $\call = \{1, 2, \ldots, n\}$ (or ${\cal J} = \{1, 2, \ldots, p\}$),
we write $A^{\call, {\cal J}}$ as  $A^{\otimes, {\cal J}}$ (or $A^{\call, \otimes}$).
Note that indices in $\call$ or ${\cal J}$ are not necessarily sorted ascendingly.

Recall that for each $1 \leq \ell \leq M$, at some point of the GS-step,  we must have considered a significance test between $Y$ and the set of design variables $\{x_j:  j \in \cg_{S, \ell}^{*, \delta} \}$.   By (\ref{graphlets}), we rewrite Model (\ref{model}) as
\[
Y =  \sum_{\ell = 1}^M  X^{\otimes, \cg_{S, \ell}^{*, \delta}}  \beta^{\cg_{S, \ell}^{*, \delta}}  +  \sigma z, \qquad z \sim N(0, I_n).
\]
The key is  the set of  matrices $\{X^{\otimes, \cg_{S, \ell}^{*, \delta}}:  1 \leq \ell \leq M\}$ are  nearly orthogonal (i.e., for any column $\xi$ of $X^{\otimes, \cg_{S, k}^{*, \delta}}$ and any column $\eta$ of
$X^{\otimes, \cg_{S, \ell}^{*, \delta}}$,   $|(\xi, \eta)|$ is
 small  when $k \neq \ell$).

When we test the significance between $Y$ and $\{x_j, j \in \cg_{S, \ell}^{*, \delta} \}$, we are testing the null hypothesis $\beta^{\cg_{S, \ell}^{*, \delta}} = 0$ against the alternative
$\beta^{\cg_{S, \ell}^{*, \delta}} \neq 0$.  By the near orthogonality aforementioned, approximately,
$(X^{\otimes, \cg_{S, \ell}^{*, \delta}})'Y$ is a sufficient statistic for $\beta^{\cg_{S, \ell}^{*, \delta}}$, and the optimal test
 is  based on the $\chi^2$-test statistic  $\| P^{\cg_{S, \ell}^{*, \delta}} Y\|^2$.

The near orthogonality also implies that significant ``signal cancellation" only happens among signals within the same graphlet. When we screen each graphlet {\it as a whole}  using the $\chi^2$-statistic above,   ``signal cancellation" between different graphlets  only has negligible effects. In this way,
GS-step is able to retain all nodes in  $\cg_{S, \ell}^{*, \delta}$ in a nearly optimal way, and so overcome the challenge of ``signal cancellation". This explains the first point.

Note that the GS-step consists of a sequence  of sub-steps, each sub-step is associated with an element of ${\cal A}(m_0)$.   When we screen $\cg_{S, \ell}^{*, \delta}$ as a whole, it is possible some of the nodes have already been retained in the previous sub-steps.  In this case, we implement the $\chi^2$-test slightly differently, but the insight is similar. See Section \ref{subsec:graph} for details.

We now discuss the GC-step.    Let $\hat{S}$ be all the surviving nodes of the GS-step, and let $\cg_{\hat{S}}^{*, \delta}$ be the subgraph of $\cg^{*, \delta}$  formed by confining all nodes to $\hat{S}$. Similarly, we have (a) the decomposition $\cg_{\hat{S}}^{*, \delta}  =  \cg_{\hat{S}, 1}^{*, \delta}  \cup \cg_{\hat{S},2}^{*, \delta}  \ldots \cup  \cg_{\hat{S}, \hat{M}}^{*, \delta}$, (b) the near orthogonality between the $\hat{M}$ different matrices, each is
formed by $\{x_j: j \in \cg_{\hat{S}, \ell}^{*, \delta} \}$.
Moreover,  a carefully tuned screening stage of the GS ensures that most of the  components  $\cg_{\hat{S}, \ell}^{*, \delta}$ are only small perturbations of their counterparts in the decomposition  of $\cg_S^{*,\delta} = \cg_{S,1}^{*, \delta} \cup \cg_{S,2}^{*, \delta}  \ldots \cup  \cg_{S, M}^{*, \delta}$ as in (\ref{graphlets}),
and the maximum size of $\cg_{\hat{S}, \ell}^{*, \delta}$ is not too much larger than  $m_0^* = m_0^*(S(\beta), G, \delta)$.     Together, these allow us to clean  $\cg_{\hat{S}, \ell}^{*, \delta}$ separately, without much loss of efficiency.
Since the maximum size of $\cg_{\hat{S}, \ell}^{*, \delta}$ is small, the computational complexity in the cleaning stage  is moderate.

\subsection{Content}
\label{subsec:notations}
The remaining part of the paper is organized as follows.
In Section \ref{sec:oldintro},  we show that GS achieves the minimax Hamming distance in the Asymptotic Rare and Weak (ARW) model, and
 use the phase diagram   to  visualize the optimality of GS, and to illustrate the   advantage of GS   over the $L^0/L^1$-penalization methods.
In Section \ref{sec:main}, we explain that GS  attributes its optimality  to the so-called {\it Sure Screening} property and the {\it  Separable  After Screening}  property, and use these two properties to prove our main result,   Theorem \ref{thm:UB}.
Section \ref{sec:Simul} contains numeric
results,   Section \ref{sec:Discu} discusses  more connections to existing literature and possible extensions of GS,  and Section \ref{sec:appen} contains technical proofs.

Below are some notations we use in this paper.
$L_p$ denotes a generic multi-$\log(p)$ term that may vary from occurrence to occurrence; see Definition 5.
For a vector $\beta \in R^p$,   $\|\beta\|_q$  denotes  the $L^q$-norm, and when $q=2$, we drop $q$ for simplicity.
For two vectors $\alpha$ and $\beta$ in $R^p$,
$\alpha \circ \beta \in R^p$ denotes the vector in $R^p$ that satisfies  $(\alpha \circ \beta)_i  = \alpha_i  \beta_i$,  $1 \leq i \leq p$; ``$\circ"$ is known as the Hadamard product.

For an $n \times p$ matrix $A$, $\|A\|_\infty$ denotes the matrix $L^\infty$-norm, and $\|A\|$ denotes the spectral norm \citep{Horn1990}.
Recall that for two  sets ${\cal I}$ and ${\cal J}$ such that ${\cal I}  \subset \{1, 2, \ldots, n\}$ and
${\cal J} \subset \{1, 2, \ldots, p\}$,
$A^{{\cal I}, {\cal J}}$ denotes the submatrix of $A$ formed by restricting
the rows and columns of $A$ to ${\cal I}$ and ${\cal J}$, respectively. Note that
the indices in ${\cal I}$ and ${\cal J}$ are not necessarily sorted in the ascending order.
In the special case where ${\cal I} = \{1, 2, \ldots, n\}$ (or ${\cal J} = \{1, 2, \ldots, p\}$), we write $A^{\call, {\cal J}}$ as $A^{\otimes, {\cal J}}$ (or $A^{{\cal I}, \otimes}$).
In the special case where $n = p$ and $A$ is positive definite, $\lambda_k^*(A)$  denotes  the minimum eigenvalue of all the size $k$ principal submatrices of $A$, $1 \leq k \leq p$.
For $X$ in \eqref{model}, $P^{\call}$ denotes the projection to the column space of $X^{\otimes, \call}$.

Recall that in Model (\ref{model}),  $Y = X \beta + \sigma z$.
Fixing a threshold $\delta > 0$.  Let $G = X'X$ be the Gram matrix, and
let  $\Omega^{*, \delta}$ be the regularized Gram matrix defined by
 $\Omega^{*, \delta}(i,j)  = G(i,j) 1\{ |G(i,j)| \geq \delta\}$, $1 \leq  i, j \leq p$.
Let $\cg^{*, \delta}$ be
the graph where each index in $\{1, 2, \ldots, p\}$ is a node,
and there is an edge between node $i$ and node $j$ if
and only if $\Omega^{*, \delta}(i,j) \neq 0$.
We let $S(\beta)$ be the support of $\beta$, and denote $\cg_S^{*, \delta}$  by the subgraph of $\cg^{*, \delta}$ formed by all nodes in $S(\beta)$.
We call $\cg^{*, \delta}$ the Graph of Strong Dependence (GOSD) and sometimes write it by $\cg$ for short.
The GOSD and $\cg$ are generic notations which depend on $(G, \delta)$ and  may vary from occurrence to occurrence.
We also denote $\cg^{\diamond}$ by the Graph of Least Favorable (GOLF). GOLF only involves the study of the information lower bound.
For an integer $K \geq 1$,  a graph $\cg$,   and one of its subgraph $\call_0$, we write $\call_0\triangleleft \cg$ if and only if $\call_0$ is a component of $\cg$ (i.e., a maximal connected subgraph of $\cg$), and  we call  $\cg$ $K$-sparse if its maximum degree  is no greater than $K$.

\section{Main results}  \label{sec:oldintro}
In Section \ref{subsec:graph},
we formally introduce GS. In Section \ref{subsec:compu},
we discuss the computational complexity of GS.
In Sections \ref{subsec:randd}-\ref{subsec:UB}, we show that GS achieves the optimal rate of convergence in the Asymptotic Rare and Weak model.
In Sections \ref{subsec:simpl}-\ref{subsec:lasso}, we introduce the
notion of phase diagram and use it to compare GS
with the  $L^0/L^1$-penalization methods.   We conclude the section with a summary in Section \ref{subsec:summary}.


\subsection{Graphlet Screening: the procedure}
\label{subsec:graph}
GS consists of a GS-step and a GC-step.  We describe two steps separately.
 Consider the $GS$-step first.  Fix $m_0 \geq 1$ and $\delta > 0$, recall that $\cg^{*, \delta}$ denotes the GOSD and ${\cal A}(m_0)$
consists of all connected subgraphs of $\cg^{*, \delta}$ with
size $\leq m_0$.
\bi
\item {\it Initial sub-step}. Let ${\cal U}_p^* = \emptyset$.  List all elements in ${\cal A}(m_0)$  in the ascending order of the number of nodes it contains,   with ties broken lexicographically.  Since a node is thought of as connected to itself, the first $p$ connected subgraphs on the list are simply the nodes $1, 2, \ldots, p$. We screen all connected subgraphs in the order they are listed.
\item {\it Updating sub-step}.  Let $\call_0$ be the connected subgraph under consideration, and let ${\cal U}_p^*$  be the current set of retained indices.  We update ${\cal U}_p^*$ with a $\chi^2$ test as follows.
Let $\he =  \call_0 \cap {\cal U}_p^*$ and $\hd =  \call_0 \setminus {\cal U}_p^*$, so that $\he$ is the set of  nodes in $\call_0$ that have already been accepted, and
 $\hd$ is  the set of  nodes in $\call_0$ that is currently  under investigation.
 Note  that no action is needed if $\hd = \emptyset$.
For a threshold $t(\hd, \he) > 0$ to be determined, we
update ${\cal U}_p^*$ by adding all nodes in $\hd$ to it if
\begin{equation}\label{DefineDF}
T(Y, \hat{D}, \he)
= \|P^{\call_0}Y\|^2 - \|P^{\hat F}Y\|^2  > t(\hd, \he),
\end{equation}
and we keep ${\cal U}_p^*$ the same otherwise (by default, $\| P^{\he} Y \| = 0$ if $\he = \emptyset$).
We continue this process until we finish screening all connected subgraphs on the list. The final set of retained indices is denoted by ${\cal U}_p^*$.
\ei
See Table $1$ for a recap of the procedure.
In the $GS$-step,  once a node is kept in any sub-stage
of the screening process, it remains there until the end of the $GS$-step (however, it may be killed in the $GC$-step). This has a similar flavor to that of the Forward regression.

In principle, the procedure depends on how the connected subgraphs of the same size are initially ordered, and different ordering could give different numeric results. However, such differences are usually negligibly small.  Alternatively,  one could  revise the procedure so that it does not depend on the ordering.  For example, in the updating sub-step, we could choose to update ${\cal U}_p^*$ only when  we finish screening all connected sub-graphs of size $k$, $1 \leq k \leq m_0$. While the theoretic results below continue to hold if
we revise GS in this way, we must note that from a numeric
perspective, the revision would not produce a very different result.  For reasons of space, we skip discussions along this line.

The GS-step uses a  set of tuning parameters:
\[
\tq  \equiv  \{t(\hd, \he):  \mbox{$(\hd, \he)$  are as defined in (\ref{DefineDF})}   \}.
\]
A convenient way to set these parameters is to let $t(\hd,\he)=2\sigma^2q\log p$ for a fixed $q>0$ and all $(\hd, \he)$.  More sophisticated
choices are given in Section \ref{subsec:UB}.

The $GS$-step has two important properties: {\it Sure Screening} and {\it  Separable   After Screening (SAS)}. With tuning parameters $\tq$ properly set, the Sure Screening property says that
${\cal U}_p^*$ retains all but a negligible fraction of the signals.
 The SAS property says that as a subgraph of
 $\cg^{*,\delta}$, ${\cal U}_p^*$ decomposes into many disconnected components,  each has a size $\leq \ell_0$ for a fixed small integer $\ell_0$.  Together, these two properties enable us to reduce the original
large-scale regression problem to many small-size  regression problems that can be solved parallelly in the $GC$-step.
See Section \ref{sec:main} for elaboration on these ideas.

We now discuss the $GC$-step.   For any $1 \leq j \leq p$, we have either $j \notin {\cal U}_p^*$,   or  that there is a unique connected subgraph  $\call_0$ such that $j \in \call_0 \lhd {\cal U}_p^*$. In the first case, we estimate $\beta_j$ as $0$.  In the second case, for two tuning parameters $u^{gs} > 0$ and $v^{gs}  > 0$,  we estimate the whole set  of variables $\beta^{\call_0}$ by minimizing the functional
\begin{equation} \label{Pstep}
\|P^{\call_0}(Y - X^{\otimes, \call_0}  \xi)\|^2 + (u^{gs})^2\|\xi\|_0
\end{equation}
over all $|\call_0| \times 1$ vectors $\xi$, each nonzero coordinate of which
$\ge v^{gs}$ in magnitude.
The resultant estimator is the final estimate of GS which we denote by $\hb^{gs} = \hb^{gs}(Y; \delta, \tq,   u^{gs}, v^{gs},  X, p, n)$. See Section \ref{subsec:notations} for notations used in this paragraph.

\begin{table}
\caption{Graphlet Screening Algorithm.}
\hspace{.3in}
{\begin{tabular}{cl}\\ \hline
$GS$-step: & List ${\cal G}^{*, \delta}$-connected submodels $\call_{0, k}$ with $|\call_{0,1}|\le |\call_{0,2}|\le\cdots\le m_0$ \\
& Initialization: ${\cal U}_p^*=\emptyset$ and $k=1$ \\
& Test $H_0: \call_{0,k}\cap {\cal U}_p^*$ against $H_1: \call_{0,k}$ with $\chi^2$ test (\ref{DefineDF})\\
& Update: ${\cal U}_p^*\leftarrow {\cal U}_p^*\cup\call_{0,k}$ if $H_0$ rejected, $k\leftarrow k+1$
\\ \hline
$GC$-step: & As a subgraph of ${\cal G}^{*, \delta}$, ${\cal U}_p^*$ decomposes into many components $\call_0$ \\
& Use the $L^0$-penalized test (\ref{Pstep}) to select a subset ${\hat\call}_0$ of each $\call_0$\\
& Return the union of ${\hat\call}_0$ as the selected model
\\ \hline
\end{tabular}}
\end{table}

Sometimes for linear models with random designs,
the Gram matrix $G$  is very noisy, and GS is more effective
if we use it iteratively for a few times ($\leq 5$). This can be implemented in a similar way as that in   \citet[Section 3]{JiJin}.
Here, the main purpose of iteration is to denoise $G$,  not for
variable selection. See \citet[Section 3]{JiJin} and  Section \ref{sec:Simul}
for more discussion.

\subsection{Computational complexity} \label{subsec:compu}
If we exclude the overhead of obtaining $\cg^{*, \delta}$, then
the computation cost of GS contains  two  parts, that of the GS-step and that of the GC-step.
In each part, the computation cost hinges on the sparsity of $\cg^{*, \delta}$.
In Section \ref{subsec:randd}, we show that with a properly chosen $\delta$, for a wide class of design matrices, $\cg^{*, \delta}$ is $K$-sparse  for some
$K = K_p \leq  C \log^{\alpha}(p)$ as $p \goto \infty$, where $\alpha > 0$ is a constant.
As a result \citep{Frieze},
\begin{equation}   \label{fz}
|{\cal A}(m_0)|  \leq  p m_0 (e K_p)^{m_0} \leq C m_0   p  \log^{m_0 \alpha}(p).
\end{equation}
We now discuss two parts separately.

In the GS-step, the computation cost comes from that of listing all elements in ${\cal A}(m_0)$,
and that of screening all connected-subgraphs in ${\cal A}(m_0)$.
Fix $1 \leq k \leq m_0$.   By (\ref{fz}) and the fact that every size $k$ ($k>1$) connected subgraph  at least contains one size $k-1$ connected subgraph, greedy algorithm can be used to list all sub-graphs
with size $k$ with computational complexity $\leq C p(K_p  k)^k  \leq C p \log^{k \alpha}(p)$,
and screening all connected subgraphs of size $k$ has computational complexity $\leq C n p \log^{k \alpha}(p)$.  Therefore, the computational complexity of the GS-step  $ \leq C n p (\log(p))^{(m_0+1) \alpha}$.

The computation cost of the $GC$-step contains the part of breaking ${\cal U}_p^*$ into
disconnected components, and that of cleaning each component by minimizing \eqref{Pstep}.
As a well-known application of the breadth-first search \citep{Hopcroft}, the first part   $\leq|{\cal U}_p^*| (K_p+1)$.  For the second part,  by the SAS property of the $GS$-step (i.e.,  Lemma \ref{lemma:SAS}), for a broad class of
design matrices, with the tuning parameters chosen properly,    there is a {\it fixed}  integer $\ell_0$  such that with overwhelming probability,  $|\call_0| \leq \ell_0$ for any
$\call_0 \lhd {\cal U}_p^*$.
As a result,  the total computational  cost of the $GC$-step is no greater than
$C (2^{\ell_0}\log^{\alpha}(p))|{\cal U}_p^*|n$,  which is moderate.

The computational complexity of GS is only moderately larger than that of Univariate Screening or  UPS \citep{JiJin}.
UPS uses univariate thresholding for screening which has a computational complexity of
$O(np)$,  and GS implements
multivariate screening for all connected subgraphs in ${\cal A}(m_0)$, which has a computational complexity $\leq Cn p  (\log(p))^{(m_0+1) \alpha}$.  The latter is only larger by a multi-$\log(p)$ term.

\subsection{Asymptotic Rare and Weak model and  Random Design model}
\label{subsec:randd}
To analyze GS, we consider the regression model $Y = X \beta + \sigma z$ as in (\ref{model}),
and use an Asymptotic Rare and Weak (ARW) model for $\beta$ and a random design model for $X$.

We introduce the ARW first.
Fix parameters $\eps \in (0,1)$, $\tau > 0$, and $a \geq  1$.
Let  $b = (b_1, \ldots, b_p)'$ be the $p \times 1$ random vector where
\begin{equation} \label{RW1}
b_i \stackrel{iid}{\sim}    \mathrm{Bernoulli}(\eps).
\end{equation}
We model  the signal vector $\beta$ in Model (\ref{model}) by
\begin{equation} \label{RW2}
\beta =  b \circ \mu,
\end{equation}
where ``$\circ$" denotes the Hadamard product (see Section \ref{subsec:notations}) and $\mu \in \Theta_p^*(\tau, a)$, with
\begin{equation} \label{RW3}
\Theta_p^*(\tau,a) = \{\mu \in \Theta_p(\tau), \|\mu\|_{\infty} \leq a \tau\}, \qquad
\Theta_p(\tau)  =  \{\mu \in \R^p: \,  |\mu_i|\geq \tau, 1 \leq i \leq p \}.
\end{equation}
In this model,  $\eps$ calibrates the sparsity level and $\tau$ calibrates the
minimum signal strength. We are primarily interested in the case where $\eps$ is small and $\tau$ is  smaller than the required signal strength for the exact recovery of the support of $\beta$, so the signals are both rare and weak. The constraint of $\|\mu\|_{\infty} \leq a \tau_p$ is mainly for  technical reasons (only needed for  Lemma \ref{lemma:SAS}); see Section \ref{subsec:UB} for more discussions.

We let $p$ be the driving asymptotic parameter, and tie $(\eps, \tau)$ to $p$ through some fixed parameters.
In detail, fixing $0 <  \vartheta < 1$, we model
\begin{equation} \label{Defineeps}
\eps = \eps_p =  p^{-\vartheta}.
\end{equation}
For any fixed $\vartheta$,  the signals become increasingly sparser as $p \goto \infty$. Also, as $\vartheta$ ranges, the sparsity level ranges from very dense to very sparse, and covers all interesting cases.

It turns out that the most interesting range for $\tau$ is $\tau = \tau_p  = O(\sqrt{\log(p)})$.  In fact,
when  $\tau_p  \ll \sigma\sqrt{\log(p)}$,  the signals are simply too rare and weak so that successful  variable selection is impossible.   On the other hand,  exact support recovery requires $\tau\gtrsim  \sigma\sqrt{2\log p}$ for orthogonal designs and possibly even larger $\tau$ for correlated designs.
In light of this, we fix $r > 0$ and calibrate $\tau$ by
\begin{equation} \label{Definetau}
\tau = \tau_p =
\sigma\sqrt{2 r \log(p)}.
\end{equation}

Next, consider the random design model.   The use of random design model is mainly for simplicity in presentation. The main results in the paper
can be translated to fixed design models with a careful modification of the notations; see Corollary \ref{cor:LB} and Section \ref{sec:Discu}.

For any positive definite matrix $A$, let $\lambda(A)$ be the smallest eigenvalue, and let
\begin{equation} \label{Definelambda*}
\lambda_k^*(\Omega) = \min\{\lambda(A):  \mbox{$A$ is a $k \times k$ principle submatrix of  $\Omega$} \}.
\end{equation}
For $m_0$ as in the GS-step,  let  $g = g(m_0, \vartheta, r)$ be the smallest integer such that
\begin{equation} \label{Defineg}
g \geq  \max\{m_0,  (\vartheta + r)^2/(2 \vartheta r) \}.
\end{equation}
Fixing a constant $c_0 > 0$,  introduce
\begin{equation}\label{def:mpc0g}
{\cal M}_p(c_0, g) = \{ \Omega:  \mbox{$p \times p$ correlation matrix, $\lambda_g^*(\Omega) \geq c_0$}   \}.
\end{equation}
Recall $X_i$ is the $i$-th row of $X$; see (\ref{columnX}). In the random design model, we fix an $\Omega  \in {\cal M}(c_0, g)$ ($\Omega$ is unknown to us), and assume
\begin{equation} \label{RD1}
X_i \stackrel{iid}{\sim} N(0,  \frac{1}{n} \Omega), \qquad 1 \leq i \leq n.
\end{equation}
In the literature, this is called the Gaussian design, which can be found in  Compressive Sensing \citep{Nowak2007},   Computer Security \citep{Nissim}, and other application areas.

At the same time, fixing   $\kappa \in (0,1)$,  we model the sample size $n$ by
\begin{equation} \label{RD2}
n = n_p = p^{\kappa}.
\end{equation}
As $p \goto \infty$, $n_p$ becomes increasingly large but is still much smaller than $p$.
We assume
\begin{equation} \label{RD3}
\kappa  >  (1 - \vartheta),
\end{equation}
so that $n_p \gg p \eps_p$. Note $p \eps_p$ is approximately the total number of signals.
Condition (\ref{RD3}) is almost necessary for successful variable
selection \citep{Donoho2006a, Donoho2006b}.
\bed\label{def:rd}
We call model (\ref{RW1})-(\ref{Definetau}) for $\beta$ the Asymptotic
Rare Weak model $\mathrm{ARW}(\vartheta, r, a,  \mu)$, and call Model (\ref{RD1})-(\ref{RD3}) for $X$ the Random Design model $\mathrm{RD}(\vartheta, \kappa, \Omega)$.
\eed

\subsection{Minimax Hamming distance}\label{subsec:hamm}
In many works on variables selection, one assesses the optimality by
the  `oracle property', where the probability of non-exact recovery  $P(\sgn(\hat{\beta}) \neq  \sgn(\beta))$ is the loss function.
When signals are rare and weak,  $P(\sgn(\hat{\beta}) \neq  \sgn(\beta)) \approx 1$ and `exact recovery' is usually impossible.
A more appropriate loss function is the Hamming distance between $\sgn(\hat{\beta})$ and $\sgn(\beta)$.

For any fixed $\beta$ and any variable selection procedure $\hb$,  we measure the performance by the Hamming distance:
\[
h_p(\hb, \beta \bigl| X) =    E \Bigl[  \sum_{j = 1}^p  1 \bigl\{  \sgn(\hb_j) \neq \sgn(\beta_j)  \bigr\}    \bigr| X  \Bigr].
\]
In the Asymptotic Rare Weak  model,  $\beta = b \circ \mu$, and $(\eps_p, \tau_p)$ depend on $p$ through $(\vartheta, r)$, so the overall Hamming distance for $\hb$ is
\[
H_p(\hb; \eps_p, n_p,   \mu,   \Omega)  = E_{\eps_p} E_{\Omega} \bigl[   h_p(\hb,  \beta \bigr| X)     \bigr] \equiv   E_{\eps_p}  E_{\Omega}  \bigl[h_p(\hb,  b \circ \mu \bigr| X ) \bigr],
\]
where $E_{\eps_p}$ is the expectation with respect to the law of $b$, and $E_{\Omega}$ is the expectation with respect to the law of $X$; see (\ref{RW1}) and (\ref{RD1}).
Finally, the minimax Hamming distance is
\begin{equation} \label{hammadd}
\hamm_p^*(\vartheta,  \kappa, r, a,  \Omega)  =    \inf_{\hb}  \sup_{\mu \in \Theta_p^*(\tau_p, a)}  \bigl\{   H_p(\hb; \eps_p, n_p,  \mu, \Omega)  \bigr\}.
\end{equation}
The Hamming distance
is no smaller than the sum of the expected number of signal components that are misclassified as noise and
the expected number of noise components that are misclassified as signal.

\subsection{Lower bound for the minimax Hamming distance, and GOLF} \label{subsec:LB}
We  first construct lower bounds for ``local risk" at different $j$,  $1 \leq j \leq p$, and
then aggregate them to construct a lower bound for the global risk. One challenge we face is  the least favorable configurations for different $j$ overlap with each other. We resolve this by exploiting the sparsity of a new graph to be introduced:
Graph of Least Favorable (GOLF).

To recap,  the model we consider is Model (\ref{model}),  where
\begin{equation}\label{eq:modelcondition}
  \beta \text{ is modeled by }  \mathrm{ARW}(\vartheta, r, a,  \mu), \qquad   \mbox{and} \qquad X \text{ is modeled by } \mathrm{RD}(\vartheta, \kappa, \Omega).
\end{equation}
Fix $1 \leq j \leq p$.  The ``local risk" at an index $j$ is the risk of estimating the set of variables
$\{\beta_k:  d(k, j) \leq g\}$, where $g$ is defined in (\ref{Defineg}) and  $d(j,k)$ denotes the geodesic distance between $j$ and $k$ in the graph $\cg^{*,\delta}$.
The  goal is to construct two subsets $V_0$ and $V_1$   and two realizations of $\beta$, $\beta^{(0)}$ and $\beta^{(1)}$ (in the special case of $V_0 = V_1$, we require $\sgn(\beta^{(0)}) \neq \sgn(\beta^{(1)})$),   such that
$j \in V_0 \cup V_1$ and
\[
\mbox{If $k \notin V_0 \cup V_1$,   $\beta^{(0)}_k = \beta_k^{(1)}$; otherwise,
$\beta^{(i)}_k \neq 0$ if and only if $k \in V_i$, $i = 0, 1$}.
\]
In the literature, it is known that how well we can estimate $\{\beta_k: d(k,j) \leq g\}$  depends on how well we can separate two hypotheses (where $\beta^{(0)}$ and $\beta^{(1)}$ are assumed as known):
\begin{equation} \label{testingprob}
H_0^{(j)}:  Y = X \beta^{(0)}   + \sigma z \qquad \mbox{vs}.  \qquad H_1^{(j)}:   Y=  X \beta^{(1)} + \sigma z, \qquad z \sim N(0, I_n).
\end{equation}
The least favorable configuration for the local risk at  index $j$ is the quadruple $(V_0, V_1, \beta^{(0)}, \beta^{(1)})$ for which two hypotheses are the
most difficult to separate.

For any $V \subset \{1, 2, \ldots, p\}$,  let $I_V$ be the
indicator vector of $V$ such that for any $1 \leq k \leq p$,  the $k$-th coordinate of $I_V$ is $1$  if $k \in V$ and is $0$ otherwise.  Define
\begin{equation}\label{defineBV}
B_{V}  =  \{ I_{V}  \circ \mu: \;  \mu \in  \Theta_p^*(\tau_p, a) \},
\end{equation}
where we recall ``$\circ$" denotes the Hadamard product (see Section \ref{subsec:notations}).
Denote for short $\theta^{(i)} = I_{V_0 \cup V_1} \circ \beta^{(i)}$, and so
$\beta^{(1)} - \beta^{(0)} = \theta^{(1)} - \theta^{(0)}$  and   $\theta^{(i)} \in B_{V_i}$, $i = 0, 1$. Introduce
\[
\alpha(\theta^{(0)}, \theta^{(1)})  =  \alpha(\theta^{(0)}, \theta^{(1)};  V_0, V_1,  \Omega, a)  =  \tau_p^{-2} (\theta^{(0)} - \theta^{(1)})' \Omega (\theta^{(0)} - \theta^{(1)}).
\]
For the testing problem in (\ref{testingprob}),
the optimal test  is to  reject $H_0^{(j)}$  if and only if $(\theta^{(1)} - \theta^{(0)})' X'  (Y - X \beta^{(0)})  \geq t \sigma \tau_p \sqrt{\alpha(\theta^{(0)}, \theta^{(1)})}$ for some threshold $t > 0$ to be determined.
In the ARW and RD models,  $P(\beta_k \neq 0, \forall k \in V_i) \sim \eps_p^{|V_i|}$, $i = 0, 1$, and $(\theta^{(0)} - \theta^{(1)})' G (\theta^{(0)} - \theta^{(1)})  \approx (\theta^{(0)} - \theta^{(1)})' \Omega (\theta^{(0)} - \theta^{(1)})$,
since the support of $\theta^{(0)} - \theta^{(1)}$ is contained in  a small-size set $V_0 \cup V_1$.
Therefore the sum of Type I and Type II error of any test associated with (\ref{testingprob})  is  no  smaller than (up to some  negligible differences)
\begin{equation} \label{testingrisk}
\eps_p^{|V_0|} \bphi(t) +  \eps_p^{|V_1|} \Phi \bigl(t - (\tau_p/\sigma) [\alpha(\theta^{(0)}, \theta^{(1)})]^{1/2}\bigr),
\end{equation}
where  $\bphi = 1 - \Phi$ is the survival function of $N(0,1)$.

For a  lower bound for the ``local risk" at $j$,  we first
optimize the quantity in (\ref{testingrisk}) over all $\theta^{(0)} \in B_{V_0}$ and
$\theta^{(1)} \in B_{V_1}$,  and then optimize over all $(V_0, V_1)$ subject to
$j \in V_0 \cup V_1$.  To this end,  define $\alpha^*(V_0, V_1) = \alpha^*(V_0, V_1; a,
\Omega)$, $\eta(V_0, V_1) = \eta(V_0, V_1; \vartheta, r, a,\Omega)$, and
$\rho_j^* = \rho_j^*(\vartheta, r, a, \Omega)$ by
\begin{equation} \label{optimalv0v1}
\alpha^*(V_0, V_1) = \min \bigl\{ \alpha(\theta^{(0)}, \theta^{(1)}; V_0, V_1, \Omega, a): \;  \theta^{(i)} \in B_{V_i} ,  i = 0, 1,     \sgn(\theta^{(0)}) \neq \sgn(\theta^{(1)})\},
\end{equation}
\[
\eta(V_0, V_1)  =  \max\{|V_0|, |V_1|\}\vartheta + \frac{1}{4}\left[ \left( \sqrt{\alpha^*(V_0,V_1) r} - \frac{ \bigl|  (|V_1|-|V_0|) \bigr|\vartheta}{\sqrt{\alpha^*(V_0,V_1) r}} \right)_+  \right]^2,
\]
and
\begin{equation}\label{eq:define-rhoj}
\rho_j^*(\vartheta, r, a, \Omega)   =  \min_{\{ (V_0, V_1):  j\in V_1\cup V_0 \}}   \eta(V_0, V_1).
\end{equation}
The  following shorthand notation is frequently used in this paper, which stands for a generic  multi-$\log(p)$ term that may vary from one occurrence to another.
\bed\label{def:Lp}
$L_p > 0$ denotes a multi-$\log(p)$ term such that when $p \goto \infty$, for any $\delta > 0$, $L_p p^{\delta} \goto \infty$ and $L_p p^{-\delta} \goto 0$.
\eed
By   (\ref{testingrisk}) and Mills' ratio \citep{Wasserman},
a lower bound for the ``local risk" at $j$ is
\begin{align}
&\qquad \qquad   \sup_{\{ (V_0, V_1): \;  j \in V_0 \cup V_1 \}}  \bigl\{  \inf_t  \bigl[ \eps_p^{|V_0|} \bphi(t) +  \eps_p^{|V_1|} \Phi \bigl(t - (\tau_p/\sigma) [\alpha^*(V_0, V_1)]^{1/2}  \bigr) \bigr] \bigr\}   \label{introLB}  \\
=  &  \sup_{\{ (V_0, V_1): \;  j \in V_0 \cup V_1 \}}  \big\{ L_p   \exp \bigl(- \eta(V_0, V_1) \cdot \log(p) \bigr) \bigr\}   =  L_p \exp( - \rho_j^*(\vartheta, r, a, \Omega)  \log(p)).
\end{align}

We now aggregate such lower bounds for ``local risk" for a global lower bound.
Since the ``least favorable" configurations of $(V_0, V_1)$ for different $j$ may overlap with each other, we need to consider a graph as follows. Revisit the optimization problem in (\ref{optimalv0v1}) and let
\begin{equation}\label{eq:defineV0V1star}
(V_{0j}^*, V_{1j}^*) = \margmin_{\{ (V_0, V_1):   j \in V_1 \cup V_0 \}}  \eta(V_0, V_1; \vartheta, r, a, \Omega).
\end{equation}
When there is a tie, pick the pair that appears first lexicographically.  Therefore,  for any  $1 \leq  j \leq p$,    $V_{0j}^* \cup V_{1j}^*$ is  uniquely defined.  In Lemma \ref{lemma:V} of the appendix, we show that $|V_{0j}^* \cup V_{1j}^*| \leq (\vartheta + r)^2/(2 \vartheta r)$ for all $1 \leq j \leq p$.

We now define a new graph,   Graph of Least Favorable (GOLF),
$\cg^{\diamond} = (V, E)$, where $V = \{1, 2, \ldots, p \}$ and there is an edge between $j$ and $k$  if and only if $(V_{0j}^* \cup V_{1j}^*)$ and $(V_{0k}^* \cup V_{1k}^*)$ have non-empty intersections. Denote the maximum degree of GOLF by $d_p(\cg^{\diamond})$.
\begin{thm} \label{thm:LB}
Fix $(\vartheta, \kappa) \in (0,1)^2$, $r > 0$, and  $a \geq  1$ such that $\kappa > (1 - \vartheta)$, and let ${\cal M}_p(c_0, g)$ be as in (\ref{def:mpc0g}). Consider Model (\ref{model}) where  $\beta$ is modeled by $ARW(\vartheta, r, a, \mu)$ and $X$ is modeled by $RD(\vartheta, \kappa, \Omega)$ and  $\Omega \in {\cal M}_p(c_0, g)$  for sufficiently large $p$.   Then as $p \goto \infty$,
$\hamm_p^*(\vartheta, \kappa,  r, a,  \Omega) \geq  L_p \bigl[d_p(\cg^{\diamond}) \bigr]^{-1}  \sum_{j = 1}^p p^{-\rho_j^*(\vartheta, r, a, \Omega)}$.
\end{thm}
A similar claim holds for deterministic design models;  the proof is similar so we omit it.
\begin{cor}\label{cor:LB}
For deterministic design models,
the parallel lower bound holds for the minimax Hamming distance with $\Omega$ replaced by $G$ in the calculation of $\rho_j^*(\vartheta, r, a, \Omega)$ and $d_p(\cg^{\diamond})$.
\end{cor}

{\bf Remark}.
The lower bounds contain a factor of $\bigl[d_p(\cg^{\diamond}) \bigr]^{-1}$.   In many cases including
that considered in our main theorem (Theorem
\ref{thm:UB}),  this factor is a multi-$\log(p)$ term so it does not have a major effect.
In some other cases, the factor  $\bigl[d_p(\cg^{\diamond}) \bigr]^{-1}$  could be much smaller, say, when
 the GOSD has one or a few hubs, the degrees of which grow algebraically fast as
 $p$ grows.  In these cases, the associated GOLF may (or may not) have large-degree hubs. As a result,  the lower bounds we derive could be  very  conservative, and can be substantially  improved if we  treat the hubs, neighboring nodes of the hubs, and other nodes separately. For  the sake  of space, we leave such discussion to   future work.

{\bf Remark}.  A similar lower bound holds if the condition $\mu\in \Theta^*_p(\tau_p,a)$ of ARW is replaced by $\mu\in \Theta_p(\tau_p)$.  In (\ref{optimalv0v1}), suppose we replace $\Theta_p^*(\tau_p, a)$ by $\Theta_p(\tau_p)$, and the minimum is achieved at $(\theta^{(0)}, \theta^{(1)}) =   (\theta_*^{(0)}(V_0, V_1; \Omega),  \theta_*^{(1)}(V_0, V_1; \Omega))$.
Let $g = g(m_0, \vartheta, r)$ be as in (\ref{Defineg}) and define
\[
a_g^*(\Omega) =   \max_{\{(V_0, V_1):  |V_0 \cup V_1| \leq g\} } \{ \| \theta_*^{(0)}(V_0, V_1; \Omega) \|_{\infty},  \|\theta_*^{(1)}(V_0, V_1; \Omega)\|_{\infty} \}.
\]
By elementary calculus, it is seen that for $\Omega \in {\cal M}_p(c_0, g)$,  there is a a constant $C = C(c_0, g)$ such that $a_g^*(\Omega) \leq C$.
If additionally we assume
\begin{equation} \label{Conditiona}
a > a_g^*(\Omega),
\end{equation}
then $\alpha^*(V_0, V_1) = \alpha^*(V_0, V_1;  \Omega, a)$,
$\eta(V_0, V_1; \Omega, a, \vartheta, r)$, and $\rho_j^*(\vartheta, r, a, \Omega)$ do not depend on $a$. Especially, we can derive an alternative formula for  $\rho_j^*(\vartheta, r, a, \Omega)$; see Lemma \ref{lemma:omega} for details.

When (\ref{Conditiona}) holds, $\Theta_p^*(\tau_p, a)$ is broad enough in the sense that  the least favorable configurations $(V_0, V_1, \beta^{(0)}, \beta^{(1)})$ for all $j$ satisfy $\|\beta^{(i)}\|_{\infty} \leq a \tau_p$, $i = 0, 1$.  Consequently, neither the minimax rate nor GS
needs to adapt to $a$.  In Section \ref{subsec:UB}, we assume (\ref{Conditiona}) holds;  (\ref{Conditiona}) is  a mild condition for it only involves small-size sub-matrices of $\Omega$.



\subsection{Upper bound and optimality of Graphlet Screening} \label{subsec:UB}
Fix constants $\gamma \in (0,1)$ and $A > 0$. Let ${\cal M}_p(c_0, g)$ be as in (\ref{def:mpc0g}).  In this section, we further restrict $\Omega$ to the following set:
\begin{equation} \label{DefinecalM}
{\cal M}_p^*(\gamma, c_0, g,  A) =  \Bigl\{\Omega \in {\cal M}_p(c_0, g): \,   \sum_{j =1}^p  |\Omega(i,j)|^{\gamma} \leq A,   \,    1 \leq i \leq p  \Bigr\}.
\end{equation}
Note that any $\Omega \in  {\cal M}_p^*(\gamma, c_0, g,  A)$  is sparse in the sense that each row of $\Omega$  has relatively few large coordinates.
The sparsity of $\Omega$ implies the sparsity of the Gram matrix $G$, since
small-size sub-matrices of $G$ approximately equal to their counterparts of $\Omega$.

In GS,
when we regularize GOSD (see (\ref{DefineOmega*})), we set the threshold $\delta$  by
\begin{equation} \label{choosedelta}
\delta = \delta_p = 1/\log(p).
\end{equation}
Such a choice for threshold is mainly for convenience, and can be replaced by any term that tends to $0$ logarithmically fast as $p \goto \infty$.

For any subsets $D$ and $F$ of $\{1, 2, \ldots, p\}$, define $\omega(D, F; \Omega) = \omega(D, F; \vartheta, r, a, \Omega, p)$ by
\begin{equation} \label{Defineomega}
 \omega(D, F; \Omega) = \min_{\xi  \in\R^{|D|}, \min_{i\in D}|\xi_i| \geq 1}
\Big\{ \xi' \bigl(\Omega^{D, D} - \Omega^{D, F} (\Omega^{F, F})^{-1} \Omega^{F, D} \bigr)  \xi \Big\},
\end{equation}
We choose the tuning parameters in the
$GS$-step in a way such that
\begin{equation} \label{Definet}
t(\hd, \he) = 2\sigma^2 q(\hd, \he)\log p,
\end{equation}
where $q=q(\hd, \he)>0$ satisfies (for short, $\omega=\omega({\hat D},{\hat F}; \Omega)$)
\begin{equation} \label{Defineq}
\left\{
\begin{array}{ll}
\sqrt{q_0} \leq \sqrt{q}  \leq \sqrt{\omega  r}  - \sqrt{\frac{(\vartheta + \omega  r)^2}{4 \omega r} -   \frac{|\hat{D}| + 1}{2} \vartheta},
&\   \mbox{$|\hat{D}|$ is odd \&  $\omega  r / \vartheta > |\hat{D}| + (|\hat{D}|^2 - 1)^{1/2}$},      \\
\sqrt{q_0} \leq \sqrt{q}  \leq \sqrt{\omega  r}  - \sqrt{\frac{1}{4} \omega  r - \frac{1}{2} |\hat{D}| \vartheta},   &\  \mbox{$|\hat{D}|$ is even \&  $\omega  r / \vartheta \geq 2 |\hat{D}|$},    \\
\mbox{$q$ is a constant such that $q \geq q_0$},  &\ \mbox{otherwise}.
\end{array}
\right.
\end{equation}
We set  the $GC$-step tuning parameters by
\begin{equation} \label{Tuneuv}
u^{gs} = \sigma \sqrt{2\vartheta \log p}, \qquad v^{gs}  = \tau_p = \sigma\sqrt{ 2 r \log p}.
\end{equation}
The main  theorem of  this paper is the following theorem.
\begin{thm} \label{thm:UB}
Fix $m_0 \geq 1$,   $(\vartheta, \gamma, \kappa) \in(0,1)^3$, $r > 0$,  $c_0>0$, $g > 0$, $a > 1$,   $A > 0$ such that $\kappa>1-\vartheta$ and   (\ref{Defineg})   is satisfied.
Consider Model (\ref{model}) where $\beta$
is modeled by $ARW(\vartheta, r, a, \mu)$, $X$ is modeled by $RD(\vartheta, \kappa, \Omega)$, and where  $\Omega \in {\cal M}_p^*(\gamma, c_0,g,  A)$ and $a > a_g^*(\Omega)$ for sufficiently large $p$.  Let $\hb^{gs}  = \hb^{gs}(Y; \delta, {\cal Q}, u^{gs}, v^{gs}, X, p, n)$ be the Graphlet Screening procedure defined as in Section \ref{subsec:graph}, where the tuning parameters $(\delta, {\cal Q}, u^{gs}, v^{gs})$ are set as in (\ref{choosedelta})-(\ref{Tuneuv}).
 Then  as $p \goto \infty$,
$\sup_{\mu \in \Theta_p^*(\tau_p, a)  } H_p(\hb^{gs}; \eps_p, n_p, \mu, \Omega)
\leq  L_p \Big[p^{1-(m_0 + 1) \vartheta} + \sum_{j = 1}^p p^{-\rho_j^*(\vartheta, r, a,  \Omega)}\Big] + o(1)$.
\end{thm}
Note that $\rho_j^* = \rho_j^*(\vartheta, r, a, \Omega)$ does not depend on $a$.
Also, note that   in the most interesting range,   $\sum_{j = 1}^p p^{-\rho_j^*} \gg 1$. So if we choose $m_0$ properly large (e.g.,  $(m_0 + 1)\vartheta > 1$), then
$\sup_{\mu \in \Theta_p^*(\tau_p, a)  } H_p(\hb^{gs}; \eps_p, n_p, \mu, \Omega) \leq  L_p   \sum_{j = 1}^p p^{-\rho_j^*(\vartheta, r, a,  \Omega)}$.
Together with  Theorem  \ref{thm:LB}, this says that GS  achieves the optimal rate of convergence, adaptively to all $\Omega$ in ${\cal M}_p^*(\gamma, c_0,g, A)$ and $\beta \in \Theta_p^*(\tau_p, a)$.
We call this property   {\it optimal adaptivity}. Note that  since the diagonals of $\Omega$ are  scaled to $1$ approximately,   $\kappa  \equiv  \log(n_p)/\log(p)$ does not have a major influence over the convergence rate,  as long as (\ref{RD3}) holds.

{\bf Remark}. Theorem \ref{thm:UB} addresses the case where (\ref{Conditiona}) holds so  $a > a_g^*(\Omega)$.
We now briefly  discuss the case  where  $a < a_g^*(\Omega)$.    In this case, the set
$\Theta_p^*(\tau_p, a)$ becomes sufficiently narrow and $a$ starts to have some influence over the optimal rate of convergence, at least for some choices of $(\vartheta, r)$.  To reflect the role of $a$,
we  modify GS as follows:
(a) in the $GC$-step (\ref{Pstep}), limit $\xi$ to the class where either $\xi_i = 0$ or $\tau_p \leq |\xi_i| \leq a \tau_p$, and (b) in the $GS$-step, replacing the $\chi^2$-screening by the likelihood based screening procedure; that is, when we screen $\call_0  = \hd \cup \he$,   we accept nodes in $\hd$ only
when $h(\he) >  h(\call_0)$, where for any subset $D \subset \{1, 2, \ldots, p\}$,
$h(D) =  \min  \bigl\{   \frac{1}{2}
\|P^D ( Y - X^{\otimes, D} \xi) \|^2 + \vartheta \sigma^2 \log(p) |D| \big\}$,
where the minimum is computed
over all $|D|\times 1$ vectors $\xi$ whose nonzero elements all have magnitudes in $[\tau_p,a\tau_p]$.
From a practical point of view, this modified procedure depends more on the underlying parameters and is harder to implement than is GS.  However, this is the price we need to pay when $a$ is small.  Since we are primarily interested in the case of relatively larger $a$ (so that $a > a_g^*(\Omega)$ holds),  we skip further discussion along this line.

\subsection{Phase diagram and examples where $\rho_j^*(\vartheta, r, a,  \Omega)$  have  simple  forms}
\label{subsec:simpl}
In general,  the exponents $\rho_j^*(\vartheta, r, a, \Omega)$ may depend on
$\Omega$ in a complicated way. Still, from time to time,  one may want to find a simple expression for $\rho_j^*(\vartheta, r, a, \Omega)$.   It turns out that in a wide class of situations, simple forms for $\rho_j^*(\vartheta, r, a, \Omega)$ are possible. The surprise is that, in many examples,  $\rho_j^*(\vartheta, r, a, \Omega)$  depends more on the
trade-off between the parameters $\vartheta$ and $r$ (calibrating the signal sparsity and signal strength, respectively), rather than on  the large coordinates of $\Omega$.

We begin with the following theorem, which is proved in \citet[Theorem 1.1]{JiJin}.
\begin{thm} \label{thm:universalLB}
Fix $(\vartheta, \kappa) \in (0, 1)$, $r > 0$, and $a > 1$ such that $\kappa > (1 - \vartheta)$.
Consider Model (\ref{model}) where $\beta$ is modeled by $ARW(\vartheta, r, a, \mu)$ and
$X$ is modeled by $RD(\vartheta, \kappa, \Omega)$.  Then as $p \goto \infty$,
\[
\frac{\hamm_p^*(\vartheta, \kappa, r, a, \Omega)}{p^{1 - \vartheta}}
\gtrsim
\left\{
\begin{array}{ll}
1,  &\qquad   0 < r <  \vartheta, \\
L_p p^{-(r - \vartheta)^2/(4r)}, &\qquad r > \vartheta.
\end{array}
\right.
\]
\end{thm}
Note that $p^{1 - \vartheta}$ is approximately the number of signals.
Therefore,  when $r < \vartheta$, the number of selection errors can not get substantially smaller than the number of signals.  This is the
most difficult case where no variable selection  method can be successful.

In this section,  we focus on the case $r > \vartheta$, so that
successful variable selection is   possible.  In this case,  Theorem
\ref{thm:universalLB} says that a {\it universal} lower bound for the Hamming distance  is
$L_p p^{1 - (\vartheta + r)^2/(4r)}$.
An interesting question is, to what extend, this lower bound is tight.

Recall  that $\lambda_k^*(\Omega)$ denotes the minimum of smallest eigenvalues across all $k \times k$ principle submatrices of $\Omega$, as defined in (\ref{Definelambda*}).
The following corollaries are proved in Section \ref{sec:appen}.
\begin{cor}  \label{cor:entrywise1}
Suppose the conditions of Theorem \ref{thm:UB} hold, and that additionally,
$1 <  r /  \vartheta <  3 + 2 \sqrt{2} \approx 5.828$, and $|\Omega(i,j)|  \leq    4 \sqrt{2} - 5 \approx 0.6569$  for all $1 \leq i, j \leq p$, $i \neq j$. Then as $p \goto \infty$,
$\hamm_p^*(\vartheta, \kappa,  r, a, \Omega) = L_p p^{1 -(\vartheta + r)^2/(4r)}$.
\end{cor}
\begin{cor} \label{cor:entrywise2}
Suppose the conditions of Theorem \ref{thm:UB} hold.  Also, suppose that
$1 < r / \vartheta <  5  + 2 \sqrt{6} \approx 9.898$,  and that $\lambda_3^*(\Omega)  \geq 2(5 - 2\sqrt{6}) \approx 0.2021$, $\lambda_4^*(\Omega)  \geq 5 - 2 \sqrt{6} \approx 0.1011$, and
$|\Omega(i,j)| \leq   8 \sqrt{6}  - 19 \approx 0.5959$ for all $1 \leq i, j \leq p$, $i \neq j$.
Then as $p \goto \infty$,
$\hamm_p^*(\vartheta, \kappa,  r, a,   \Omega) = L_p p^{1 -(\vartheta + r)^2/(4r)}$.
\end{cor}
In these corollaries, the conditions on $\Omega$  are rather relaxed. Somewhat surprisingly,
the off-diagonals of $\Omega$ do not necessarily have a major influence on the optimal rate of convergence,  as one might have expected.

\begin{figure}[htb]
\centering
\includegraphics[width = 5.5 in,height = 2  in]{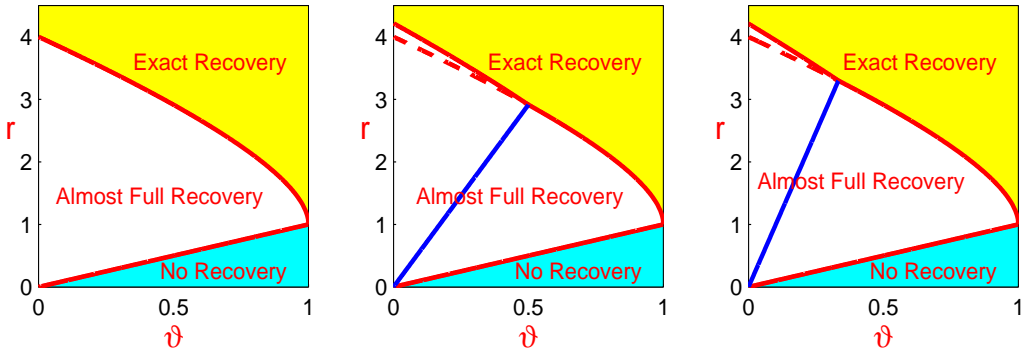}
\caption{Phase diagram for $\Omega = I_p$ (left),  for $\Omega$ satisfying conditions of Corollary \ref{cor:entrywise1} (middle), and  for $\Omega$ satisfying  conditions of Corollary \ref{cor:entrywise2} (right).  Red line: $r = \vartheta$.  Solid red curve: $r = \rho(\vartheta,  \Omega)$. In each of the last two panels, the blue line intersects with the red curve at $(\vartheta, r) = (1/2, [3 + 2 \sqrt{2}]/2)$ (middle) and $(\vartheta, r) = (1/3, [5 + 2 \sqrt{6}]/3)$ (right), which splits  the red solid curve
into two parts; the part  to the left is illustrative for it depends on $\Omega$ in a complicated way;   the part  to the right, together with the dashed red curve, represent $r = (1 + \sqrt{1 - \vartheta})^2$ (in the left panel, this is illustrated by  the red curve).}
\label{fig:phase}
\end{figure}


Note also that  by Theorem \ref{thm:UB}, under the condition of either Corollaries  \ref{cor:entrywise1} or Corollary \ref{cor:entrywise2},
GS achieves the optimal rate in that
\begin{equation} \label{gsadd}
\sup_{\mu \in \Theta_p^*(\tau_p, a)}  H_p(\hat{\beta}^{gs}; \eps_p, n_p, \mu, \Omega) \leq L_p p^{1 - (\vartheta + r)^2/(4r)}.
\end{equation}
Together, Theorem \ref{thm:universalLB},  Corollaries  \ref{cor:entrywise1}-\ref{cor:entrywise2}, and (\ref{gsadd})    have  an interesting implication on the so-called {\it  phase diagram}.
Call the two-dimensional {\it parameter space}  $\{(\vartheta, r): 0 < \vartheta < 1, r > 0\}$ the phase space.
There are two curves $r = \vartheta$ and $r = \rho(\vartheta,  \Omega)$  (the latter can be thought of as the solution of
$\sum_{j=1}^p p^{-\rho^*_j(\vartheta,r,a,\Omega)}=1$; recall that  $\rho_j^*(\vartheta, r, a, \Omega)$ does not depend on $a$)    that partition
the whole phase space into three different regions:
\begin{itemize}
\item {\it Region of No Recovery}.  $\{(\vartheta, r):     0 <  r < \vartheta,  0  < \vartheta < 1\}$.
In this region,  as $p \goto \infty$,  for any $\Omega$ and any procedures, the minimax Hamming error equals approximately to the total expected number of signals.
This is the most difficult region, in which no procedure can be successful in the minimax sense.
\item {\it Region of Almost Full Recovery}. $\{(\vartheta, r):   \vartheta < r < \rho(\vartheta,  \Omega) \}$.   In this region,  as $p \goto \infty$,   the minimax Hamming distance satisfies
$1 \ll \hamm_p^*(\vartheta, \kappa,  r, a, \Omega)  \ll p^{1 - \vartheta}$, and
 it is possible to recover most of the signals, but it is impossible to recover all of them.
\item {\it Region of Exact Recovery}.  In this region, as $p \goto \infty$,   the minimax Hamming distance $\hamm_p^*(\vartheta, \kappa, r, a, \Omega) = o(1)$, and it is possible to exactly recover all signals with overwhelming probability.
\end{itemize}
In general, the function $\rho(\vartheta, \Omega)$ depends on $\Omega$ in a complicated way.  However, by Theorem \ref{thm:universalLB} and Corollaries  \ref{cor:entrywise1}-\ref{cor:entrywise2}, we have the following conclusions.
First, for  all $\Omega$ and $a > 1$,
$\rho(\vartheta,  \Omega) \geq (1 + \sqrt{1 - \vartheta})^2$  for all $0 < \vartheta < 1$.
Second, in the simplest case where $\Omega = I_p$,   $\hamm_p^*(\vartheta, \kappa, r, a, \Omega) = L_p p^{1 - (\vartheta + r)^2/(4r)}$, and
$\rho(\vartheta,  \Omega) = (1 + \sqrt{1 - \vartheta})^2$  for all  $0 < \vartheta < 1$.
Third, under the conditions of Corollary \ref{cor:entrywise1},
$ \rho(\vartheta,  \Omega) =  (1 + \sqrt{1 - \vartheta})^2$ if  $1/2 < \vartheta < 1$.
Last,  under the conditions of Corollary \ref{cor:entrywise2},
$\rho(\vartheta,  \Omega) =  (1 + \sqrt{1 - \vartheta})^2$ if  $1/3 < \vartheta < 1$.
The phase diagram for the last three cases are illustrated in
Figure \ref{fig:phase}. The blue lines are  $r / \vartheta = 3 + 2 \sqrt{2}$ (middle) and $r/ \vartheta = 5 + 2\sqrt{6}$ (right).

Corollaries \ref{cor:entrywise1}-\ref{cor:entrywise2} can be extended to more general situations, where $r/ \vartheta$ may get arbitrary large, but
consequently, we need stronger conditions on $\Omega$.
Towards this end, we note that  for any $(\vartheta, r)$ such that $r > \vartheta$, we can find a unique integer $N  = N(\vartheta, r)$   such that
$2N - 1 \leq (\vartheta /r + r /\vartheta)/2 < 2N  + 1$.
Suppose that  for any $2 \leq k \leq 2N -1$,
\begin{equation} \label{cor:main1}
\lambda_k^*(\Omega)  \geq \max_{\{ (k+1)/2 \leq j \leq \min\{k, N\} \}}
\Big\{ \frac{(r/\vartheta + \vartheta/r)/2 - 2 j + 2+ \sqrt{[(r/\vartheta + \vartheta/r)/2 - 2j +2]^2 - 1}} { (2k - 2j  + 1)(r/\vartheta)} \Big\},
\end{equation}
and that for any $2 \leq k \leq 2 N $,
\begin{equation} \label{cor:main2}
\lambda_k^*(\Omega)   \geq  \max_{\{ k/2 \leq j \leq \min\{ k-1, N\} \}  }  \Big\{   \frac{(r/\vartheta + \vartheta/r)/2 + 1 - 2 j}{(k-j)(r/\vartheta)}  \Big\}.
\end{equation}
Then we have the following corollary.
\begin{cor}  \label{cor:main}
Suppose the conditions in Theorem \ref{thm:UB} and that in  (\ref{cor:main1})-(\ref{cor:main2}) hold. Then as $p \goto \infty$,
$\hamm_p^*(\vartheta, \kappa, r, a,   \Omega) = L_p p^{1-(\vartheta + r)^2/(4r)}$.
\end{cor}
The right hand sides of (\ref{cor:main1})-(\ref{cor:main2}) decrease with $(r/\vartheta)$. For a constant $s_0>1$, (\ref{cor:main1})-(\ref{cor:main2}) hold for all $1<r/\vartheta\leq s_0$ as long as they hold for $r/\vartheta= s_0$. Hence Corollary \ref{cor:main} implies a similar partition of the phase diagram as do Corollaries \ref{cor:entrywise1}-\ref{cor:entrywise2}.

{\bf Remark}.  Phase diagram can be viewed as a new criterion for
assessing the optimality, which is especially appropriate for
rare and weak signals.
The phase diagram is a partition of the phase space $\{(\vartheta, r): 0 < \vartheta < 1, r > 0\}$  into different regions
where statistical inferences are distinctly different.
In general, a phase diagram has the
following four regions:
\begin{itemize}
\item An ``exact recovery'' region corresponding to the ``rare and strong'' regime
in which high probability of completely correct variable selection is feasible.
\item An  ``almost full recovery'' region as a part of the ``rare and weak'' regime
in which completely correct variable selection is not
achievable with high probability but variable selection is still feasible in the sense that
with high probability, the number of incorrectly selected variables is a small fraction
of the total number of signals.
\item A ``detectable'' region in which variable selection is infeasible but the detection
of the existence of a signal (somewhere) is feasible (e.g.,  by the Higher Criticism method).
\item An ``undetectable'' region where  signals are so rare and weak that nothing can be sensibly done.
\end{itemize}

In the sparse signal detection \citep{DonohoJin2004} and classification \citep{PNAS} problems,
the main interest is to find the detectable region, so that the exact recovery and almost full recovery
regions were lumped into a single ``estimable'' region (e.g., \citet[Figure 1]{DonohoJin2004}).
For variable selection, the main interest is to find the boundaries of the almost full
discovery region so that the detectable and non-detectable regions are lumped
into a single ``no recovery'' region as in \citet{JiJin} and Figure 1 of this paper.

Variable selection in the ``almost full recovery'' region is a new and challenging problem. It was studied in \citet{JiJin} when the effect of signal cancellation is negligible,  but the hardest part of the problem was unsolved in \citet{JiJin}.
This paper (the second  in this area)  deals with the important issue of signal
cancellation, in hopes of gaining a much deeper insight on variable selection in much broader context.


\subsection{Non-optimaility of  subset selection and the lasso} \label{subsec:lasso}
Subset selection (also called the $L^0$-penalization method) is a well-known method for variable selection, which selects variables
by minimizing the following functional:
\begin{equation} \label{Definess}
\frac{1}{2} \| Y  - X \beta\|^2 + \frac{1}{2} (\lambda_{ss})^2 \| \beta\|_0,
\end{equation}
where $\|\beta\|_q$ denotes the $L^q$-norm, $q \geq 0$, and $\lambda_{ss}  >  0$ is a tuning parameter.    The AIC,
BIC, and RIC are methods of this type \citep{AIC, BIC, RIC}.
Subset selection is believed to have good ``theoretic property", but
the main drawback of this method is that it is computationally NP hard.
To overcome the computational challenge,
many {\it relaxation} methods are proposed, including but are not limited to the lasso  \citep{Chen, Tibshirani}, SCAD \citep{FanLi}, MC+ \citep{Zhang10}, and Dantzig selector \citep{CandesTao}.   Take the lasso for example. The method selects  variables by minimizing
\begin{equation} \label{Definelasso}
\frac{1}{2} \| Y  - X \beta\|^2 + \lambda_{lasso}  \| \beta\|_1,
\end{equation}
where the $L^0$-penalization is replaced by the $L^1$-penalization, so
the functional is convex and the optimization problem is solvable in polynomial time under proper conditions.

Somewhat surprisingly, subset selection is generally {\it rate non-optimal} in terms of selection errors. This sub-optimality of subset selection is due to its lack of flexibility in adapting to the ``local" graphic structure of the design variables. Similarly, other global relaxation methods are sub-optimal as well,  as the subset selection is the ``idol" these methods try to mimic.
To save space, we only discuss subset selection and the lasso, but a similar conclusion can be drawn for SCAD, MC+, and  Dantzig selector.

For mathematical simplicity, we illustrate the point with an idealized regression model where the Gram matrix $G = X'X$ is diagonal block-wise and  has  $2\times 2$ blocks
\begin{equation} \label{eq:2by2omega}
G(i,j) =  1 \{ i =  j\} + h_0 \cdot 1\{\mbox{$|j - i| = 1$,  $max(i,j)$ is even}  \}, \;\;   |h_0|  < 1, \; 1 \leq i, j \leq p.
\end{equation}
Using an idealized model is mostly for technical convenience, but the non-optimality
of  subset selection or  the lasso holds much more broadly than what is considered
here.  On the other hand,  using a simple model is sufficient here: if a procedure is
non-optimal in an idealized case, we can not expect it to be optimal in a more
general context.

At the same time, we continue to model $\beta$ with the Asymptotic Rare and Weak model ARW$(\vartheta, r, a,  \mu)$, but where we relax the assumption of $\mu \in \Theta_p^*(\tau_p, a)$ to that of  $\mu \in \Theta_p(\tau_p)$ so that the strength of each signal $\geq \tau_p$  (but there is no upper bound on the strength).  Consider a variable selection procedure $\hb^{\star}$, where $\star = gs, ss, lasso$, representing GS, subset selection, and the lasso (where the tuning parameters for each method  are ideally set; for the worst-case risk considered below,  the ideal tuning parameters depend on $(\vartheta, r, p, h_0)$ but do not depend on $\mu$).   Since the index groups $\{2j-1,2j\}$ are exchangeable in (\ref{eq:2by2omega}) and the ARW models,
the Hamming error of $\beta^{\star}$ in its worst case scenario has the form of $\sup_{\{ \mu \in \Theta_p(\tau_p) \}} H_p(\hb^{\star}; \eps_p,  \mu, G) = L_p p^{1 -\rho_{\star}(\vartheta, r, h_0)}$.

We now study $\rho_{\star}(\vartheta, r, h_0)$.  Towards this end, we first  introduce
$\rho_{lasso}^{(3)}(\vartheta,r, h_0) =  \bigl\{ (2 | h_0|)^{-1}  [(1- h_0^2)\sqrt{r}-\sqrt{  (1- h_0^2)(1-| h_0|)^2r-4| h_0|(1-| h_0|)\vartheta    }]  \bigr\}^2$
and
$\rho_{lasso}^{(4)}(\vartheta,r, h_0) = \vartheta + \frac{(1 -  |h_0|)^3(1+|h_0|)}{16   h_0^2}   \bigl[(1+ |h_0|)\sqrt{r}-\sqrt{(1- |h_0|)^2r-4| h_0|\vartheta/(1- h_0^2)   } \bigr]^2$.
We then let
\[
\rho_{ss}^{(1)}(\vartheta,r, h_0)= \left\{
 \begin{array}{ll}
 2\vartheta, & \qquad \quad r/\vartheta\leq2/(1- h_0^2) \\
\, [2\vartheta + (1- h_0^2) r]^2/[4 (1- h_0^2) r ], & \qquad \quad r/\vartheta>2/(1- h_0^2)
\end{array},
 \right.
 \]
 \[
 \rho_{ss}^{(2)}(\vartheta,r, h_0)= \left\{
 \begin{array}{ll}
 2\vartheta, & \quad r/\vartheta\leq2/(1-| h_0|) \\
  2[\sqrt{2 (1-| h_0|) r}-\sqrt{ (1-| h_0|) r -\vartheta}]^2, & \quad r/\vartheta>2/(1-| h_0|)
 \end{array},
 \right.
 \]
 \[
  \rho_{lasso}^{(1)}(\vartheta,r, h_0)=\left\{
  \begin{array}{ll}
  2\vartheta, & \qquad \qquad  r/\vartheta\leq2/(1-| h_0|)^2   \\
\rho_{lasso}^{(3)}(\vartheta,r, h_0), & \qquad \qquad  r/\vartheta>2/(1-| h_0|)^2
\end{array},
 \right.
\]
and
 \[
 \rho_{lasso}^{(2)}(\vartheta,r, h_0)=\left\{
 \begin{array}{ll}
 2\vartheta, & \qquad \qquad r/\vartheta\leq(1+| h_0|)/(1-| h_0|)^3   \\
  \rho_{lasso}^{(4)}(\vartheta,r, h_0), & \qquad \qquad r/\vartheta>(1+| h_0|)/(1-| h_0|)^3
 \end{array}.
 \right.
 \]
The following theorem is proved in Section \ref{sec:appen}.
\begin{thm} \label{thm:compareGsSsLasso}
Fix $\vartheta \in (0,1)$ and $r > 0$ such that $r >  \vartheta$.
Consider Model (\ref{model}) where $\beta$ is modeled by  $ARW(\vartheta, r, a, \mu)$ and
$X$ satisfies (\ref{eq:2by2omega}).  For GS, we set the tuning parameters
$(\delta, m_0) = (0,2)$, and  set $({\cal Q}, u^{gs}, v^{gs})$  as in (\ref{Definet})-(\ref{Tuneuv}). For subset selection as in (\ref{Definess}) and the lasso as in (\ref{Definelasso}), we  set their tuning parameters ideally given that $(\vartheta, r)$ are known. Then as $p \goto \infty$,
\begin{equation}\label{eq:blockGSrate}
\rho_{gs}(\vartheta, r, h_0) =   \min \bigl\{\frac{(\vartheta +r)^2}{4r},  \,   \vartheta + \frac{(1 - | h_0|)}{2} r,   \,  2 \vartheta +   \frac{\{[(1 - h_0^2) r  - \vartheta]_+\}^2}{4 (1 -  h_0^2) r}  \bigr\},
\end{equation}
\begin{equation}\label{eq:blockSSrate}
\rho_{ss}(\vartheta, r,  h_0) =   \min
\bigl\{\frac{(\vartheta +r)^2}{4r},  \,   \vartheta + \frac{(1 - | h_0|)}{2} r, \,   \rho_{ss}^{(1)}(\vartheta, r,  h_0), \rho_{ss}^{(2)}(\vartheta, r,  h_0)  \bigr\},
\end{equation}
and
\begin{equation}\label{eq:blockLassorate}
\hspace{-0.3 em} \rho_{lasso}(\vartheta, r,  h_0) =  \min\{ \frac{(\vartheta +r)^2}{4r},  \,   \vartheta + \frac{(1 - | h_0|)r}{2(1+\sqrt{1- h_0^2})} , \,  \rho_{lasso}^{(1)}(\vartheta, r,  h_0), \, \rho_{lasso}^{(2)}(\vartheta, r,  h_0)  \bigr\}.
\end{equation}
\end{thm}

It can be shown that $\rho_{gs}(\vartheta, r,  h_0)\geq \rho_{ss}(\vartheta, r,  h_0)\geq \rho_{lasso}(\vartheta, r,  h_0)$, where depending on the choices of $(\vartheta, r,  h_0)$, we may
have equality or strict inequality (note that a larger exponent  means a better error rate).  This fits well with our expectation, where as far as the convergence rate is concerned,   GS is optimal for all  $(\vartheta, r,  h_0)$, so it  outperforms  the subset selection, which in turn  outperforms  the lasso.
Table \ref{tab:optimalRates} summarizes the exponents for some representative $(\vartheta, r,  h_0)$.   It is seen that differences between these exponents become increasingly prominent when $ h_0$ increase  and $\vartheta$ decrease.
\begin{table}[h]
\begin{tabular}{|c|c|c|c|c|c|c|c|c|}
 \hline
$\vartheta/r / h_0$  &.1/11/.8 &.3/9/.8 & .5/4/.8& .1/4/.4  & .3/4/.4  & .5/4/.4&.1/3/.2 & .3/3/.2
 \\
 \hline
$\star = gs$ &     1.1406  &  1.2000  &  0.9000  &  0.9907  &  1.1556  &  1.2656   & 0.8008  &  0.9075\\
$\star = ss$  &  0.8409  &  0.9047  &  0.9000  &  0.9093  &  1.1003 &   1.2655  &  0.8007  &  0.9075 \\
$\star=lasso$  &  0.2000  &  0.6000  &  0.7500  &  0.4342   & 0.7121   & 1.0218  &  0.6021  &  0.8919\\
 \hline
\end{tabular}
  \caption{The exponents $\rho_{\star}(\vartheta, r,  h_0)$ in Theorem \ref{thm:compareGsSsLasso}, where $\star = gs, ss, lasso$.}
  \label{tab:optimalRates}
\end{table}

\begin{figure}[htb]
\begin{center}
\includegraphics[width=4.6 in, height = 4 in]{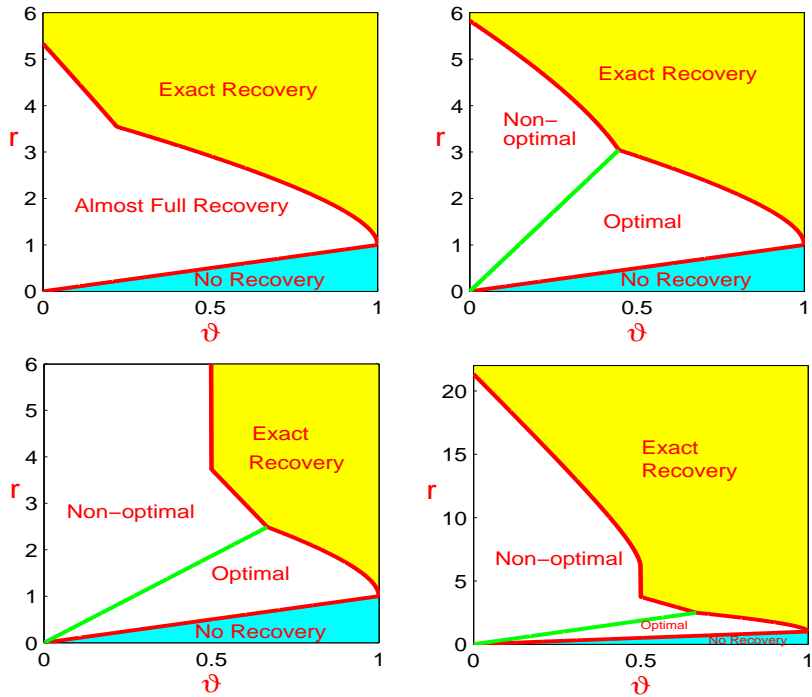}
\caption{Phase diagrams for GS (top left), subset selection (top right),  and the lasso (bottom; zoom-in on the left and zoom-out on the right), where $h_0=0.5$.}
\label{fig:VS}
\end{center}
\end{figure}

 As in Section \ref{subsec:simpl}, each of these methods has a phase diagram plotted in Figure \ref{fig:VS},
where the phase space partitions into three regions:  {\it Region of Exact Recovery},
{\it Region of Almost Full Recovery}, and {\it Region of No Recovery}.
Interestingly, the separating boundary for the last two regions are the same
for three methods, which is the line $r = \vartheta$.   The boundary that separates the first two regions,  however, vary significantly for different methods.   For any $ h_0\in (-1, 1)$ and $\star = gs, ss, lasso$,  the equation for this
boundary can be obtained by setting $\rho_{\star}(\vartheta, r,  h_0) = 1$ (the calculations are elementary so we omit them). Note that
the lower the boundary is, the better the method is,  and that the boundary corresponding to the lasso is discontinuous at $\vartheta = 1/2$.   In  the non-optimal region of either subset selection or the lasso, the Hamming errors of the procedure are much smaller than $p \eps_p$, so the procedure gives ``almost full recovery"; however, the rate of  Hamming errors is slower  than that of the  optimal procedure,  so subset selection or the lasso is non-optimal in such regions.

Subset selection and the lasso are rate non-optimal for they are so-called {\it one-step} or {\it non-adaptive} methods \citep{JiJin},
which use only one tuning parameter, and which do  not adapt to the local  graphic structure. The non-optimality can be best illustrated with the diagonal block-wise model presented here, where each block is a $2 \times 2$ matrix.
Correspondingly, we can partition the vector $\beta$ into many size $2$ blocks, each of which is of the following three types (i)  those have  no signal, (ii) those have exactly one signal, and (iii)
those have two signals.  Take the subset selection for example. To best separate
(i) from (ii), we need to set the tuning parameter ideally. But such a tuning parameter
may not be the ``best" for separating (i) from (iii).  This explains  the non-optimality of subset selection.

Seemingly, more complicated penalization methods that use multiple tuning parameters may have better performance than the subset selection and the lasso.
However, it remains open how to design such extensions
to achieve the optimal rate for general cases. To save space,
we leave the study along this line to the future.

\subsection{Summary}
\label{subsec:summary}
We propose GS as a new approach to variable selection. The key methodological
innovation is to use the GOSD to guide the multivariate screening.
While a brute-force $m$-variate screening has a computation cost of $O(p^m+np)$,
GS only has a computation cost of $L_p np$ (excluding the
overhead of obtaining the GOSD), by utilizing graph sparsity. 
Note that when the design matrix $G$ is approximately banded, say,  all its
large entries are confined to a diagonal band with bandwidth $\leq K$, the 
overhead of GS can  be reduced to $O(npK)$. One such example is in Genome-Wide Association Study (GWAS),
where $G$ is the empirical Linkage Disequilibrium (LD)
matrix, and $K$  can be as small as a few tens.  
We remark that the lasso has a computational complexity of $O(npk)$, where $k$, dominated by the number steps requiring re-evaluation of the correlation between design vectors and updated residuals, could be smaller than the $L_p$ term for GS \citep{Wang}.

We use {\it asymptotic minimaxity of the Hamming distance} as the criterion for
assessing optimality.  Compared with existing literature on variable selection where we use the {\it oracle property} or {\it probability of exact support recovery} to assess optimality,
our approach is mathematically more demanding, yet  scientifically more relevant  in the rare/weak paradigm.

We have proved that GS achieves the optimal rate of convergence
of Hamming errors, especially
when signals are rare and weak, provided that the Gram matrix is sparse.
Subset
selection and the lasso  are not rate optimal, even
 with very simple Gram matrix $G$ and even when the tuning parameters
are ideally set.   The sub-optimality of these  methods is due to that they do not   take advantage of the
`local' graphical structure as GS does.

GS has three key tuning parameters:   $q$
for the threshold level $t(\hd,\he) = 2\sigma^2 q \log p$ in the $GS$-step,  and
$(u^{gs}, v^{gs})  = (\sigma\sqrt{2\vartheta\log p},   \sigma\sqrt{2r\log p})$ in the $GC$-step.
While the choice of $q$ is reasonably flexible and a sufficiently small fixed $q>0$ is usually adequate, the choice of $u^{gs}$ and $v^{gs}$
are more directly tied to the signal sparsity and signal strength.
Adaptive choice of these tuning parameters is a challenging direction of further research.
One of our ideas to be developed in this direction is a subsampling scheme similar to
the Stability Selection \citep{MeinB2010}. On the other hand, as shown in our numeric results in Section \ref{sec:Simul}, the performance of  GS  is relatively insensitive to mis-specification of $(\epsilon_p,\tau_p)$; see details therein.

\section{Properties of Graphlet Screening,  proof of Theorem \ref{thm:UB}}   \label{sec:main}
GS attributes the success to two important properties: the Sure Screening
property  and the {\it Separable After Screening}  (SAS) property.

The Sure Screening property
means that in the $m_0$-stage $\chi^2$ screening, by picking an appropriate threshold, the set ${\cal U}_p^*$ (which is the set of retained indices after the GS-step)  contains  all but a small fraction of true signals.
Asymptotically, this fraction is comparably smaller than the minimax Hamming errors, and so negligible.
The SAS property means that except for a negligible probability,   as a subgraph of the GOSD,  ${\cal U}_p^*$ decomposes into many disconnected components of the GOSD, where the size of each component does not exceed  a fixed integer.
These two properties ensure that the original regression problem reduces to
many small-size regression problems, and thus pave the way for the $GC$-step.

Below,
we explain these ideas in detail, and conclude the section by the proof of Theorem \ref{thm:UB}. Since the only place we need the knowledge of $\sigma$ is in setting the tuning parameters,  so without loss of generality, we assume $\sigma = 1$  throughout this section.

First, we discuss the $GS$-step.  For short,  write $\hb = \hb^{gs}(Y; \delta, \tq, u^{gs}, v^{gs}, X, p, n)$ throughout this section. We first discuss the computation cost of the $GS$-step. As in Theorem \ref{thm:UB}, we take the threshold $\delta$ in $\cg^{*,\delta}$  to be $\delta  = \delta_p = 1/\log(p)$.
The proof of the following lemma  is similar to  that of \citet[Lemma 2.2]{JiJin}, so we omit it.
\begin{lemma} \label{lemma:sparseGOSD}
Suppose the conditions of Theorem \ref{thm:UB} hold, where we recall
$\delta = 1/\log(p)$, and $\Omega^{*,\delta}$ is defined as in \eqref{DefineOmega*}. As $p \goto \infty$,
with probability $1 - o(1/p^2)$,  $\|\Omega - \Omega^{*,\delta} \|_{\infty} \leq C (\log(p))^{-(1 - \gamma)}$, and $\cg^{*,\delta}$ is $K$-sparse, where $K \leq C (\log(p))^{1/\gamma}$.
\end{lemma}
Combining Lemma \ref{lemma:sparseGOSD} and  \cite{Frieze}, it follows that
with probability $1 - o(1/p^2)$,  $\cg^{*,\delta}$ has at most
$p  (C e (\log(p))^{1/\gamma})^{m_0}$
connected subgraphs of size $\leq m_0$. Note that the second factor is at most logarithmically large, so  the computation cost in the $GS$-step is at most  $L_p p$ flops.

Consider the performance of the $GS$-step.
The goal of this step is two-fold: on one hand, it tries to retain as many signals as possible
during the screening; on the other hand, it tries to minimize the computation cost of the $GC$-step by  controlling the maximum size of all   components of ${\cal U}_p^*$.   The key in the $GS$-step is to set the collection of thresholds ${\cal Q}$. The tradeoff is that,  setting the thresholds too high may miss too many signals during the screening, and setting the threshold too low may increase the maximum size of the components in ${\cal U}_p^*$, and
so increase the computational burden of the $GC$-step.
The following lemma  characterizes   the Sure Screening property of GS, and  is  proved in Section \ref{sec:appen}.
\begin{lemma} \label{lemma:SS} ({\it Sure Screening}).
Suppose the settings and conditions are as in Theorem \ref{thm:UB}.
In the $m_0$-stage $\chi^2$ screening of  the  GS-step, if we set the thresholds $t(\hd, \he)$ as in (\ref{Definet}),  then as $p \goto \infty$,  for any $\Omega \in {\cal M}_p^*(\gamma,  c_0, g, A)$, $\sum_{j = 1}^p P(\beta_j \neq 0,j\notin{\cal U}_p^*) \leq L_p [p^{1-(m_0+1)\vartheta} +  \sum_{j = 1}^p p^{-\rho_j^*(\vartheta, r, a, \Omega)}] + o(1)$.
\end{lemma}

Next,  we formally state the SAS property. Viewing it  as a subgraph of $\cg^{*,\delta}$,    ${\cal U}_p^*$ decomposes  into many  disconnected components $ \call^{(k)}$,   $1 \leq k \leq N$, where $N$ is an integer that may depend on the data.
\begin{lemma}   \label{lemma:SAS}
({\it SAS}).
Suppose the settings and conditions are as in Theorem \ref{thm:UB}.
In the $m_0$-stage $\chi^2$ screening in the GS-step, suppose
we set  the thresholds $t(\hd, \he)$ as in (\ref{Definet}) such that $q(\hd, \he) \geq q_0$  for some constant $q_0 = q_0(\vartheta, r) > 0$. As  $p \goto \infty$, under the conditions of Theorem \ref{thm:UB},   for any
$\Omega \in {\cal M}_p^*(\gamma,c_0,g, A)$,     there is a constant $\ell_0  = \ell_0(\vartheta, r, \kappa, \gamma, A, c_0, g) > 0$ such that with probability at least $1 - o(1/p)$,
$| \call^{(k)} | \leq \ell_0$, $1 \leq k \leq N$.
\end{lemma}
We remark that a more convenient way of picking  $q$ is to let
\begin{equation} \label{Defineq2}
\left\{
\begin{array}{ll}
q_0 \leq q \leq  (\frac{\omega  r  +  \vartheta}{2 \omega  r})^2 \omega  r,
&\qquad  \mbox{$|\hat{D}|$ is odd \&  $\omega  r / \vartheta > |\hat{D}| + (|\hat{D}|^2 - 1)^{1/2}$},      \\
q_0 \leq q  \leq
\frac{1}{4} \omega  r,    &\qquad  \mbox{$|\hat{D}|$ is even \&  $\omega  r / \vartheta \geq 2 |\hat{D}|$},
\end{array}
\right.
\end{equation}
and let $q$ be any other number otherwise, with which
both lemmas continue to hold with this choice of $q$. Here, for short, $\omega = \omega(\hd, \he; \Omega)$.  Note that numerically
this choice is comparably more conservative.

Together, the above two  lemmas say that the $GS$-step makes only negligible false non-discoveries,    and decomposes ${\cal U}_p^*$ into
many disconnected components, each has a size not exceeding a fixed integer.
As a result, the computation cost of the following $GC$-step is moderate,
at least in theory.

We now discuss the $GC$-step. The key to understanding the $GC$-step is that
the original regression
problem reduces to many disconnected small-size regression problems.
To see the point,  define  $\ty = X' Y$ and recall that $G= X'X$.  Let $\call_0 \lhd {\cal U}_p^*$
be a component,  we limit our attention to $\call_0$ by considering  the following regression problem:
\begin{equation} \label{sec21}
\ty^{\call_0} =  G^{\call_0, \otimes}   \beta + (X' z)^{\call_0},
\end{equation}
where  $(X'z)^{\call_0} \sim N(0, G^{\call_0, \call_0})  \approx  N(0,   \Omega^{\call_0, \call_0})$,  and $G^{\call_0,\otimes}$ is a $|\call_0| \times p$ matrix according to our notation.  What is non-obvious here is that,  the regression problem still involves the whole vector $\beta$,  and is still  high-dimensional.  To see the point, letting $V = \{1, 2, \ldots, p\} \setminus {\cal U}_p^*$, we  write $G^{\call_0,\otimes}  \beta = G^{\call_0, \call_0} \beta^{\call_0} +  I + II$, where $I = \sum_{{\cal J}_0: {\cal J}_0 \lhd {\cal U}_p^*,   {\cal J}_0 \neq \call_0} G^{\call_0, {\cal J}_0} \beta^{{\cal J}_0}$ and $II =  G^{\call_0, V} \beta^{V}$.
First, by Sure Screening property, $\beta^{V}$ contains only a negligible number of signals,   so we can think $II$  as negligible. Second, for any ${\cal J}_0 \neq \call_0$ and ${\cal J}_0 \lhd {\cal U}_p^*$,   by the SAS property, $\call_0$ and ${\cal J}_0$ are disconnected and so the matrix $G^{\call_0, {\cal J}_0}$ is a small size matrix whose coordinates are uniformly small.
This heuristic is  made precise in the proof of Theorem \ref{thm:UB}.
It is now seen  that the regression problem in (\ref{sec21}) is indeed low-dimensional:
\begin{equation} \label{sec22}
\ty^{\call_0} \approx  G^{\call_0, \call_0} \beta^{\call_0} + (X' z)^{\call_0} \approx N(\Omega^{\call_0, \call_0} \beta^{\call_0}, \Omega^{\call_0, \call_0}),
\end{equation}
The above argument is made precise in Lemma \ref{lemma:A.4}, see details therein.
Finally,  approximately,  the $GC$-step is to minimize
$\frac{1}{2}(\ty^{\call_0} - \Omega^{\call_0, \call_0} \xi)'  (\Omega^{\call_0, \call_0})^{-1} (\ty^{\call_0} - \Omega^{\call_0, \call_0} \xi) +\frac{1}{2} (u^{gs})^2 \|\xi\|_0$,
where each coordinate of $\xi$ is either $0$ or $\geq v^{gs}$ in magnitude.
Comparing   this with (\ref{sec22}), the procedure is nothing but the penalized MLE of  a low dimensional normal model, and the main result follows by exercising basic statistical inferences.

We remark that in the $GC$-step,  removing the constraints on the coordinates of $\xi$ will not
give the optimal rate of convergence. This is one of  the reasons why the classical
subset selection procedure is rate non-optimal. Another reason why the subset selection
is non-optimal is that, the procedure has only one tuning parameter,  but GS
has the flexibility of using different tuning parameters in the $GS$-step and the $GC$-step.
See Section \ref{subsec:lasso} for more discussion.

We are now ready for the proof of  Theorem \ref{thm:UB}.
\subsection{Proof of Theorem  \ref{thm:UB}}
For notational simplicity, we write $\rho_j^* = \rho_j^*(\vartheta, r, a, \Omega)$.
By Lemma \ref{lemma:SS},
\begin{equation}
\sum_{j = 1}^p P(\beta_j \neq 0,j\notin{\cal U}_p^*) \leq L_p [p^{1-(m_0+1)\vartheta} +  \sum_{j = 1}^p p^{-\rho_j^*}] + o(1).
\end{equation}
So to show the claim, it is sufficient to show
\begin{equation} \label{thmUB2}
\sum_{j=1}^p P(j \in {\cal U}_p^*,\sgn(\beta_j)\neq \sgn(\hb_j) ) \leq      L_p [ \sum_{j = 1}^p p^{-\rho_j^*} +  p^{1-(m_0 + 1) \vartheta}] + o(1).
\end{equation}
Towards this end, let $S(\beta)$ be the support of $\beta$,
$\Omega^{*,\delta}$ be as in \eqref{DefineOmega*}, and $\cg^{*,\delta}$ be the GOSD.
Let ${\cal U}_p^*$ be the set of retained indices after the GS-step.
Note that when $\sgn(\hb_j) \neq 0$,   there is a unique  component $\call_0$ such that $j \in \call_0 \lhd  {\cal U}_p^*$. For any connected subgraph  $\call_0$ of $\cg^{*,\delta}$, let
$B(\call_0)  = \{\mbox{$k$:  $k \notin \call_0$,   $\Omega^{*,\delta}(k, \ell) \neq 0$ for some $\ell \in \call_0$, $1 \leq k \leq p$} \}$.
Note that when $\call_0$ is a component of ${\cal U}_p^*$, we must have $B(\call_0) \cap {\cal U}_p^*=\emptyset$ as
for any node in $B(\call_0)$, there is at least one edge between it and some nodes in the component $\call_0$.
As a result,
\begin{equation} \label{thmUB3}
P(j \in \call_0 \lhd {\cal U}_p^*, B(\call_0) \cap S(\beta) \neq \emptyset) \leq \sum_{\call_0: j \in  \call_0} \sum_{k \in B(\call_0)}  P(k \notin {\cal U}_p^*, \beta_k \neq 0),
\end{equation}
where the first summation is over all connected subgraphs that contains node $j$.
By Lemma \ref{lemma:SAS}, with probability at least $1 - o(1/p)$, $\cg^{*,\delta}$ is $K$-sparse with $K = C (\log(p))^{1/\gamma}$, and there is a finite integer $\ell_0$ such that $|\call_0| \leq \ell_0$.
As a result, there are at most finite $\call_0$ such that the event $\{j \in \call_0 \lhd {\cal U}_p^*\}$ is non-empty, and for each of such $\call_0$, $B(\call_0)$ contains at most $L_p$ nodes.
Using (\ref{thmUB3}) and Lemma \ref{lemma:SS}, a direct result is
\begin{equation} \label{thmUB4}
\sum_{j = 1}^p P( j \in \call_0 \lhd {\cal U}_p^*, B(\call_0) \cap S(\beta) \neq \emptyset)  \leq L_p [\sum_{j = 1}^p p^{-\rho_j^*} + p^{1-(m_0+1) \vartheta}]  + o(1).
\end{equation}
Comparing  (\ref{thmUB4}) with   (\ref{thmUB2}),  to show the claim, it is sufficient to show that 
\begin{equation} \label{thmUB5}
\sum_{j = 1}^p P(\sgn(\beta_j)\neq \sgn(\hb_j), j \in \call_0 \lhd {\cal U}_p^*, B(\call_0) \cap S(\beta) =  \emptyset)  \leq L_p [  \sum_{j = 1}^p p^{-\rho_j^*} + p^{1-(m_0 + 1) \vartheta} ] + o(1).
\end{equation}
Fix $1 \leq j \leq p$ and a connected subgraph $\call_0$ such that $j \in \call_0$.  For short,
let $S$ be the support of $\beta^{\call_0}$ and $\hat{S}$ be the support of $\hb^{\call_0}$. The event  $\{\sgn(\beta_j)\neq \sgn(\hb_j), j \in \call_0 \lhd {\cal U}_p^* \}$ is identical to the event of $\{\sgn(\beta_j)\neq \sgn(\hb_j), j \in S\cup \hat{S} \}$. Moreover,
Since $\call_0$ has a finite size, both $S$ and $\hat{S}$ have finite possibilities.
So to show (\ref{thmUB5}), it is sufficient to show that for any fixed $1 \leq j \leq p$,  connected subgraph $\call_0$,  and  subsets $S_0, S_1 \subset \call_0$ such that $j \in S_0 \cup S_1$,
\begin{equation} \label{thmUB6}
P(\sgn(\beta_j)\neq\sgn(\hb_j),S = S_0, \hat{S} = S_1,j \in\call_0 \lhd {\cal U}_p^*, B(\call_0) \cap S(\beta) =  \emptyset) \leq L_p [ p^{-\rho_j^*} + p^{-(m_0 + 1) \vartheta}].
\end{equation}

We now show (\ref{thmUB6}).
The following lemma is proved in \citet[A.4]{JiJin}.
\begin{lemma} \label{lemma:A.4}
Suppose the conditions of Theorem \ref{thm:UB} hold. Over the event $\{j \in \call_0 \lhd {\cal U}_p^*\} \cap \{ B(\call_0) \cap S(\beta) = \emptyset \}$,
$\|(\Omega \beta)^{\call_0} - \Omega^{\call_0, \call_0} \beta^{\call_0} \|_{\infty} \leq C \tau_p (\log(p))^{-(1 - \gamma)}$.
\end{lemma}
Write for short
$\hat{M} = G^{\call_0, \call_0}$ and $M = \Omega^{\call_0, \call_0}$.   By definitions, $\hb^{\call_0}$ is the minimizer of the following functional
$Q(\xi) \equiv  \frac{1}{2}  (\ty^{\call_0} - \hat{M} \xi)' \hat{M}^{-1} (\ty^{\call_0} - \hat{M} \xi) +  \frac{1}{2} (u^{gs})^2 \|\xi\|_0$,
where $\xi$ is an $|\call_0| \times 1$ vector whose coordinates are either $0$ or $\geq v^{gs}$ in magnitude,  $u^{gs}  = \sqrt{2 \vartheta \log(p)}$, and
$v^{gs} = \sqrt{2 r \log(p)}$.  In particular,
$Q(\beta^{\call_0} ) \geq  Q(\hb^{\call_0})$,
or equivalently
\begin{equation} \label{thmUB7}
(\hb^{\call_0}  - \beta^{\call_0})'  (\ty^{\call_0} - \hat{M} \beta^{\call_0}) \geq  \frac{1}{2} (\hb^{\call_0} - \beta^{\call_0})' \hat{M} (\hb^{\call_0}  - \beta^{\call_0})  +  (|S_1| - |S_0|) \vartheta \log(p).
\end{equation}
Now, write for short $\delta  =  \tau_p^{-2} (\hb^{\call_0} - \beta^{\call_0})' M (\hb^{\call_0} - \beta^{\call_0})$.
First, by Schwartz inequality,  $[(\hb^{\call_0}  - \beta^{\call_0})'  (\ty^{\call_0} - \hat{M} \beta^{\call_0})]^2  \leq \delta \tau_p^2  (\ty^{\call_0} -   \hat{M} \beta^{\call_0}  )'  M^{-1} (\ty^{\call_0} - \hat{M} \beta^{\call_0})$.
Second, by Lemma \ref{lemma:A.4}, $\ty^{\call_0} = w + M \beta^{\call_0} + rem$, where
$w \sim N(0, M)$ and with probability $1 - o(1/p)$,
$|rem| \leq C (\log(p))^{-(1 - \gamma)} \tau_p$.   Last,  with probability at least $(1 - o(1/p))$, $\parallel\hat{M} - M\parallel_{\infty} \leq C \sqrt{\log(p)} p^{-[\kappa - (1 - \vartheta)]/2}$.   Inserting these into (\ref{thmUB7}) gives that with probability at least $(1 - o(1/p))$, $w' M^{-1} w  \geq  \frac{1}{4} \biggl[ \bigl(\sqrt{\delta r } + \frac{(|S_1| - |S_0|)  \vartheta}{\sqrt{\delta r }} \bigr)_+ \biggr]^2  (2 \log(p)) +  O((\log(p))^{ \gamma})$.
Since $\gamma<1$, $O((\log(p))^\gamma)$ is negligible. We note that $w' M^{-1} w \sim \chi_{|\call_0|}^2(0)$.  Inserting this back to (\ref{thmUB6}), the left hand side
$\leq   \eps_p^{|S_0|} P(\chi_{|\call_0|}^2(0) \geq [(\sqrt{\delta r } + (|S_1| - |S_0|)  \vartheta/\sqrt{\delta r })_+]^2  (\log(p)/2)) + o(1/p)$.
Assume $\sgn(\beta_j) \neq \sgn(\hb_j)$, and fix all parameters except $\delta$, $S_0$ and $S_1$. By arguments similar to the proof of Lemma \ref{lemma:omega}, the above quantity cannot achieve its maximum in the cases where $S_0= S_1$. Hence we only need to consider the cases where $S_0\neq S_1$. We also only need to consider the cases where $\max(|S_0|,|S_1|)\leq m_0$, since the sum of the probabilities of other cases is  controlled by $p^{1-(m_0+1)\vartheta}$. The claim follows by the definitions of $\rho_j^*$.
\qed

\section{Simulations}  \label{sec:Simul}
We  conduct a small-scale simulation study  to investigate  the numerical performance of Graphlet Screening and compare it  with the lasso  and the UPS.
The subset selection is not included for comparison since  it is computationally
NP hard.   We consider the experiments for both random design and fixed design, where as before, the parameters $(\eps_p, \tau_p)$ are tied to $(\vartheta, r)$ by $\eps_p = p^{-\vartheta}$ and $\tau_p = \sqrt{2 r \log(p)}$ (we assume $\sigma = 1$ for simplicity in this section).

In random design settings where $p$ is not very large, we follow the spirit of the refined UPS in \cite{JiJin} and  propose the  iterative Graphlet Screening algorithm where
we iterate Graphlet Screening for a few times ($\leq 5$).  The main purpose
for the iteration is to denoise the Gram matrix; see \citet[Section 3]{JiJin} for
more discussion.

Even with the refinement as in \citet[Section 3]{JiJin},  UPS behaves
poorly for most examples presented below.
Over close investigations, we find out that this is due to the threshold choice
in the initial $U$-step is too low,  and increasing
the threshold largely increases the performance.
Note that the purpose of this step is to denoise
the Gram matrix \citet[Section 3]{JiJin}, not for signal retainment, and
so a larger threshold helps.

In this section, we use this improved version of refined UPS, but for
simplicity, we still call it the refined UPS.
With that being said, recall that  UPS is unable to resolve the problem
of signal cancellation, so it usually performs poorer than
GS, especially when the effect of signal cancellation is strong.
For this reason, part of the comparison is between GS  and the lasso only.

The experiments with random design  contain the following steps.
\begin{enumerate}
\item Fix $(p, \vartheta, r,  \mu, \Omega)$ such that $\mu \in \Theta_p(\tau_p)$.
Generate a vector $b = (b_1, b_2, \ldots, b_p)'$ such that $b_i \stackrel{iid}{\sim}  \mathrm{Bernoulli}(\eps_p)$, and set $\beta = b \circ \mu$.
\item Fix $\kappa$ and let $n = n_p = p^{\kappa}$. Generate an $n \times p$ matrix with $iid$ rows from  $N(0, (1/n) \Omega)$.
  \item Generate $Y \sim  N(X \beta, I_n)$, and apply the iterative Graphlet Screening, the refined UPS and the lasso.
  \item Repeat 1-3  independently, and  record the average Hamming distances or the Hamming ratio, the ratio of the Hamming distance and the number of the signals.
\end{enumerate}
The steps for fixed design experiments are similar, except for that  $n_p = p$, $X = \Omega^{1/2}$ and we apply GS  and UPS directly.

GS  uses tuning parameters $(m_0, \tq, u^{gs}, v^{gs})$.  We set $m_0 = 3$ for our experiments, which is usually large enough due to signal sparsity. The choice of $\tq$ is not critical, as long as the
corresponding parameter $q$ satisfies (\ref{Defineq}), and we use the maximal $\tq$ satisfying \eqref{Defineq} in most experiments. Numerical studies below (e.g.,  Experiment $5a$) support this point.
In principle,  the optimal choices of $(u^{gs}, v^{gs})$ depend on the unknown parameters $(\eps_p, \tau_p)$,
and how to estimate them in general settings is a lasting open problem (even for linear models with orthogonal designs). Fortunately,  our studies (e.g., Experiment 5b-5d)
show that  mis-specifying  parameters $(\eps_p, \tau_p)$ by a reasonable amount  does not significantly affect the performance of the procedure.
For this reason, in most experiments below, assuming  $(\eps_p, \tau_p)$  are known, we set $(u^{gs}, v^{gs})$ as  $(\sqrt{2 \log(1/\eps_p)}, \tau_p)$. For the iterative Graphlet Screening, we use the same tuning parameters in each iteration.

For the UPS and the refined UPS, we use the tuning parameters $(u^{ups}, v^{ups})=(u^{gs}, v^{gs})$.
For both the iterative Graphlet Screening and the refined UPS,
we use the following as the initial estimate:  $\hat{\beta}_i = \sgn(\tilde{Y}_i)   \cdot 1 \{ |\tilde{Y}_i| \geq \tau_p\}$, $1 \leq i \leq p$, where $\tilde{Y} = X'Y$.
The main purpose of initial estimate is to denoise the Gram matrix, not for screening.
We use {\it glmnet} package \citep{Friedman2010} to perform lasso. To be fair in comparison, we apply the lasso with all  tuning parameters, and we report the  Hamming error associated with the ``best"  tuning parameter.

The simulations contain $6$ different experiments which we now describe separately.

{\it Experiment 1}.
The goal of this experiment is two-fold. First, we compare  GS  with UPS and the lasso in the fixed design setting. Second, we investigate the minimum signal strength levels $\tau_p$ required by these three methods to yield exact recovery, respectively.

Fixing $p=0.5\times10^4$, we let  $\eps_p  = p^{ - \vartheta}$ for $\vartheta \in\{ 0.25,0.4, 0.55\}$,
and $\tau_p \in \{6,7,8,9,10\}$.  We use a fixed design model where  $\Omega$ is a symmetric diagonal block-wise matrix, where each block is a $2 \times 2$ matrix, with $1$ on the diagonals, and $\pm 0.7$ on the
off-diagonals (the signs alternate across different blocks). Recall the $\beta = b\circ \mu$. For
each pair of $(\eps_p, \tau_p)$, we generate $b$ as $p$ iid samples from $Bernoulli(\eps_p)$, and  we let $\mu$ be the vector where the signs of $\mu_i  = \pm 1$ with equal probabilities,
and $|\mu_i|  \stackrel{iid}{\sim}    0.8 \nu_{\tau_p} + 0.2 h$, where $\nu_{\tau_p}$ is the point mass  at   $\tau_p$ and $h(x)$ is the density of  $\tau_p(1 + V/6)$ with $V \sim \chi_1^2$.
The average Hamming errors across  $40$ repetitions  are tabulated in Table \ref{tab:e1}.  For all $(\vartheta, \tau_p)$ in this experiment, GS  behaves more
satisfactorily than the UPS, which in turn behaves more satisfactorily
than the lasso.

Suppose we say a method yields `exact recovery' if the average  Hamming error $\leq 3$. Then, when $\vartheta = 0.25$, the minimum $\tau_p$ for  GS  to yield exact recovery is $\tau_p \approx 8$, but that for UPS and the lasso are much larger ($\geq 10$). For larger $\vartheta$, the differences are less prominent, but the pattern is similar.

The comparison between GS  and UPS is particularly interesting.
Due to the block structure of $\Omega$, as $\vartheta$ decreases, the signals become increasingly less sparse, and the effects of  signal cancellation   become increasingly stronger. As a result, the advantage of GS over the UPS becomes increasingly more prominent.

\begin{table}[h]
\begin{center}
\begin{tabular}{|c|c|c|c|c|c|c|}
 \hline
 & $\tau_p$  & 6 & 7 & 8 & 9 & 10 \\
 \hline
\multirow{3}{*}{$\vartheta=0.25$} & Graphic Screening  & 24.7750  &  8.6750  &  2.8250  &  0.5250  &  0.1250 \\
& UPS    & 48.5500 & 34.6250 &  36.3500 &  30.8750 &  33.4000 \\
    & lasso & 66.4750 &  47.7000 &  43.5250 &  35.2500  & 35.0500 \\
    \hline
\multirow{3}{*}{$\vartheta=0.40$} & Graphic Screening&  6.9500  &  2.1500  &  0.4000  &  0.0750  &  0.0500
  \\
& UPS   &  7.7500  &  4.0000  &  2.2000  &  2.7750 &   2.4250  \\
   &    lasso &  12.8750  &  6.8000  &  4.3250  &  3.7500  &  2.6750 \\
   \hline
\multirow{3}{*}{$\vartheta=0.55$} & Graphic Screening  &  1.8750  &  0.8000  &  0.3250  &  0.2250  &  0.1250 \\
& UPS    &  1.8750  &  0.8000  &  0.3250  &  0.2250  &  0.1250  \\
 &   lasso &   2.5000  &  1.1000  &  0.7750  &  0.2750  &  0.1250  \\
  \hline
\end{tabular}
\end{center}
  \caption{Comparison of average  Hamming errors (Experiment 1).}
  \label{tab:e1}
\end{table}

{\it  Experiment 2}.    In this experiment, we compare  GS, UPS and the lasso in the random design setting, and investigate the effect of signal cancellation on their performances. We fix $(p, \kappa, \vartheta, r) = (0.5\times10^4, 0.975, 0.35, 3)$, and assume $\Omega$ is blockwise diagonal. We generate $\mu$ as in Experiment $1$, but to better illustrate the difference between UPS and GS in the presence of signal cancellation, we generate the vector $b$ differently and allow it to depend on $\Omega$. The experiment contains $2$ parts, $2a$ and $2b$.

In Experiment 2a, $\Omega$ is the block-wise matrix where
each block is $2$ by $2$ matrix with $1$ on the diagonals and $\pm .5$ on the off diagonals (the signs alternate on adjacent blocks). According to the blocks in $\Omega$, the  set of indices $\{1, 2, \ldots, p\}$ are also partitioned into blocks accordingly.
For any fixed $\vartheta$ and $\eta \in \{0, 0.01, 0.02, 0.03, 0.04, 0.05, 0.06, 0.1, 0.2\}$, we randomly choose  $(1 - 2p^{-\vartheta})$ fraction of the blocks (of indices)
where $b$ is $0$ at both indices, $2(1 - \eta) p^{-\vartheta}$ fraction of the blocks where $b$ is $0$ at one index and $1$ at the other (two indices are equally likely to be $0$), $2\eta p^{-\vartheta}$ faction
of the blocks where $b$ is $1$ on both indices.

Experiment 2b has similar settings, where the difference is that (a) we choose
$\Omega$ to be a diagonal  block  matrix where each block is a $4$ by $4$ matrix (say, denoted by $A$) satisfying
$A(i,j) = 1\{i = j\} + 0.4  \cdot  1\{|i - j| = 1\} \cdot \sgn(6 - i -j)  + 0.05   \{|i -j| \geq 2\} \cdot \sgn(5.5 -   i -  j)$, $1 \leq i, j \leq 4$,
and (b) $(1 - 4p^{-\vartheta})$ is the fraction of blocks where $b$ is nonzero in $k=0$ indices,  $4(1-\eta)  p^{-\vartheta}$ is that for $k = 1$, and $4\eta p^{-\vartheta}$ is that for $k\in\{2,3,4\}$ in total. In a block where $\beta$ is nonzero at $k$ indices, all configurations with $k=1$ are equally likely, and all those with $k \in \{2,3,4\}$ are equally likely.

\begin{figure}[h]
        \centering
        \begin{subfigure}[h]{0.4\textwidth}
                \centering
  \includegraphics[width=\textwidth]{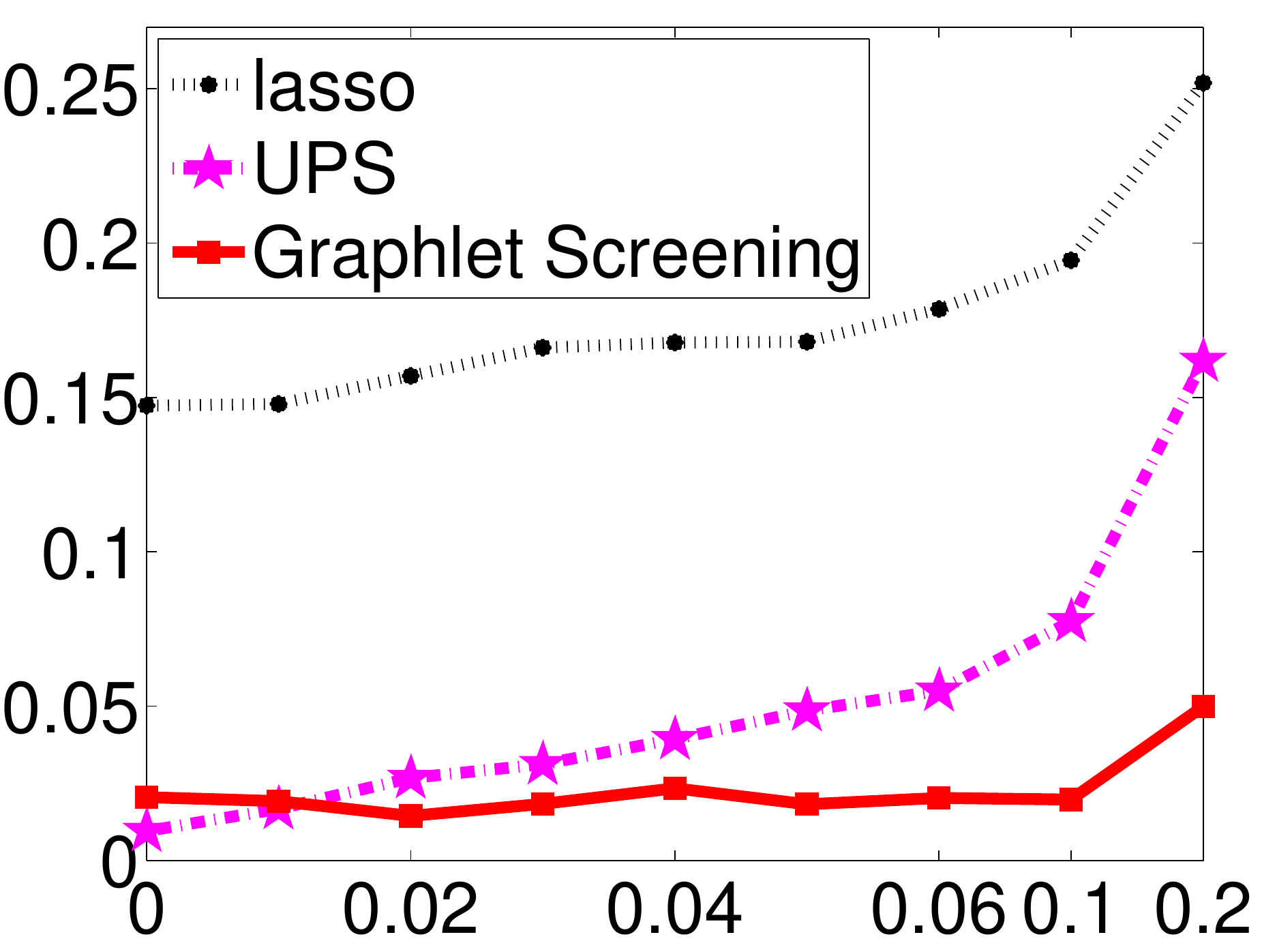}\\
  \caption{Blockwise diagonal $\Omega$ in $2 \times 2$ blocks}\label{fig:block2}
        \end{subfigure}%
         \begin{subfigure}[h]{0.4\textwidth}
                \centering
  \includegraphics[width=\textwidth]{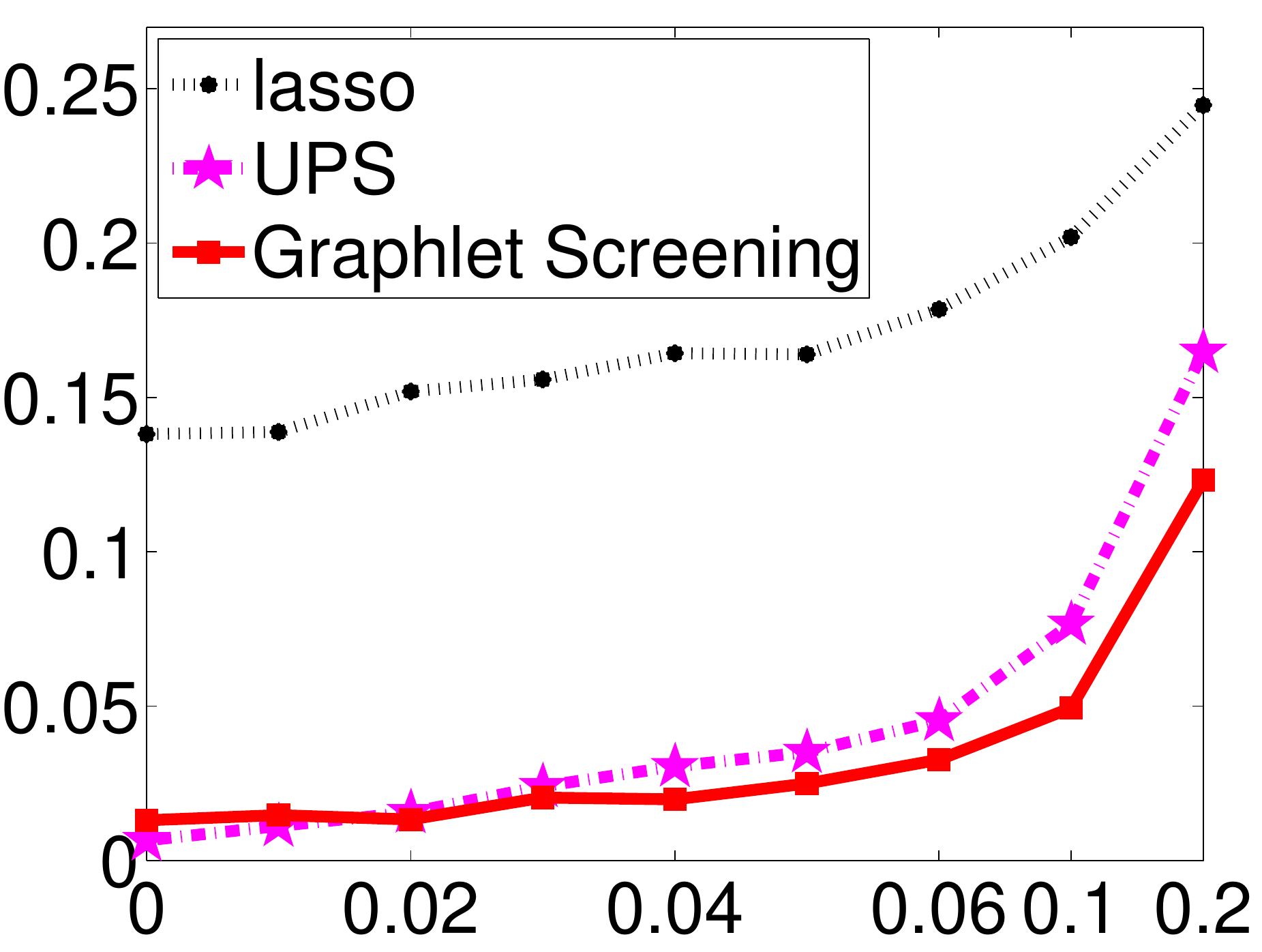}\\
  \caption{Blockwise diagonal $\Omega$ in $4 \times 4$ blocks}\label{fig:block4}
        \end{subfigure}
                \caption{Hamming ratio results in Experiment 2}\label{fig:e2}
\end{figure}

The average Hamming ratio results across $40$ runs for two Experiment 2a and 2b are reported in Figure \ref{fig:e2}, where UPS and GS   consistently outperform the lasso.  Additionally, when $\eta$ is small, the effect of signal cancellation
is negligible, so UPS and GS  have  similar performances. However, when $\eta$ increases, the effects of signal cancellation grows, and the advantage of GS  over UPS becomes increasingly more prominent.

Through  Experiment 1-2, the comparison of UPS and
GS  is more or less understood. For this reason,  we do not include
UPS for study in Experiment 3-5, but we include UPS for study in Experiment 6 where we investigate robustness of all three methods.

{\it  Experiment 3}. In this experiment, we investigate how different choices of signal vector $\beta$ affect the comparisons of  GS  and the lasso.
 We use a random design model, and $\Omega$ is a symmetric
 tri-diagonal correlation matrix  where the vector on each sub-diagonal consists of blocks of $(.4, .4, -.4)'$. Fix  $(p,\kappa)=(0.5\times10^4,0.975)$ (note $n  = p^{\kappa} \approx 4,000$).  We let $\eps_p = p^{-\vartheta}$ with $\vartheta \in \{0.35, 0.5 \}$ and let $\tau_p\in\{6,8,10\}$. For each combination of $(\eps_p, \tau_p)$,
we consider two choices of $\mu$. For the first choice, we
let $\mu$ be the vector where all coordinates equal to $\tau_p$ (note $\beta$ is still sparse).
For  the second one, we let $\mu$ be as in Experiment $1$.  The average Hamming ratios for both procedures across $40$ repetitions are tabulated in Table \ref{tab:e3}.
\begin{table}[h]
\begin{center}
\begin{tabular}{|c|c|c|c|c|c|c|c|}
  \hline
  \multicolumn{2}{|c|}{$\tau_p$} & \multicolumn{2}{|c|}{6} & \multicolumn{2}{|c|}{8} & \multicolumn{2}{|c|}{10}\\
  \hline
 \multicolumn{2}{|c|}{Signal Strength} & Equal & Unequal & Equal & Unequal & Equal & Unequal\\
 \hline
\multirow{2}{*}{$\vartheta=0.35$} & Graphic Screening &   0.0810  &  0.0825 &  0.0018 &   0.0034 &   0 &   0.0003  \\
    & lasso & 0.2424 &   0.2535  &   0.1445  &  0.1556  &   0.0941   &  0.1109 \\
    \hline
\multirow{2}{*}{$\vartheta=0.5$} & Graphic Screening & 0.0315  &  0.0297   &  0.0007   &   0.0007&   0 &   0 \\
   &    lasso & 0.1107 &  0.1130  &  0.0320  &  0.0254  &  0.0064  &  0.0115 \\
   \hline
\end{tabular}
\end{center}
  \caption{Hamming ratio results of Experiment $3$, where ``Equal" and ``Unequal" stand  for the first and the second choices of $\mu$, respectively.}
  \label{tab:e3}
\end{table}

{\it Experiment 4}. In this experiment,  we  generate $\beta$ the same way as in Experiment $1$,  and
investigate how different choices of design matrices affect the performance of the two methods. Setting $(p,\vartheta,\kappa)=(0.5\times10^4,0.35,0.975)$ and  $\tau_p\in\{6,7,8,9,10,11,12\}$,
we use Gaussian random design model for the study.
The experiment contains $3$ sub-experiments $4a$-$4c$.

\begin{figure}[htb]
\centering
\includegraphics[width=5.8 in, height = 2 in]{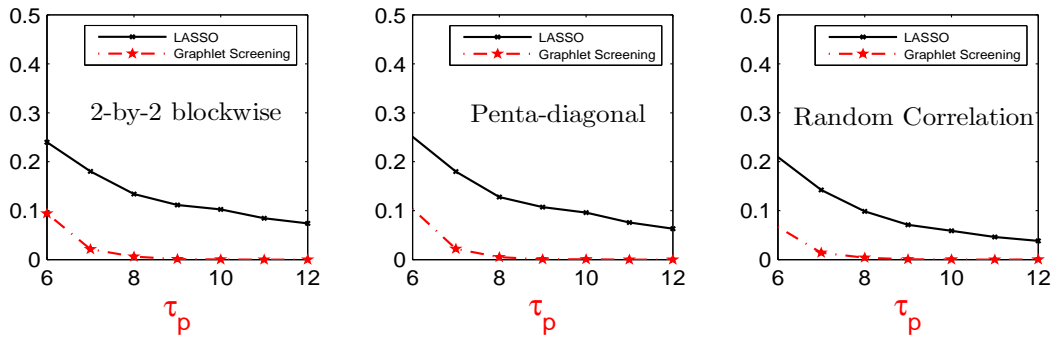}
\caption{$x$-axis: $\tau_p$. $y$-axis: Hamming ratios. Left to right:  Experiment $4a$, $4b$, and $4c$.  }\label{fig:4abc}
\end{figure}

In Experiment $4a$, we set $\Omega$ as
the symmetric diagonal block-wise matrix, where each block
is a $2 \times 2$ matrix, with $1$ on the diagonals, and $\pm 0.5$ on the
off-diagonals (the signs alternate across different blocks). The average Hamming   ratios  of $40$ repetitions are reported in Figure \ref{fig:4abc}.

In Experiment $4b$, we set $\Omega$ as a symmetric penta-diagonal correlation matrix,
where the main diagonal are ones, the first sub-diagonal consists of blocks of $(.4, .4, -.4)'$, and the second sub-diagonal consists of blocks of $(.05, -.05)'$. The average Hamming ratios across $40$ repetitions are reported in Figure \ref{fig:4abc}.

In Experiment $4c$, we generate $\Omega$ as follows. First, we generate $\Omega$ using the function \emph{sprandsym(p,K/p)} in \emph{matlab}. We then set the diagonals of $\Omega$  to be zero, and remove some of entries so that $\Omega$ is $K$-sparse for a pre-specified $K$.  We then normalize each non-zero entry by the sum of the absolute values in that row or that column, whichever is larger, and multiply each entry by a pre-specified positive constant $A$.  Last,  we set the diagonal elements to be 1.
We choose $K=3$ and $A=0.7$, draw $5$ different $\Omega$ with this method, and for each of them,
 we draw $(X, \beta, z)$
$10$ times independently.
The average Hamming ratios are reported in Figure \ref{fig:4abc}.
The results suggest that  GS  is consistently better than the lasso.

{\it Experiment 5.} In this experiment, we investigate how sensitive GS  is with respect to the tuning parameters.
The experiment contains $4$ sub-experiments, $5a$-$5d$. In Experiment $5a$,
we investigate how sensitive the procedure is with respect to the tuning parameter $q$ in
$\tq$ (recall that the main results hold as long as $q$ fall into the range  given in (\ref{Defineq})), where we assume $(\eps_p, \tau_p)$ are known.  In Experiment $5b$-$5d$,
 we mis-specify  $(\eps_p, \tau_p)$ by a reasonably small amount, and
investigate how the mis-specification affect the performance of the procedure.
For the whole experiment, we choose $\beta$ the same as in Experiment $1$, and $\Omega$ the same as in Experiment $4b$. We use a fixed design model in Experiment $5a$-$5c$, and a random design model  in Experiment $5d$. For each sub-experiment, the results are based on $40$ independent
repetitions. We now describe the sub-experiments with details.

In Experiment $5a$, we choose $\vartheta\in\{0.35,0.6\}$ and $r\in \{1.5,3\}$.
In GS, let $q_{max} = q_{max}(\hat{D}, \hat{F})$ be the maximum
value of $q$ satisfying (\ref{Defineq}).
For each combination of $(\vartheta,r)$ and $(\hat{D}, \hat{F})$,    we choose $q(\hat{D},\hat{F})=q_{max}(\hat{D},\hat{F})  \times  \{0.7,0.8,0.9,1,1.1,1.2\}$ for our experiment. The results are tabulated in Table \ref{tab:e5a}, which suggest that different choices of $q$
have little influence over the variable selection errors. We must note that the larger we set $q(\hat{D}, \hat{F})$, the faster the algorithm runs.
\begin{table}[h]
\begin{center}
\begin{tabular}{|c|c|c|c|c|c|c|}
 \hline
  $q(\hat{F},\hat{D})/q_{max}(\hat{F},\hat{D})$ & 0.7 & 0.8 & 0.9 & 1 & 1.1 & 1.2\\
 \hline
$(\vartheta,r)=(0.35,1.5)$  &  0.0782  &  0.0707  &  0.0661   & 0.0675  &  0.0684  &  0.0702\\
\hline
$(\vartheta,r)=(0.35,3)$  &  0.0066  &  0.0049  &  0.0036  &  0.0034  &  0.0033  &  0.0032\\
\hline
$(\vartheta,r)=(0.6,1.5)$  &  0.1417  &  0.1417  &  0.1417  &  0.1417   & 0.1417  &  0.1417\\
\hline
$(\vartheta,r)=(0.6,3)$  &  0.0089  &  0.0089  &  0.0089  &  0.0089  &  0.0089  &  0.0089 \\
 \hline
\end{tabular}
\end{center}
  \caption{Hamming ratio results in Experiment $5a$.}
  \label{tab:e5a}
\end{table}

In Experiment $5b$, we use the same settings as in Experiment $5a$, but
we assume $\vartheta$ (and so $\eps_p$)  is unknown (the parameter $r$ is assumed as known, however), and let $\vartheta^*$ is the misspecified value of $\vartheta$.
We take $\vartheta^*  \in   \vartheta\times\{0.85,0.925,1,1.075,1.15,1.225\}$ for the experiment.

In Experiment $5c$,  we use the same settings as in Experiment $5b$, but
we assume $r$ (and so $\tau_p$)  is unknown (the parameter $\vartheta$ is assumed as known, however), and let $r^*$ is the misspecified value of $r$.  We take $r^*=r\times\{0.8,0.9,1,1.1,1.2,1.3\}$ for the experiment.

In Experiment $5b$-$5c$,
we run GS  with tuning parameters set as in Experiment $1$, except
$\vartheta$ or $r$ are replaced by  the misspecified counterparts $\vartheta^*$ and $r^*$, respectively.
The results   are reported in Table \ref{tab:e5bc}, which suggest that
the mis-specifications have little effect as long as $r^*/r$ and $\vartheta^*/\vartheta$ are
reasonably close to $1$.
\begin{table}[h]
\begin{center}
\begin{tabular}{|c|c|c|c|c|c|c|}
 \hline
  $\vartheta^*/\vartheta$ & 0.85 & 0.925 & 1 & 1.075 & 1.15 & 1.225\\
 \hline
$(\vartheta,r)=(0.35,1.5)$  &  0.0799  &  0.0753  &  0.0711  &  0.0710  &  0.0715   & 0.0746\\
\hline
$(\vartheta,r)=(0.35,3)$  &  0.0026  &  0.0023  &  0.0029  &  0.0030  &  0.0031   & 0.0028\\
\hline
$(\vartheta,r)=(0.6,1.5)$  &  0.1468  &  0.1313  &  0.1272  &  0.1280 &   0.1247  &  0.1296\\
\hline
$(\vartheta,r)=(0.6,3)$  &  0.0122  &  0.0122  &  0.0139   & 0.0139  &  0.0130  &  0.0147 \\
 \hline
  \hline
  $r^*/r$ & 0.8 & 0.9 & 1 & 1.1 & 1.2 & 1.3\\
 \hline
$(\vartheta,r)=(0.35,1.5)$  &  0.0843  &  0.0731  &  0.0683  &  0.0645  &  0.0656   & 0.0687\\
\hline
$(\vartheta,r)=(0.35,3)$  &  0.0062  &  0.0039  &  0.0029  &  0.0030 &   0.0041  &  0.0054\\
\hline
$(\vartheta,r)=(0.6,1.5)$  &  0.1542  &  0.1365  &  0.1277  &  0.1237  &  0.1229  &  0.1261\\
\hline
$(\vartheta,r)=(0.6,3)$  &   0.0102  &  0.0076  &  0.0085 &   0.0059  &  0.0051  &  0.0076 \\
 \hline
\end{tabular}
\end{center}
  \caption{Hamming ratio results in Experiment $5b$ (top) and in Experiment $5c$ (bottom).}
   \label{tab:e5bc}
\end{table}

In Experiment $5d$, we re-examine the mis-specification issue in the random design setting. We use the same settings  as in Experiment $5b$ and Experiment $5c$, except for (a) while we use the same $\Omega$ as in Experiment $5b$, the design matrix  $X$ are generated according to the random design model as in Experiment $4b$,   and (b)
we only investigate for the case of $r=2$ and $\vartheta\in\{0.35,0.6\}$.   The results are summarized in Table \ref{tab:e5d}, which is consistent with the results in $5b$-$5c$.
\begin{table}[h]
\begin{center}
\begin{tabular}{|c|c|c|c|c|c|c|}
 \hline
  $\vartheta^*/\vartheta$ & 0.85 & 0.925 & 1 & 1.075 & 1.15 & 1.225\\
 \hline
$(\vartheta,r)=(0.35,2)$  &  0.1730  &  0.1367  &  0.1145 &   0.1118   & 0.0880  &  0.0983\\
\hline
$(\vartheta,r)=(0.6,2)$  &  0.0583  &  0.0591  &  0.0477  &  0.0487  &  0.0446  &  0.0431 \\
 \hline
  \hline
  $r^*/r$ & 0.8 & 0.9 & 1 & 1.1 & 1.2 & 1.3\\
 \hline
$(\vartheta,r)=(0.35,2)$  &  0.1881  &  0.1192  &  0.1275  &  0.1211  &  0.1474  &  0.1920\\
\hline
$(\vartheta,r)=(0.6,2)$  &   0.0813  &  0.0515 &   0.0536  &  0.0397  &  0.0442  &  0.0510 \\
 \hline
\end{tabular}
\end{center}
  \caption{Hamming ratio results in Experiment $5d$. }
  \label{tab:e5d}
\end{table}

{\it Experiment $6$.}  In this experiment, we investigate the robustness of all three methods for the mis-specification of the linear model \eqref{model}.  We use the random design setting as in Experiment $4b$, except that we fix $(\vartheta, r) = (0.35, 3)$.
The experiment contains $3$ sub-experiments, $6a$-$6c$,
where we consider three scenarios where the linear model
  (\ref{model}) is in question: the presence of nonGaussianity, the presence of
  missing predictors, and the presence of non-linearity, correspondingly.

\begin{figure}
        \centering
        \begin{subfigure}[h]{0.33\textwidth}
                \centering
  \includegraphics[width=\textwidth]{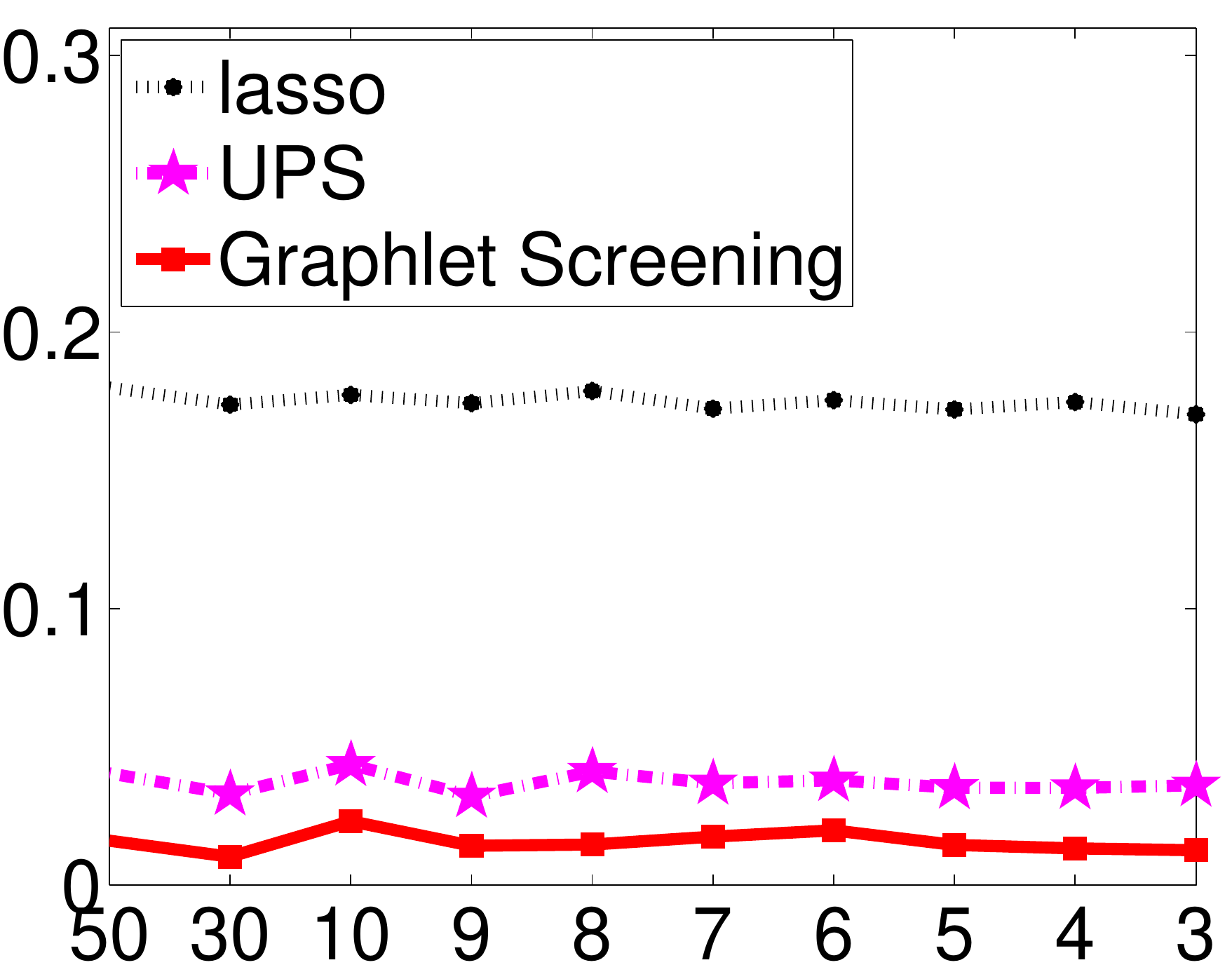}\\
  \caption{nonGaussian noise}\label{fig:tscaled}
        \end{subfigure}%
         \begin{subfigure}[h]{0.33\textwidth}
                \centering
  \includegraphics[width=\textwidth]{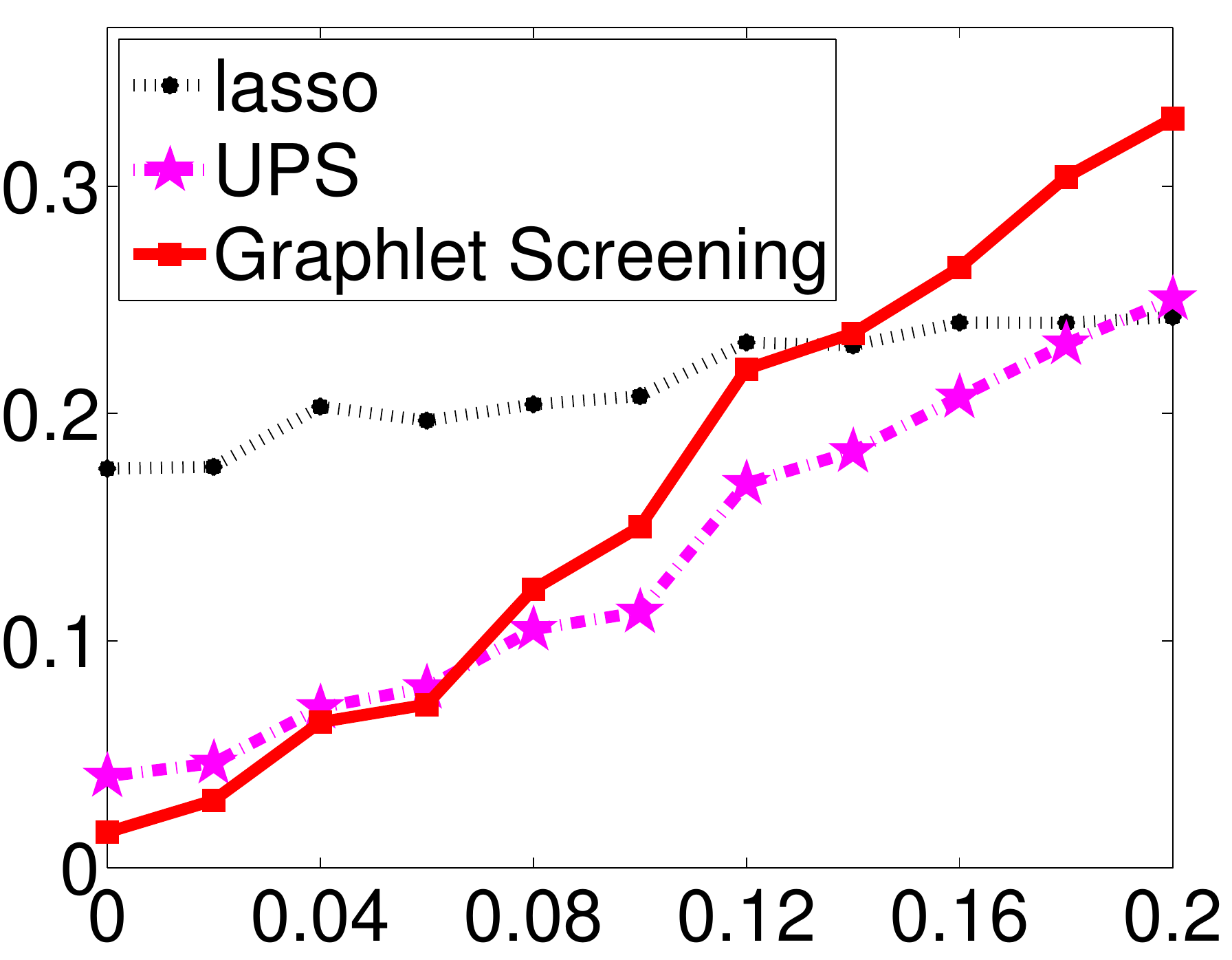}\\
  \caption{Missing predictors}\label{fig:mp}
        \end{subfigure}
        \begin{subfigure}[h]{0.33\textwidth}
                \centering
  \includegraphics[width=\textwidth]{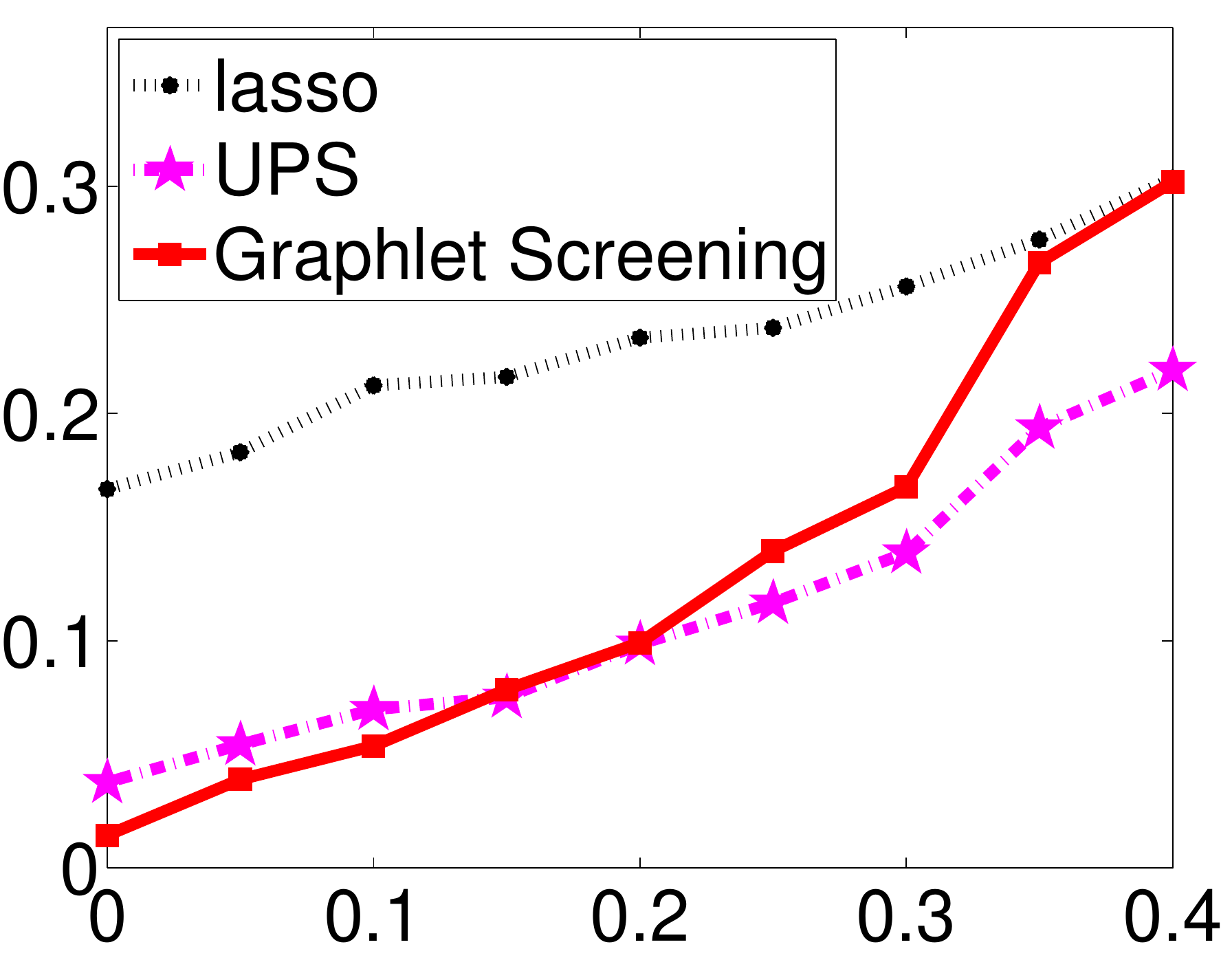}\\
  \caption{Non-linear predictors}\label{fig:nlp}
        \end{subfigure}%
                \caption{Hamming ratio results in Experiment $6$}\label{fig:e6}
\end{figure}

In Experiment $6a$, we assume the noise vector $z$ in Model (\ref{model})  is nonGaussian, where the coordinates are  iid samples from a $t$-distribution with the same degree of freedom (df)  (we assume that  $z$ is normalized so each coordinate has unit variance), where the df
range in  $\{ 3,4,5,6,7,8,9,10,30, 50\}$. Figure \ref{fig:tscaled} shows how the Hamming ratios (based on $40$ independent
repetitions) change when the $df$ decreases.   The results suggest that all three
methods are reasonably robust against nonGaussianity, but
GS continues to have the best performance.

In Experiment $6b$, we assume that the  true model is $Y = X \beta + z$ where $(X, \beta, z)$ are generated as in $4b$, but the model that is accessible to us
is a  misspecified model where the some of the true predictors are missing.
Fix $\eta \in (0,1)$, and let $S(\beta)$ be the  support of $\beta$.
For each $i \in S(\beta)$, we flip a coin that lands on head with probability $\eta$,
and we retain $i$ if and only if the coin lands on tail.
Let $S^* \subset S(\beta)$ be the set of retained indices, and let $R = S^* \cup S^c$.  The misspecified model we consider is then $Y = X^{\otimes, R} \beta^{R} + z$.

For the experiment, we let $\eta$ range in $0.02\times\{0,1,2,3,4,5,6,7,8,9,10\}$.
The average Hamming ratios (based on $40$ independent repetitions)
are reported in Figure \ref{fig:mp}. The results suggest that all three results are
reasonably robust to missing predictors, with the lasso being the most robust.
However,  as long as the proportion of true predictors that are missing is reasonably small  (say, $\eta \leq .1$),  GS continues to outperform UPS and the lasso.

In Experiment $6c$,  for $i = 1,\ldots,n$, the true model is an additive model in the form of $Y_i = \sum_{j=1}^pf_j(X_{ij})\beta_j+z_i$, but what is accessible to us is the linear model  $Y_i = \sum_{j=1}^pX_{ij}\beta_j+z_i$  (and thus misspecified; the true model is non-linear). For experiment,
we let  $(X, \beta, z)$ be generated as in $4b$, and $S(\beta)$ be the support of $\beta$.  Fixing $\eta \in (0,1)$, for each $i \in S(\beta)$, we flip a coin that lands on head with probability $\eta$, and let $S_{nl} \subset S(\beta)$ be all indices of the heads. We then randomly split $S_{nl}$ into two sets $S_1$ and $S_2$ evenly.
 For $j = 1,\ldots,p$, we define $f_j(x) = [\sgn(x)x^2\cdot1\{j\in S_{1}\}+(e^{\sqrt{n}x}-a_j)\cdot1\{ j\in S_{2}\}+x \cdot1\{j\in S_{nl}^c\}]/c_j$, where $a_j$ and $c_j$ are constants such that $\{f_j(X(i,j))\}_{i=1}^n$ has mean $0$ and variance $1/n$.

For the experiment, we let $\eta$ range in $ .05\times\{0, 1, 2,3,4,5,6,7,8\}$. The  average Hamming ratios (based on $40$ independent repetitions) are reported in Figure \ref{fig:nlp}.
The results suggest that all three methods are reasonably robust to the presence of  nonlinearity,
and GS continues to outperform UPS and the lasso when the degree of nonlinearly is moderate (say, $\eta <  .2$).

\section{Connection to existing literature and possible extensions} \label{sec:Discu}

Our idea of utilizing graph sparsity is related to the graphical lasso \citep{Meinshausen, glasso},
which also attempts to exploit
 graph structure. However, the setting we consider
here is different from that in \citet{Meinshausen,glasso},  and  our emphasis on precise optimality and
calibration is also very different.   Our method  allows nearly optimal detection
 of very rare and weak effects, because they are based on careful
analysis that has revealed a number of subtle high-dimensional
effects (e.g., phase transitions)
that we properly exploit.
Existing methodologies are not able to exploit or capture these phenomena, and can be shown to fail at the levels of rare and weak effects where we are successful.

The paper is closely related to the recent work by Ji and Jin \citep{JiJin} (see also \cite{FanLv, GJW2011}), and the two papers
use a similar rare and weak signal framework and a similar random design model.
However, they are different in important ways, since the
technical devise developed in \citet{JiJin}  can not be extended to
the current study.  For example, the lower bound derived in this paper
is different and sharper than that in \citet{JiJin}. Also,
the procedure in \citet{JiJin} relies on marginal regression  for screening.  The limitation of marginal regression is that
it neglects  the graph structure of GOSD for the regularized Gram matrix (1.5), so that it is
incapable of picking variables that have {\it  weak marginal}  correlation but {\it significant joint}  correlation to $Y$.
Correct selection of such hidden significant variables,
termed as the challenge of {\it signal cancellation} \citep{Wasserman},
is the difficulty at the heart of the variable selection problem.
One of the main innovation of GS  is that it uses  the graph structure to guide the screening,
so that it is able to successfully overcome the challenge of signal cancellation.

Additionally,  two papers have very different
objectives, and consequently the underlying analysis are very different. The main results of each of these two papers
can not be deduced from the other. For example, to assess optimality, \citet{JiJin} uses the criterion of the partition of the phase
diagram, while the current paper uses the minimax Hamming distance. Given the complexity
of the high dimensional variable selection, one type of
optimality does not imply the other, and vice versa.
Also, the main result in \citet{JiJin}
focuses on conditions under which the optimal rate of convergence is $L_p p^{1 - (\vartheta + r)^2/(4r)}$ for the {\it whole} phase space.  While this overlaps with our Corollaries  \ref{cor:entrywise1} and \ref{cor:entrywise2},
we must note that \citet{JiJin} deals with the much more difficult cases where $r / \vartheta$ can get arbitrary large; and to ensure the success in that case, they assume very strong conditions on the design matrix and the range of the signal strength.  On the other hand, the
main focus of the current paper is on optimal variable selection under conditions (of the Gram matrix $G$ as well as the signal vector $\beta$) that are as general as possible.

While the study in this paper has been focused on  the  Random Design model
RD$(\vartheta,\kappa,\Omega)$, extensions to deterministic design models
are straightforward (in fact, in Corollary \ref{cor:LB},  we have already stated some results on deterministic design models), and the omission  of discussion on the latter is largely for
technical simplicity and the sake of space.  In fact, for models with deterministic designs,
since the likelihood ratio test in the derivation of the
lower bound matches the penalized MLE in the cleaning step of GS, the optimality of GS
follows from the {\it Sure Screening and Separable After Screening} properties of GS.
The proof of these properties, and therefore the optimality of GS, follows the same line as those
for random design as long as $\max_j | \sum_i \beta_iG(i,j) I\{\Omega^{*,\delta}(i,j)=0\}|/\tau_p$ is small.
This last condition on $G$ holds when $p^{1 - \vartheta}\|G-\Omega\|_\infty =o(1)$ with a certain
$\Omega\in{\cal M}^*_p(\gamma,c_0,g,A)$.
Alternatively, this condition holds when $p^{1 - \vartheta}\|G-\Omega\|_\infty^2\log p =o(1)$
with $\Omega\in{\cal M}^*_p(\gamma,c_0,g,A)$,
provided that $\sgn(\beta_j)$ are iid symmetric random variables as in \cite{Candes}.

In this paper, we assume the signal vector $\beta$ is independent of the design matrix $X$, and that $\beta$ is modeled by a Bernoulli model through $\beta = b \circ \mu$. Both assumptions
can be relaxed. In fact,  in order for GS  to work, what we really need is some {\it decomposability} condition similar to that  in  Lemma \ref{lemma:sparsegraph}, where except for negligible probabilities,  the maximum size of the graphlets $m_0^*= m_0^*(S(\beta), G, \delta)$ is small.
In many situations, we can show that $m_0^*$ does not exceed a fixed integer. One of such examples is as follows.
Suppose for any fixed integer $m \geq 1$ and size-$m$ subset $S$ of $\{1, 2, \ldots, p\}$,  there are constants $C > 0$ and $d > 0$  such that the conditional probability $P(\beta_j \neq 0, \forall j \in S |X) \leq  C p^{- d m}$.
In fact, when such a condition holds, the claim follows since $\cg^{*, \delta}$  has no more than $C (e K)^m$ size-$m$ connected subgraphs if it is $K$-sparse.
See the proof of Lemma \ref{lemma:sparsegraph} for details.
Note that when $\eps_p = p^{-\vartheta}$ as in the ARW,  then the condition holds for the Bernoulli model in Lemma \ref{lemma:sparsegraph},   with $d  = \vartheta$. Note also that the Bernoulli model can be replaced by some Ising models.

Another interesting direction of future research is the extension of GS  to
more general models such as logistic regression. The extension of the lower bound
in Theorem \ref{thm:LB} is relatively simple since the degree of GOLF can be bounded
using the true $\beta$. This indicates the optimality of GS  in logistic and other
generalized linear models as long as proper generalized likelihood ratio or Bayes tests
are used in both the $GS$- and $GC$-steps.

\section{Proofs}  \label{sec:appen}
In this section, we provide all technical proofs. We assume $\sigma = 1$ for simplicity.
\subsection{Proof of Lemma \ref{lemma:sparsegraph}}
When $\cg_S^{*,\delta}$ contains a connected subgraph of size $\geq m_0+1$, it must
contain a connected subgraph with size $m_0+1$. By \cite{Frieze}, there are
$\leq p (e K)^{m_0+1}$ connected subgraph of size $m_0 +1$.
Therefore, the probability that $\cg_S^{*,\delta}$ has a connected subgraph of size $(m_0 + 1)$
$\leq p (e K)^{m_0+1} \eps_p^{m_0+1}$. Combining these gives the claim. \qed

\subsection{Proof of Theorem \ref{thm:LB}}   \label{subsec:proofLB}
Write for short $\rho_j^* = \rho_j^*(\vartheta, r, a, \Omega)$.
Without loss of generality, assume $\rho_1^* \leq \rho_2^* \leq  \ldots \leq  \rho_p^*$.
We construct indices $i_1 < i_2 < \ldots < i_m$ as follows.  (a)  start with  $B = \{1, 2, \ldots, p \}$ and let $i_1 = 1$, (b) updating $B$ by removing $i_1$ and all nodes $j$  that are   neighbors of  $i_1$ in GOLF,  let $i_2$ be the smallest index, (c) defining $i_3, i_4, \ldots, i_m$ by  repeating (b), and terminates the process when no indices is left in $B$.  Since each time
we remove at most $d_p(\cg^{\diamond} )$ nodes, it follows that
\begin{equation} \label{lemmagraph2pf1}
\sum_{j = 1}^p p^{-\rho_j^*} \leq  d_p(\cg^{\diamond} )  \sum_{k = 1}^m  p^{- \rho_{i_k}^*}.
\end{equation}

For each $1 \leq j \leq p$, as before, let $(V_{0j}^*,  V_{1j}^*)$ be the least favorable configuration, and let
$(\theta_{*j}^{(0)}, \theta_{*j}^{(1)}) = \margmin_{ \{ \theta^{(0)} \in B_{V_{0j}^*},
\theta^{(1)} \in B_{V_{1j}^*},\sgn(\theta^{(0)})\neq \sgn(\theta^{(1)})\}}  \alpha(\theta^{(0)}, \theta^{(1)}; \Omega)$.  By our notations, it is seen that
\begin{equation} \label{testprob0}
\rho_j^* = \eta(V_{0j}^*, V_{1j}^*; \Omega), \qquad \alpha^*(V_{0j}^*, V_{1j}^*;\Omega) = \alpha(\theta_{*j}^{(0)}, \theta_{*j}^{(1)}; \Omega).
\end{equation}
In case $(\theta_{*j}^{(0)}, \theta_{*j}^{(1)})$ is not unique, pick one arbitrarily. We construct a $p \times 1$ vector $\mu^*$ as follows.
Fix $j \in \{i_1,\cdots,i_m\}$.   For all indices in $V_{0j}^*$, set the constraint of $\mu^*$ on these indices to be $\theta_{*j}^{(0)}$.  For any index $i \notin \cup_{k=1}^m V_{0 i_k}^* $, set $\mu^*_i = \tau_p$. Since
\begin{equation} \label{lemmagraph2pf2a}
\hamm_p^*(\vartheta,\kappa, r,a, \Omega) \geq \inf_{\hb} H_p(\hb;\eps_p,n_p, \mu^*, \Omega) = \inf_{\hb} \sum_{i = 1}^p P(\sgn(\hb_j) \neq \sgn(\beta_j)),
\end{equation}
it follows that
\begin{equation} \label{lemmagraph2pf2b}
\hamm_p^*(\vartheta,\kappa, r,a, \Omega) \geq   \sum_{k = 1}^m  \sum_{j \in V_{0 i_k}\cup V_{1 i_k}}    P(\sgn(\hb_j) \neq \sgn(\beta_j)),
\end{equation}
where $\beta  = b \circ \mu^*$ in (\ref{lemmagraph2pf2a})-(\ref{lemmagraph2pf2b}). Combining (\ref{lemmagraph2pf1}) and (\ref{lemmagraph2pf2b}), to show the claim, we only need to show that for any $1 \leq k \leq m$ and any procedure $\hb$,
\begin{equation} \label{lemmagraph2pf9}
\sum_{j \in V_{0 i_k}\cup V_{1 i_k}}    P(\sgn(\hb_j) \neq \sgn(\beta_j))  \geq  L_p p^{-\rho_{i_k}^*}.
\end{equation}

Towards this end, we write for short $V_0 = V_{0 i_k}$, $V_1 = V_{1 i_k}$, $V = V_0 \cup V_1$,
$\theta^{(0)} = \theta_{*i_k}^{(0)}$, and $\theta^{(1)} = \theta_{*i_k}^{(1)}$.  Note that by Lemma \ref{lemma:V},
\begin{equation} \label{eqv}
|V| \leq (\vartheta + r)^2/(2 \vartheta r).
\end{equation}
Consider a test setting where under the null  $H_0$,    $\beta =  \beta^{(0)} = b \circ \mu^*$ and $I_V \circ \beta^{(0)} = I_V \circ \theta^{(0)}$,   and under the alternative $H_1$,  $\beta = \beta^{(1)}$ which is constructed by keeping  all coordinates of $\beta^{(0)}$  unchanged,  except those coordinates in $V$ are perturbed in a way so that  $I_V \circ \beta^{(1)} = I_{V} \circ\theta^{(1)}$.   In this construction, both $\beta^{(0)}$ and $\beta^{(1)}$ are assumed as known, but we don't know which of $H_0$ and $H_1$ is true.
In the literature, it is known that $\inf_{\hat{\beta}} \sum_{j \in  V}    P(\sgn(\hb_j) \neq \sgn(\beta_j))$ is not smaller than the minimum  sum of Type I and Type II errors associated with this testing problem.

Note that  by our construction and (\ref{testprob0}),  the right hand side is $\alpha^*(V_0,V_1;\Omega)$.
At the same time, it is seen the optimal test statistic is $Z\equiv(\theta^{(1)} - \theta^{(0)})' X'  (Y - X \beta^{(0)})$. It is seen that up to some negligible terms,  $Z \sim  N(0, \alpha^*(V_0,V_1;\Omega) \tau_p^2)$ under $H_0$, and
$Z \sim   N(\alpha^*(V_0,V_1;\Omega) \tau_p^2, \alpha^*(V_0,V_1;\Omega) \tau_p^2)$ under $H_1$.
 The optimal test is to reject $H_0$ when $Z \geq t[\alpha^*(V_0,V_1;\Omega)]^{1/2}  \tau_p$ for some threshold $t$, and the minimum sum of Type I and  Type II error is
\begin{equation} \label{testprob4}
\inf_t \bigl\{ \eps_p^{|V_0|} \bphi(t)  + \eps_p^{|V_1|} \Phi(t -  [\alpha^*(V_0,V_1;\Omega)]^{1/2}  \tau_p) \bigr\}.
\end{equation}
Here, we have used  $P(H_0) \sim  \eps_p^{|V_0|}$ and $P(H_1) \sim
\eps_p^{|V_1|}$, as a result of the Binomial structure in $\beta$.
It follows that
$\sum_{j \in V}    P(\sgn(\hb_j) \neq \sgn(\beta_j))  \gtrsim  \inf_t \bigl\{ \eps_p^{|V_0|} \bphi(t)  + \eps_p^{|V_1|} \Phi(t -  [\alpha^*(V_0,V_1;\Omega)]^{1/2}  \tau_p) \bigr\}$.
Using Mills' ratio and definitions,  the right hand side $\geq L_p p^{-\eta(V_0, V_1;\Omega)}$, and (\ref{lemmagraph2pf9})  follows by recalling (\ref{testprob0}). \qed


\subsection{Proof of Corollaries  \ref{cor:entrywise1}, \ref{cor:entrywise2}, and \ref{cor:main}}  \label{subsec:entrywise1}
When $a > a_g^*(\Omega)$,  $\rho_j^*(\vartheta, r, a, \Omega)$
does not depend on $a$, and have an alternative expression as follows.
For any subsets $D$ and $F$ of $\{1, 2, \ldots, p\}$, let $\omega(D, F; \Omega)$ be as in (\ref{Defineomega}).
Introduce  $\rho(D, F; \Omega) = \rho(D, F; \vartheta, r, a,  \Omega, p)$ by
\begin{equation} \label{Definerho}
\rho(D,F; \Omega)  = \frac{(|D|+2|F|)\vartheta}{2} +
\left\{
\begin{array}{ll}
\frac{1}{4}\omega(D,F; \Omega) r,  &      |D| \text{ is even},  \\
\frac{\vartheta}{2}  + \frac{1}{4}\bigl[ (\sqrt{\omega(D,F; \Omega)r}- \frac{\vartheta}{\sqrt{\omega(D,F; \Omega) r} } )_+\bigr]^2, &  |D| \text{ is odd}.
\end{array}\right.
\end{equation}
The following lemma is proved in Section \ref{subsubsec:omega}.
\begin{lemma} \label{lemma:omega}
Fix $m_0 \geq 1$,    $(\vartheta, \kappa) \in (0,1)^2$,  $r > 0$,  $c_0 > 0$, and $g > 0$ such that
$\kappa > (1 - \vartheta)$.
Suppose the conditions of Theorem \ref{thm:LB} hold, and that  for sufficiently large $p$,
(\ref{Conditiona}) is satisfied.
Then as $p \goto \infty$,
$\rho_j^*(\vartheta, r, a, \Omega)$ does not depend on $a$, and satisfies
$\rho_j^*(\vartheta, r, a, \Omega)  = \min_{\{ (D, F):
{ j \in D \cup F},  D \cap F = \emptyset, D \neq \emptyset, |D\cup F|\leq g \}}  \rho(D,F; \Omega)$.
\end{lemma}

We now show Corollaries \ref{cor:entrywise1}-\ref{cor:main}.  Write for short $\omega = \rho(D,F; \Omega)$, $T = r / \vartheta$,  and $\lambda_k^* = \lambda_k^*(\Omega)$.  The following inequality is frequently used below, the proof of which is elementary so we omit it:
\begin{equation} \label{lemmaentry0}
\omega\geq \lambda_k^* |D|, \qquad  \mbox{where $k = |D| + |F|$}.
\end{equation}
To show these corollaries,
it is sufficient to show for all subsets  $D$ and $F$ of $\{1, 2, \ldots, p\}$,
\begin{equation} \label{lemmaentry1a}
\rho(D,F; \Omega) \geq (\vartheta + r)^2/(4r),  \qquad |D| \geq 1,
\end{equation}
where $\rho(D,F; \Omega)$ is as in (\ref{Definerho}). By basic algebra,   (\ref{lemmaentry1a}) is equivalent to
\begin{equation} \label{lemmaentry1b}
\left\{
\begin{array}{ll}
(\omega T + 1/(\omega T) -2) 1 \{ \omega T \geq 1\} \geq (T + 1/T - 2 (|D| + 2|F|)),  &\qquad \mbox{$|D|$ is odd},  \\
\omega \geq  \frac{2}{T}[(T + 1/T)/2 + 1 - (|D| + 2 |F|)], &\qquad \mbox{$|D|$ is even}.
\end{array}
\right.
\end{equation}
Note that when  $(|D|, |F|) = (1, 0)$, this claim holds trivially,  so it is sufficient to consider the case where
\begin{equation} \label{lemmaentry1c}
|D| + |F| \geq 2.
\end{equation}
We now show that (\ref{lemmaentry1b})  holds under the conditions of each of
corollaries.

\subsubsection{Proof of  Corollary \ref{cor:entrywise1}}
In this corollary,  $1 < (T + 1/T)/2 \leq 3$, and if
either (a) $|D| + 2 |F|  \geq 3$ and $|D|$ is odd  or (b) $|D| + 2 |F| \geq 4$ and $|D|$ is even,
the right hand side of (\ref{lemmaentry1b}) $\leq 0$, so the claim holds trivially. Therefore,
all we need to show is the case where $(|D|, |F|) = (2, 0)$.   In this case,
since each off-diagonal coordinate $\leq  4 \sqrt{2} - 5 \equiv \rho_0$, it follows from definitions  and basic algebra that $\omega \geq 2 (1 - \rho_0) = 4(3 - 2 \sqrt{2})$, and (\ref{lemmaentry1b}) follows by noting that $\frac{2}{T}[(T + 1/T)/2 + 1 - (|D| + 2 |F|)] = (1 - 1/T)^2 \leq 4(3 - 2 \sqrt{2})$.   \qed

\subsubsection{Proof of  Corollary \ref{cor:entrywise2}}
In this corollary,  $1 < (T + 1/T)/2 \leq 5$.
First, we consider the case where $|D|$ is odd. By similar argument,  (\ref{lemmaentry1b}) holds trivially when $|D| + 2 |F| \geq 5$, so all we need to consider is  the case $(|D|, |F|)  = (1,1)$ and the case $(|D|, |F|) = (3, 0)$. In both cases, $|D| + 2 |F| = 3$.   By \eqref{lemmaentry0}, when $\omega T<1$, there must be $T<1/\min(\lambda_2^*,3\lambda_3^*)$. By the conditions of this corollary, it follows $T<(5+2\sqrt{6})/4<3+2\sqrt{2}$. When $1<T<3+2\sqrt{2}$, there is $T + 1/T - 6<0$, and thus \eqref{lemmaentry1b} holds for $\omega T<1$. When $\omega T \geq 1$, (\ref{lemmaentry1b}) holds if and only if $\omega T  + \frac{1}{\omega T} - 2 \geq T + 1/T - 6$.
By basic algebra, this holds if
\begin{equation} \label{entry2a}
\omega  \geq \frac{1}{4}\bigl[(1 - 1/T) + \sqrt{(1 - 1/T)^2 - 4/T}\bigr]^2.
\end{equation}
Note that the right hand side of (\ref{entry2a}) is a monotone in  $T$ and has a  maximum of  $(3 + 2 \sqrt{2}) (5 - 2\sqrt{6})$ at $T = (5 + 2 \sqrt{6})$. Now, on one other hand,
when $(|D|, |F|) = (1, 0)$,  by (\ref{lemmaentry0}) and conditions of the corollary, $\omega \geq 3 \lambda_3^* > (3 + 2 \sqrt{2}) (5 - 2 \sqrt{6})$. On the other hand, when $(|D|, |F|) = (1,1)$, by  basic algebra and that each off-diagonal coordinate of $\Omega  \leq \sqrt{1 + (\sqrt{6} - \sqrt{2})}/(1 + \sqrt{3/2}) \equiv \rho_1$ in magnitude,
$\omega \geq 1 - \rho_1^2 =  (3 + 2 \sqrt{2})(5 - 2 \sqrt{6})$.  Combining these gives (\ref{lemmaentry1b}).

We now consider the case where $|D|$ is even. By similar argument, (\ref{lemmaentry1b}) holds when
$|D| + 2 |F|  \geq  6$, so all we need is to show is that (\ref{lemmaentry1b}) holds for the following three cases: $(|D|, |F|) = (4, 0), (2, 1),(2, 0)$. Equivalently, this is to show that  $\omega \geq \frac{2}{T} [(T + 1/T)/2 - 3]$ in the first two cases and  that $\omega \geq \frac{2}{T} [(T + 1/T)/2 - 1]$ in the last case. Similarly, by the monotonicity of the right hand side of these  inequalities, all we need to show is $\omega \geq 4(5 - 2 \sqrt{6})$ in the first two cases, and
$\omega \geq 8 (5 - 2 \sqrt{6})$ in the last case.
Now, on one hand, using (\ref{lemmaentry0}),  $\omega \geq 4 \lambda_4^*$ in the first case, and $\omega \geq 2 \lambda_3^*$ in the second case, so by the conditions of the corollary, $\omega \geq 4(5 - 2 \sqrt{6})$ in the first two cases. On  the other hand,
in the last case, since all  off-diagonal coordinates of $\Omega \leq
8 \sqrt{6} - 19 \equiv \rho_0$ in magnitude,   and $\omega \geq 2 (1 - \rho_0) = 8 (5 - 2 \sqrt{6})$. Combining these gives (\ref{lemmaentry1b}). \qed

\subsubsection{Proof of Corollary \ref{cor:main}}
Let $N$ be the unique integer such that $2N-1 \leq (T + 1/T)/2 < 2N + 1$. First, we consider the case where $|D|$ is odd. Note that when $|D| + 2 |F| \geq 2N + 1$,  the right hand side of (\ref{lemmaentry1b}) $\leq 0$, so all we need to consider is the case $|D| + 2 |F| \leq 2N-1$.
Write for short $k = k(D, F) = |D| + |F|$ and
$j = j (D, F) = (|D| + 2|F| + 1)/2$.  By (\ref{lemmaentry1c}), definitions, and that $|D| + 2 |F| \leq 2 N -1$,  it is seen that $2 \leq k \leq 2 N - 1$ and  $(k+1)/2 \leq j \leq \min\{k, N \}$. By the condition of the corollary,
$\lambda_k^* \geq        \frac{(T + 1/T)/2 - 2 j + 2+ \sqrt{[(T + 1/T)/2 - 2j +2]^2 - 1}} {T (2k - 2j  + 1) }$.
Note that $|D| = 2k - 2j + 1$.  Combining these with (\ref{lemmaentry0})  gives
$\omega  T  \geq (2k - 2j +1) \lambda_k^* T  \geq (T + 1/T)/2 - 2j +  2+ \sqrt{[(T + 1/T)/2 - 2j + 2]^2 -1}  \geq 1$.
and (\ref{lemmaentry1b}) follows by basic algebra.

We now consider  the case where $|D|$ is even. Similarly, the right hand side of (\ref{lemmaentry1b}) is negative when  $|D| + 2 |F| \geq 2(N +1)$,  so we only need to consider the case where $|D| + 2 |F| \leq 2 N$.  Similarly, write for short $k = k(D, F)  = |D| + |F|$ and
$j = (|D| + 2 |F|)/2$. It is seen that $2 \leq k \leq 2N$ and $k/2 \leq j \leq \min\{k-1, N \}$.
By the conditions of the corollary,
$\lambda_k^* \geq \frac{(T + 1/T)/2 + 1 - 2 j}{T(k-j)}$.
Note that $|D| = k-j$. It follows from (\ref{lemmaentry0}) that
$\omega \geq 2(k - j) \lambda_k^* \geq \frac{2}{T} [(T + 1/T)/2 + 1 - 2 j]$, and
(\ref{lemmaentry1b}) follows. \qed

\subsubsection{Proof of Lemma \ref{lemma:omega}}
\label{subsubsec:omega}
Let sets $V_0$ and $V_1$ and vectors $\theta^{(0)}$ and $\theta^{(1)}$ be as  in Section \ref{subsec:LB}, and let $V=V_0\cup V_1$.
By  the definition of $\rho_j^*(\vartheta, r, a, \Omega)$,
$\rho_j^*(\vartheta, r, a, \Omega)=\min(I,II)$,
where
$I=\min_{\{ (V_0, V_1):  j\in V_1\cup V_0,V_0\neq V_1 \}}\eta(V_0, V_1;\Omega)$ and $II=\min_{\{ V_0:  j \in V_0 \cup  V_1, V_0 = V_1 \}}   \eta(V_0, V_1;\Omega)$.
So to show the claim, it is sufficient to show
\begin{equation} \label{proveIII}
I = \min_{\{ (D, F): j \in D \cup F, D \cap F = \emptyset, D \neq \emptyset, |D\cup F|\leq g \}} \rho(D,F; \Omega), \qquad II \geq I.
\end{equation}

Consider the first claim in (\ref{proveIII}).  Write for short
$F= F(V_0, V_1) =  V_0\cap V_1$ and $D= D(V_0, V_1) = V \setminus F$.  By the definitions,  $D\neq\emptyset$. The key is to show that when $|V_0 \cup V_1| \leq g$,
\begin{equation} \label{lemmaomega1}
\alpha^*(V_0,V_1;\Omega)   =    \omega(D,  F; \Omega).
\end{equation}
Towards this end, note that by definitions,  $\alpha^*(V_0,V_1;\Omega)  =  \alpha(\theta_*^{(0)}, \theta_*^{(1)})$, where
$(\theta_*^{(0)}, \theta_*^{(1)}) =    \margmin_{\{ \theta^{(0)} \in B_{V_0}, \theta^{(1)} \in B_{V_1}
\}}     \alpha(\theta^{(0)}, \theta^{(1)})$.
By $a > a_g^*(\Omega)$ and the way $a_g^*(\Omega)$ is defined,  $(\theta_*^{(0)}, \theta_*^{(1)})$  remains
as the solution of  the optimization problem if we relax the conditions $\theta^{(i)} \in B_{V_i}$ to that of $\theta^{(i)} =  I_{V_i} \circ \mu^{(i)}$, where $\mu^{(i)}  \in  \Theta_p(\tau_p)$ (so that upper bounds on the signal strengths are removed), $i = 0, 1$.  As a result,
\begin{equation} \label{lemmaomega2}
\alpha^*(V_0,V_1;\Omega)  =  \min_{\{ \theta^{(i)} \in I_{V_i} \circ \mu^{(i)},  \mu^{(i)} \in \Theta_p(\tau_p),  i = 0, 1,\}} \alpha(\theta^{(0)}, \theta^{(1)}).
\end{equation}
We now study (\ref{lemmaomega2}). For short, write $\xi = \tau_p^{-1} (\theta^{(1)} - \theta^{(0)})^{V}$, $\Omega_{VV}  = \Omega^{V, V}$,  $\xi_D  = \tau_p^{-1} (\theta^{(1)} - \theta^{(0)})^D$,  and similarly for $\Omega_{DD}$,  $\Omega_{DF}$, $\Omega_{FD}$, $\Omega_{FF}$, and $\xi_F$.
Without loss of generality, assume the indices in $D$ come first in $V$. It follows
\[
\Omega_{VV}  = \begin{pmatrix}
\Omega_{DD} & \Omega_{DF}\\
\Omega_{FD} & \Omega_{FF}
\end{pmatrix},
\]
and
\begin{equation} \label{lemmaomega3}
\alpha(\theta^{(0)}, \theta^{(1)}) = \xi' \Omega_{VV}
\xi  = \xi_D' \Omega_{DD} \xi_D + 2 \xi_D'\Omega_{DF} \xi_F + \xi_F'\Omega_{FF} \xi_F.
\end{equation}
By definitions, it is seen that there is no  constraint on the coordinates of $\xi_F$, so to optimize the quadratic form in (\ref{lemmaomega1}), we need to choose $\xi$ is a way such that
$\xi_F = -\Omega_{FF}^{-1}\Omega_{FD} \xi_D$, and that  $\xi_D$ minimizes
$\xi_D'  ( \Omega_{DD} - \Omega_{DF} \Omega_{FF}^{-1} \Omega_{FD})  \xi_D$,
where every coordinate  of $\xi_D \geq 1$ in magnitude.
Combining these with  (\ref{lemmaomega2}) gives (\ref{lemmaomega1}).

At the same time, we rewrite
\begin{equation} \label{lemmaomega4}
I=\min_{\{ (D, F): j \in D \cup F, D \neq \emptyset, D \cap F = \emptyset\}}    \biggl\{  \min_{\{(V_0, V_1): V_0\cup V_1=D\cup F,   V_0\cap V_1=F \} }   \eta(V_0, V_1;\Omega) \biggr\}.
\end{equation}
By similar arguments as in the proof of Lemma \ref{lemma:V}, the subsets $(V_0, V_1)$ that achieve the minimum of $\eta(V_0, V_1;\Omega)$ must satisfy $|V_0 \cup V_1| \leq g$.
Using (\ref{lemmaomega1}), for any fixed $D$ and $F$ such that $|D \cup F| \leq g $,  $D \neq \emptyset$ and $D \cap F = \emptyset$,   the term in the big bracket on the right hand side is
$
\min_{\{(V_0, V_1):  V_0\cup V_1=D\cup F,  V_0\cap V_1=F\}}  \{   \frac{(2|F|+|D|)\vartheta}{2} +  \frac{\bigl||V_1| - |V_0| \bigr|\vartheta}{2}  + \frac{1}{4} [(\sqrt{\omega(D,F; \Omega)r} - \frac{\bigl||V_1| - |V_0| \bigr|\vartheta}{\sqrt{\omega(D,F; \Omega)r}})_+]^2\}$.
It is worth noting that for fixed $D$ and $F$, the above quantity is monotone increasing with $\bigl||V_1| - |V_0| \bigr|$.  When $|D|$ is even, the minimum is achieved at $(V_0, V_1)$ with  $|V_0| = |V_1|$, and when $|D|$ is odd, the minimum is achieved at $(V_0, V_1)$ with
$\bigl| |V_1|  -  |V_0| \bigr|  = 1$, and in both cases, the minimum is $\rho(D,F; \Omega)$. Inserting this to (\ref{lemmaomega4}), it is seen that
\begin{equation} \label{lemmaomega5}
I=\min_{\{ (D, F): j \in D \cup F, D \cap F = \emptyset, D \neq \emptyset,  |D \cup F| \leq g \}} \rho(D,F; \Omega),
\end{equation}
which is the first claim in (\ref{proveIII}).

Consider the second claim of (\ref{proveIII}).  In this case, by definitions,
$V_0=V_1$ but $\sgn(\theta^{(0)})\neq\sgn(\theta^{(1)})$. Redefine $D$ as  the subset of $V_0$ where the signs of the coordinates of $\theta^{(0)}$ do not equal to those of $\theta^{(1)}$, and let $F =V\setminus  D$. By definitions, it is seen that  $\alpha^*(V_0,V_0;\Omega)=4\alpha^*(F,V_0;\Omega)$, where we note $D \neq \emptyset$ and $F \neq V_0$.   By the definition of $\eta(V_0, V_1;\Omega)$, it follows that $\eta(V_0, V_0;\Omega)\geq \eta(F,V_0;\Omega)$, and the claim follows.   \qed


\subsection{Proof of Lemma \ref{lemma:SS}}
Write for short $\rho_j^* = \rho_j^*(\vartheta, a,  r, \Omega)$.  To show the claim, it is sufficient to show that for any
fixed $1 \leq j \leq p$,
\begin{equation} \label{lemmaSS1}
P(j\notin{\cal U}_p^*, \beta_j \neq 0) \leq L_p [ p^{-\rho_j^*} + p^{- (m_0+1) \vartheta} + o(1/p)].
\end{equation}
Using Lemma \ref{lemma:sparseGOSD} and  \cite[Lemma 3.1]{JiJin},
there is an event $A_p$ that depends on $(X, \beta)$ such that $P(A_p^c) \leq o(1/p)$ and that
over the event, $\Omega^{*,\delta}$ is $K$-sparse with $K = C (\log(p))^{1/\gamma}$,
$\|\Omega^{*,\delta} - \Omega\|_\infty \leq (\log(p))^{-(1 - \gamma)}$, $\|(X'X - \Omega) \beta \| _{\infty} \leq C \|\Omega\| \sqrt{2 \log(p)} p^{-[(\kappa - (1 - \vartheta)]/2}$,
and for all subset $B$ with size $\leq m_0$,
$\|G^{B, B} - \Omega^{B, B} \|_\infty \leq L_p p^{-\kappa/2}$.
Recall that $\cg^{*,\delta}$ is the GOSD and $\cg_S^{*,\delta}$ is the subgraph of the GOSD formed by the nodes in the support of $\beta$, $S(\beta) = \{1 \leq j \leq p: \beta_j \neq 0\}$.
When $\beta_j  \neq 0$, there is a unique component $\call_0$ such that $j \in \call_0 \lhd \cg_S^{*,\delta}$ ($A \lhd B$ means that $A$ is component  or maximal connected subgraph  of $B$).  Let $B_p$ be the event $|\call_0| \leq m_0$.  By Frieze \cite{Frieze}, it is seen that
$P(B_p^c \cap A_p) \leq L_p p^{-(m_0+1) \vartheta}$.
So to show (\ref{lemmaSS1}), it is sufficient to show that
\begin{equation} \label{lemmaSS2}
P(j\notin{\cal U}_p^*, j \in \call_0 \lhd \cg_S^{*,\delta},   A_p \cap B_p)  \leq  L_p p^{-\rho_j^*}.
\end{equation}

Now, in the screening procedure, when we screen $\call_0$, we have $\call_0 = \hat{D} \cup \hat{F}$ as in (\ref{DefineDF}).  Since
the event $\{j\notin{\cal U}_p^*, j \in \call_0 \lhd \cg_S^{*,\delta}  \}$ is contained in the event
$\{ T(Y, \hd, \hat{F}) < t(\hd, \hat{F})\}$,
$P(j\notin{\cal U}_p^*, j \in \call_0 \lhd \cg_S^{*,\delta},  A_p \cap B_p)  \leq  P(T(Y, \hat{D}, \hat{F}) \leq t(\hat{D}, \hat{F}), j \in \call_0 \lhd \cg_S^{*,\delta},   A_p \cap B_p)$,
where the right hand side does not exceed
\begin{equation} \label{lemmaSSadd}
\sum_{(\call_0,  D, F): \mbox{$j  \in \call_0$ \& $\call_0 = D \cup F$ is a partition} }   P(T(Y, D, F) \leq t(D, F), j \in \call_0 \lhd \cg_S^{*,\delta},  A_p \cap B_p);
\end{equation}
note that $(\call_0, D, F)$ do not depend on $z$ (but may still depend on $(X, \beta)$).
First, note that over the event $A_p$,  there are at most $(e K)^{m_0+1}$ $\call_0$ such that $j \in \call_0$ and $|\call_0| \leq m_0$. Second, note that
 for each $\call_0$, there are only finite ways to partition it to $D$ and $F$.   Last,   note that for any fixed $j$ and $\call_0$,
$P(j \in \call_0 \lhd \cg_S^{*,\delta}) \leq  \eps_p^{|\call_0|}$.   Combining these observations,  to show (\ref{lemmaSS2}), it is sufficient to show that for  any such triplet $(\call_0, D, F)$,
\begin{equation} \label{lemmaSS3}
\eps_p^{|\call_0|} P \bigl(T(Y, D, F) \leq t(D, F)  \bigr|    \{j \in \call_0 \lhd \cg_S^{*,\delta} \}  \cap A_p \cap B_p \bigr)  \leq L_p p^{-\rho_j^*}.
\end{equation}

We now show (\ref{lemmaSS3}).  Since $\lambda^*_{m_0}(\Omega) \geq C > 0$, it follows from the definition of $A_p$ and basic algebra  that  for any realization of $(X, \beta)$ in $A_p \cap B_p$,
\begin{equation} \label{lemmaSS7a}
\|(G^{\call_0, \call_0})^{-1} \|_{\infty} \leq C.
\end{equation}
Recall that $\ty = X'Y$ and denote for short $y = (G^{\call_0, \call_0})^{-1} \ty^{\call_0}$.
It is seen that
\begin{equation} \label{lemmaSS9a}
y =   \beta^{\call_0} + w + rem, \qquad w \sim N(0, (G^{\call_0, \call_0})^{-1}), \qquad rem \equiv  (G^{\call_0, \call_0})^{-1}   G^{\call_0, \call_0^c} \beta^{\call_0^c}.
\end{equation}
Since $\call_0$ is a component of $\cg_S^{*,\delta}$,   $(\Omega^{*,\delta})^{\call_0, \call_0^c} \beta^{\call_0^c} = 0$. Therefore, we can write $rem = (G^{\call_0, \call_0})^{-1} (I + II)$, where $I = (G^{\call_0, \call_0^c} - \Omega^{\call_0, \call_0^c}) \beta^{\call_0^c}$ and
$II = [\Omega^{\call_0, \call_0^c} - (\Omega^{*,\delta})^{\call_0, \call_0^c}]  \beta^{\call_0^c}$.
By the definition of $A_p$,
$\| I\|_\infty \leq C \sqrt{2 \log(p)} p^{-[\kappa - (1 - \vartheta)]/2}$, and $\| II \|_\infty \leq   \| \Omega - \Omega^{*,\delta} \|_{\infty} \| \beta^{\call_0^c}\|_{\infty} \leq C \tau_p (\log(p))^{-(1 - \gamma)}$. Combining these  with (\ref{lemmaSS7a}) gives $\|rem\|_{\infty} \leq C \tau_p (\log(p))^{-(1 - \gamma)}$.

At the same time,  let $y_1$, $w_1$, and $rem^1$ be the restriction of $y$,  $w$, and $rem$  to indices in $D$, correspondingly,  and let $H = [G^{D, D} -  G^{D, F} ( G^{F, F})^{-1} G^{F, D}]$.
By (\ref{lemmaSS9a}) and  direct calculations,
$T(Y,   D, F) =  y_1'  H  y_1$, $y_1 \sim N(\beta^D + rem^1, H^{-1})$, and so $T(Y, D, F)$ is distributed as $\chi_{|D|}^2(\delta)$,
where the non-central parameter is
$(\beta^{D} + rem^1)' H (\beta^D + rem^1) = \delta   + O( (\log(p))^{ \gamma})$ and $ \delta \equiv  (\beta^D)' H \beta^D$.
Since $\lambda_{m_0}^*(\Omega) \geq C$,   $\delta \geq C \tau_p^2$ and is the dominating term. It follows that
\begin{equation} \label{lemmaSS11a}
P(T(Y, D, F) \leq t(D, F)  \bigr|  \{j \in \call_0 \lhd \cg_S^{*,\delta}\}  \cap A_p \cap B_p \bigr)  \lesssim P\bigl( \chi_{|D|}^2(\delta)  \leq t(D, F) \bigr).
\end{equation}
Now, first, by definitions, $\delta \geq 2 \omega(D,F; \Omega) r \log(p)$, so by basic knowledge on non-central $\chi^2$,
\begin{equation} \label{lemmaSS11b}
P(\chi_{|D|}^2(\delta) \leq t(D, F) ) \leq P(\chi_{|D|}^2(2 \omega(D,F; \Omega) r \log(p)) \leq t(D, F)).
\end{equation}
Second, recalling $t(D, F) = 2 q \log(p)$, we have
\begin{equation} \label{lemmaSS11c}
P( \chi_{|D|}^2(2 \omega(D,F; \Omega) r \log(p)) \leq t(D, F)) \leq L_p p^{-[(\sqrt{\omega(D,F; \Omega) r}-\sqrt{q})_+]^2}.
\end{equation}
Inserting (\ref{lemmaSS11b})-(\ref{lemmaSS11c}) into (\ref{lemmaSS11a}) and recalling $\eps_p = p^{-\vartheta}$,
\begin{equation} \label{lemmaSS12}
\eps_p^{|\call_0|} P(T(Y, D, F) \leq t(D, F)  \bigr|  \{j \in \call_0 \lhd \cg_S^{*,\delta} \}  \cap A_p \cap B_p \bigr)  \leq L_p  p^{-(|\call_0|\vartheta + [(\sqrt{\omega(D,F; \Omega) r}-\sqrt{q})_+]^2)}.
\end{equation}
By the choice of $q$ and direct calculations,
\begin{equation} \label{lemmaSS13}
|\call_0|  \vartheta + [(\sqrt{\omega(D,F; \Omega) r}-\sqrt{q})_+]^2 \geq \rho(D,F; \Omega) \geq \rho_j^*,
\end{equation}
where $\rho(D,F; \Omega)$ as in \eqref{Definerho}.
Combining (\ref{lemmaSS12})-(\ref{lemmaSS13}) gives (\ref{lemmaSS3}).    \qed

\subsection{Proof of Lemma \ref{lemma:SAS}}
In the screening stage, suppose we pick the threshold $t(\hd, \he) = 2 q \log(p)$ in a way such that there is a constant $q_0(\vartheta, r, \kappa)  > 0$ such that
$q = q(\hd, \hat{F}) \geq q_0(\vartheta, r, \kappa) > 0$.
Recall that $\cg^{*,\delta}$ denotes the GOSD.  Let ${\cal U}_p^*$ be the set of retained indices. Viewing it as a subgraph of $\cg^{*,\delta}$,  ${\cal U}_p^*$ decomposes into many components
${\cal U}_p^* = \call^{(1)} \cup \call^{(2)} \ldots \cup \call^{(N)}$.
Recall that $\ty = X'Y$. The following lemma is proved below.
\begin{lemma}  \label{lemma:SASLB}
Given that the  conditions of Lemma \ref{lemma:SAS} hold,   there is a constant $c_1 = c_1(\vartheta, r, \kappa, \gamma, A) > 0$ such that with probability at least $1 - o(1/p)$,  for any component $\call_0 \lhd {\cal U}_p^*$,
$\| \ty^{\call_0} \|^2 \geq 2 c_1 |\call_0| \log(p)$.
\end{lemma}
The remaining part of the proof is similar to that of \citet[Lemma 2.3]{JiJin} so we omit it. We note that however Lemma \ref{lemma:SASLB} is  new and needs a much harder proof.  \qed

\subsubsection{Proof of Lemma \ref{lemma:SASLB}}
First, we  need some notations.   Let $\call_0$ be a component of ${\cal U}_p^*$, and  let $\call_0^{(i)}$, $1 \leq i \leq N_0$, be all connected subgraphs with size $\leq m_0$, listed in the order as in the GS-step, where $N_0$ is an integer that may depend on $(X, Y)$.
For each $1 \leq i \leq N_0$,  let $\call_0^{(i)} = \hat{D}^{(i)} \cup \hat{F}^{(i)}$
be the {\it exactly the same}  partition when we screen $\call_0^{(i)}$ in the $m_0$-stage $\chi^2$-screening of the GS-step. In out list, we only keep $\call_0^{(i)}$ such that $\hd^{(i)} \cap  \call_0 \neq \emptyset$.  Since $\call_0$ is a component of ${\cal U}_p^*$ and $\call_0^{(i)}$ is a connected subgraph, it follows from the way that the $\chi^2$-screening
is designed and the definition of  $\hd^{(i)}$   that
\begin{equation} \label{lemmaSAS3a}
\call_0^{(i)} \subset \call_0,   \qquad \mbox{and} \qquad  \hat{D}^{(i)}  =  \call_0^{(i)} \setminus (\cup_{j = 1}^{i-1} \call_0^{(j)}),     \qquad 1 \leq i \leq N_0,
\end{equation}
and
\begin{equation} \label{lemmaSAS3b}
\mbox{$\call_0 = \hd^{(1)} \cup \hd^{(2)} \ldots \cup \hd^{(N_0)}$ is a partition},
\end{equation}
where $\he^{(1)}$ is empty.

Now, for each $1 \leq i \leq N_0$, recall that as long as $G^{\call_0^{(i)}, \call_0^{(i)}}$ is non-singular,  the $\chi^2$-test score in  GS  is $T(Y, \hd^{(i)}, \he^{(i)}) = T(Y, \hd^{(i)}, \he^{(i)};  \call_0^{(i)}, X, p, n)  = (\ty^{\call_0^{(i)}})' (G^{\call_0^{(i)}, \call_0^{(i)}})^{-1}  \ty^{\call_0^{(i)}} -  (\ty^{\he^{(i)}} )' (G^{\he^{(i)}, \he^{(i)} })^{-1}  \ty^{\he^{(i)}}$.
By basic algebra and direct calculations, it can be verified that
$T(Y, \hd^{(i)}, \he^{(i)})   =  \|W_i\|^2$,
where $W_i = W(\ty, \hd^{(i)}, \he^{(i)};  \call_0^{(i)}, X, p, n)$ is defined as $W_i = V_i^{-1/2} y_i$, and for short,
$V_i = G^{\hd^{(i)}, \hd^{(i)}} -  G^{\hd^{(i)}, \he^{(i)} } (G^{\he^{(i)}, \he^{(i)}})^{-1} G^{\he^{(i)}, \hd^{(i)}}$,  $y_i = \ty^{\hd^{(i)}} -  G^{\hd^{(i)}, \he^{(i)} } (G^{\he^{(i)}, \he^{(i)}})^{-1} \ty^{\he^{(i)}}$.
At the same time, for a constant $\delta > 0$ to be determined, define
$\tilde{\Omega}$ by
$\tilde{\Omega}(i,j)  = G(i,j) \cdot 1\{ |G(i,j)| \geq \delta \}$.
The definition of $\tilde{\Omega}$ is the same as that of $\Omega^{*,\delta}$, except for that  the threshold $\delta$ would be selected differently.
We introduce a counterpart of $W_i$ which we call  $W_i^*$,
\begin{equation} \label{DefineW2}
W_i^* =   V_i^{-1/2}   y_i^*.
\end{equation}
where
$y_i^* = \ty^{\hd^{(i)}} -  \tilde{\Omega}^{\hd^{(i)}, \he^{(i)} } (\tilde{\Omega}^{\he^{(i)}, \he^{(i)}})^{-1} \ty^{\he^{(i)}}$.
Let $W^* = ((W_1^*)', (W_2^*)', \ldots, (W_{N_0}^*)')'$,  and define $|\call_0| \times |\call_0|$ matrices $H_1$ and $H_2$ as follows:  $H_1$ is a diagonal block-wise matrix where the $i$-th block is $V_i^{-1/2}$,
and
$H_2 = \tilde{H}_2^{\call_0,\call_0}$, where $\tilde{H}_2$ is a $p \times p$ matrix such that  for every component $\call_0$ of ${\cal U}_p^*$, and $\hd^{(i)}$ and $\he^{(i)}$ defined on each component,
$\tilde{H}_2^{\hd^{(i)},  \he^{(i)} } =  -  (\tilde{\Omega})^{\hd^{(i)}, \he^{(i)} } [(\tilde{\Omega})^{\he^{(i)}, \he^{(i)}}]^{-1}$, $\tilde{H}_2^{\hd^{(i)}, \hd^{(i)} } = I_{|\hd^{(i)}|}$,
and that the coordinates of $\tilde{H}_2$ are zero elsewhere. Here $I_k$ stands for $k\times k$ identity matrix.
From the definitions, it is seen that
\begin{equation} \label{lemmaSAS3c}
W^* =  H_1 H_2 \ty^{\call_0}.
\end{equation}

Compared with $W_i$, $W_i^*$ is relatively easier to study, for it induces column-sparsity  of $H_2$. In fact, using \cite[Lemma 2.2, 3.1]{JiJin}, there is an event $A_p$ that depends on $(X, \beta)$ such that $P(A_p^c) \leq o(1/p^2)$ and that
over the event,
for all subset $B$ with size $\leq m_0$,
\begin{equation} \label{lemma:SAS1b}
\|G^{B, B} - \Omega^{B, B} \|_\infty \leq L_p p^{-\kappa/2}.
\end{equation}
The following lemma is proved below.
\begin{lemma} \label{lemma:eigenbound}
Fix $\delta > 0$ and suppose the conditions in Lemma \ref{lemma:SASLB} hold.  Over the event $A_p$,  there is a constant $C > 0$ such that each row and column of $\tilde{H}_2$ has no more than $C$ nonzero coordinates.
\end{lemma}

We are now ready to show Lemma \ref{lemma:SASLB}. To begin with, note that since we accept
$\hd^{(i)}$ when we graphlet-screen $\call_0^{(i)}$ and  $|\hd^{(i)}| \leq m_0$,
\begin{equation} \label{lemmaSAS2a}
\| W_i\|^2  \geq 2 (q_0/ m_0)  |\hd^{(i)}|\log(p).
\end{equation}
At the same time, by basic algebra,  $\|W_i - W_i^*\| \leq \|V_i^{-1/2}\| \| y_i - y_i^*\|$,
and
$\| y_i - y_i^*\| \leq   \| G^{\hd^{(i)}, \he^{(i)} } (G^{\he^{(i)}, \he^{(i)}})^{-1} -( \tilde{\Omega})^{\hd^{(i)}, \he^{(i)} } ((\tilde{\Omega})^{\he^{(i)}, \he^{(i)}})^{-1} \|_\infty \cdot \|\ty^{\he^{(i)}}\|$.
First,  since $\lambda_{m_0}^*(\Omega) \geq C$, it is seen that  over the event $A_p$,
$\|V_i^{-1/2}\| \leq C$. Second, by similar reasons,  it is not hard to see that except for probability $o(p^{-2})$, $\| G^{\hd^{(i)}, \he^{(i)} } (G^{\he^{(i)}, \he^{(i)}})^{-1} -( \tilde{\Omega})^{\hd^{(i)}, \he^{(i)} } ((\tilde{\Omega})^{\he^{(i)}, \he^{(i)}})^{-1} \|_\infty  \leq C \delta^{1-\gamma}$, and $\|\ty^{\he^{(i)}} \|\leq C \sqrt{\log(p)} \leq C \tau_p$. Combining these gives
\begin{equation} \label{lemmaSAS2b}
   \|  W_i - W_i^* \| \leq C \delta^{1-\gamma} \tau_p,
\end{equation}
Inserting this to  (\ref{lemmaSAS2a}), if we choose $\delta$ to be a sufficiently small constant,
$\| W_i^* \|^2\geq \frac{1}{2} \|W_i\|^2 \geq (q_0 / m_0) |\hat{D}^{(i)}| \log(p)$.

At the same time, by definitions, it follows from  $\|V_i^{-1/2}\| \leq C$  that $\|H_1\| \leq  C$.  Also, since over the event $A_p$, each coordinate of $H_2$ is bounded from above by a constant in magnitude, it follows from Lemma \ref{lemma:eigenbound} that $\|H_2\| \leq C$.
Combining this with (\ref{lemmaSAS3b})-(\ref{lemmaSAS3c}), it follows from basic algebra that except for probability $o(p^{-2})$,
$(q_0 / m_0)  |\call_0|  \log(p)  \leq \|W^*\|^2 \leq \| H_1 H_2 \ty^{\call_0}\|^2 \leq C \|\ty^{\call_0}\|^2$,
and the claim follows since $m_0$ is a fixed integer.   \qed

\subsubsection{Proof of Lemma \ref{lemma:eigenbound}}
By definitions, it is equivalent to show that over
the event $A_p$,  each  row and column of   $\tilde{H}_2$  has finite nonzero
coordinates. It is seen that each row of $\tilde{H}_2$ has $\leq m_0$ nonzeros,
so all we need to show is that each column of $\tilde{H}_2$ has finite nonzeros.

 Towards this end, we introduce a new graph $\tilde{\cg} = (V, E)$, where $V= \{1, 2, \ldots, p\}$ and nodes $i$ and $j$ are connected if and only if $|\tilde{\Omega}(i,j)| \neq 0$. This definition is the same as GOSD, except that $\Omega^{*,\delta}$ is substituted by $\tilde{\Omega}$.  It is seen that over the event $A_p$, for any $\Omega \in {\cal M}_p^*(\gamma,c_0,g, A)$,   $\tilde{\cg}$ is $K$-sparse with $K \leq C\delta^{-1/\gamma}$. The key for the proof is to show that
for any $k \neq \ell$ such that $\tilde{H}_2(k,\ell) \neq 0$, there is a path with length $\leq (m_0-1)$ in $\tilde{\cg}$ that connects $k$ and $\ell$.

To see the point, we note that when   $\tilde{H}_2(k, \ell) \neq 0$,  there must be an $i$ such that $k \in \hd^{(i)}$ and $\ell \in \he^{(i)}$.
We claim that there is a path in $\call_0^{(i)}$ (which is regarded as a subgraph of $\tilde{\cg}$) that connects $k$ and $\ell$.
In fact, if $k$ and $\ell$ are not connected in $\call_0^{(i)}$, we can partition $\call_0^{(i)}$ into two separate sets of nodes  such that one contains $k$ and the other contains $\ell$, and two
sets are disconnected. In effect, both the matrix $\tilde{\Omega}^{\hd^{(i)}, \hd^{(i)}}$ and $\tilde{\Omega}^{\hd^{(i)}, \he^{(i)}}$ can be visualized as two by two blockwise matrix, with off-diagonal
blocks being $0$. As a result, it is seen that $\tilde{H}_2(k, \ell)  = 0$. This contradiction shows that whenever $\tilde{H}_2(k, \ell) \neq 0$,
$k$ and $\ell$ are connected by a path in $\call_0^{(i)}$. Since $|\call_0^{(i)}| \leq m_0$, there is a path $\leq m_0-1$ in $\tilde{\cg}$ that connects $k$ and $\ell$ where $k \neq \ell$.

Finally, since $\tilde{\cg}$ is $K$-sparse with $K = C\delta^{-1/\gamma}$, for any fixed $\ell$,  there are at most finite $k$ connecting to $\ell$ by a path with length $\leq (m_0-1)$. The claim follows.  \qed

\subsection{Proof of Theorem \ref{thm:compareGsSsLasso}}
\label{subsubsec:proof:thm:compareGsSsLasso}
Since $\sigma$ is known, for simplicity, we assume $\sigma=1$.
First, consider (\ref{eq:blockGSrate}).   By Theorem \ref{thm:UB} and (\ref{Definerho}), $\rho_{gs}  = \min_{\{(D, F): D \cap F = \emptyset, D \neq \emptyset, D \cup F \subset \{1, 2\}\}}  \rho(D,F; \Omega)$, where we have used that $G$ is a diagonal block-wise matrix, each block is the same $2 \times 2$ matrix. To calculate $\rho(D,F; \Omega)$, we consider
 three cases (a) $(|D|, |F|) = (2,0)$,  (b) $(|D|, |F|)  = (1, 0)$, (c) $(|D|, |F|) = (1, 1)$.  By definitions and direct calculations, it is seen that $\rho(D,F; \Omega) = \vartheta + [(1 - |h_0|) r]/2$ in case (a), $\rho(D,F; \Omega) = (\vartheta + r)^2/(4r)$ in case (b), and   $\rho(D,F; \Omega) = 2 \vartheta +  [(\sqrt{(1 - h_0^2) r} - \vartheta/\sqrt{(1 - h_0^2) r})_+]^2/4$ in case (c). Combining these gives the claim.

Next,  consider (\ref{eq:blockSSrate}).  Similarly, by the block-wise structure of $G$, we can restrict our attention to the first two coordinates of $\beta$, and apply the subset selection to the size $2$ subproblem where the Gram matrix is the $2 \times 2$ matrix
with $1$ on the diagonals and $h_0$ on the off-diagonals.
Fix $q > 0$, and let the tuning parameter $\lambda_{ss}=\sqrt{2q_{ss}\log(p)}$.
Define
$f_{ss}^{(1)}(q)  =  \vartheta+[(\sqrt{r}-\sqrt{q})_{+}]^2$, $f_{ss}^{(2)}(q) =  2\vartheta+[(\sqrt{r(1-h_0^2)}-\sqrt{q})_+]^2$,
and
$f_{ss}^{(3)}(q) = 2 \vartheta + 2 [(\sqrt{r (1 - |h_0|)} - \sqrt{q})_+]^2$,
where $x_+ = \max\{x, 0\}$.
The following lemma is proved below, where the key is to use \citet[Lemma 4.3]{JiJin}.
\begin{lemma}
\label{lemma:RateGivenQ}
Fix $q > 0$ and suppose the conditions in Theorem \ref{thm:compareGsSsLasso} hold.   Apply the subset selection to the aforementioned size $2$ subproblem with $\lambda_{ss} = \sqrt{2 q \log(p)}$.  As $p \goto \infty$, the worst-case Hamming error rate is $L_p p^{-f_{ss}(q)}$, where
$f_{ss}(q) = f_{ss}(q, \vartheta, r, h_0) =   \min
\bigl\{\vartheta +  (1-|h_0|) r /2,     q,   f_{ss}^{(1)}(q), f_{ss}^{(2)}(q), f_{ss}^{(3)}(q)  \bigr\}$.
\end{lemma}
By direct calculations,  $\rho_{ss}(\vartheta, r, h_0) = \max_{\{q > 0\}} f_{ss}(\vartheta, r, h_0)$ and the claim follows.

Last, consider (\ref{eq:blockLassorate}).
The proof is very similar to that of the subset selection, except for that
we need to use   \citet[Lemma 4.1]{JiJin}, instead of \citet[Lemma 4.3]{JiJin}.
For this reason, we omit the proof. \qed

\subsubsection{Proof of Lemma \ref{lemma:RateGivenQ}}
By the symmetry in \eqref{Definess}-\eqref{Definelasso} when $G$ is given by  \eqref{eq:2by2omega}, we only need to consider the case where $h_0  \in [ 0,1)$ and $\beta_1\geq 0$.  Introduce events, $A_0=\{\beta_1=\beta_2=0\}$, $A_1=\{\beta_1 \geq \tau_p,\beta_2=0\}$, $A_{21}=\{\beta_1 \geq \tau_p,\beta_2 \geq \tau_p \}$, $A_{22}=\{\beta_1 \geq \tau_p,\beta_2 \leq - \tau_p\}$, $B_0=\{\hb_1= \hb_2=0\}$, $B_1=\{\hb_1>0,\hb_2=0\}$, $B_{21}=\{\hb_1>0,\hb_2>0\}$ and $B_{22}=\{\hb_1>0,\hb_2<0\}$.
It is seen that the Hamming error
\begin{equation} \label{eq:prf:RategivenQ}
=L_p( I + II + III),
\end{equation}
where $I=P(A_0\cap B_0^c)$, $II=P(A_1\cap B_1^c)$ and $III=P(A_{21}\cap B_{21}^c)+P(A_{22}\cap B_{22}^c)$.

Let  $H$ be the  $2\times2$ matrix
with ones  on the diagonals and $h_0$ on the off-diagonals,   $\alpha = (\beta_1, \beta_2)'$, and
$w = (\ty_1, \ty_2)$, where we recall $\ty=X'Y$.
It is seen that  $w \sim N(H \alpha, H)$.  Write for short $\lambda   =  \sqrt{2 q \log(p)}$.   Define regions on the plane of $(\ty_1,\ty_2)$,  $D_0=\{\max(|\ty_1|,|\ty_2|)>\lambda \ or \ \ty_1^2+\ty_2^2-2h_0\ty_1\ty_2>2\lambda^2(1-h_0^2)\}$, $D_1=\{|\ty_1|<\lambda \ , \ \ty_1<\ty_2 \ or \ |\ty_2-h_0\ty_1|>\lambda\sqrt{1-h_0^2}\}$, $D_{21}=\{\ty_2-h_0\ty_1<\lambda\sqrt{1-h_0^2} \ or \ \ty_1-h_0\ty_2<\lambda\sqrt{1-h_0^2}\}$ and $D_{22}=\{\ty_2-h_0\ty_1>-\lambda\sqrt{1-h_0^2}\ or \ \ty_1-h_0\ty_2>\lambda\sqrt{1-h_0^2} \ or \ \ty_1^2+\ty_2^2-2h_0\ty_1\ty_2<2\lambda^2(1-h_0^2)\}$. Using \cite[Lemma 4.3]{JiJin}, we have $B_0^c=\{(\ty_1,\ty_2)'\in D_0\}$, $B_1^c=\{(\ty_1,\ty_2)'\in D_1\}$, $ B_{21}^c=\{(\ty_1,\ty_2)'\in D_{21}\}$, and $B_{22}^c=\{(\ty_1,\ty_2)'\in D_{22}\}$.
By direct calculation and Mills' ratio, it follows that for all $\mu \in \Theta_p(\tau_p)$,
\begin{equation} \label{lass1}
I =  L_p \cdot (P(N(0,1)>\lambda)+P(\chi_2^2>2\lambda^2))=L_p \cdot p^{-q},
\end{equation}
\begin{equation} \label{lass2}
II \leq L_p \cdot P(N((\tau_p,h_0\tau_p)',H)\in D_1)=L_p \cdot p^{-\vartheta-\min[(\sqrt{r}-\sqrt{q})^2,(1-h_0)r/2,q]},
\end{equation}
and when $\beta_1=\tau_p$ and $\beta_2=0$, the equality holds in \eqref{lass2}.
At the same time, note that over the event $A_{21}$,  the worst case scenario, is where $\beta_1=\beta_2=\tau_p$. In such a  case, $(\ty_1,\ty_2)'  \sim N(((1+h_0)\tau_p,(1+h_0)\tau_p)',H)$. Combining this with Mills' ratio, it follows that for  all $\mu \in \Theta_p(\tau_p)$,
\begin{equation} \label{lass3a}
P(A_{21}\cap B_{21}^c)  =  P((\ty_1,\ty_2)'\in D_{21})
\leq   L_p\cdot p^{-2\vartheta-(\sqrt{r(1-h_0^2)}-\sqrt{q})_+^2},
\end{equation}
and the equality holds when $\beta_1=\beta_2=\tau_p$.
Similarly, note that over the event $A_{22}$, in the worst case scenario, $\beta_1=-\beta_2=\tau_p$. In such a case,  $(\ty_1,\ty_2)'\sim N(((1-h_0)\tau_p,-(1-h_0)\tau_p)',H)$. Combining this with Mills' ratio, it follows that   for all $\mu \in \Theta_p(\tau_p)$,
\begin{equation} \label{lass3b}
P(A_{22}\cap B_{22}^c)  =P((\ty_1,\ty_2)'\in D_{22}) \leq    L_p\cdot p^{-2\vartheta-\min([(\sqrt{r(1-h_0^2)}-\sqrt{q})_+]^2,2\{[\sqrt{r(1-h_0)}-\sqrt{q}]_+\}^2)},
\end{equation}
and the equality holds when $\beta_1=-\beta_2=\tau_p$.
Inserting (\ref{lass1})-(\ref{lass3b}) into (\ref{eq:prf:RategivenQ}) gives the claim.   \qed

\subsection{Lemma \ref{lemma:V} and the proof}
\begin{lemma} \label{lemma:V}
Let $(V_{0j}^*, V_{1j}^*)$ be defined as in \eqref{eq:defineV0V1star}. If the conditions of Theorem \ref{thm:LB} hold,  then
$\max\{|V_{0j}^* \cup V_{1j}^*| \}  \leq  (\vartheta + r)^2/(2 \vartheta r)$.
\end{lemma}
{\bf Proof}. Let $V_0 = \emptyset$ and $V_1 = \{j\}$. It is seen that
$\alpha^*(V_0,V_1;\Omega) = 1$, and $\eta(V_0, V_1;\Omega)   \leq (\vartheta + r)^2/(4 r)$.
Using this and the definitions of $V_{0j}^*$ and $V_{1j}^*$,   $\max\{|V_{0j}^*|, |V_{1j}^*|\} \vartheta  \leq (\vartheta + r)^2/(4r)$ and the claim follows.  \qed




\begin{thebibliography}{99}
\bibitem[Akaike(1974)]{AIC}
{\sc Akaike}, H. (1974). A new look at the statistical model identification. {\it IEEE Trans. Automatic Control},  {\bf 19}(6),   716--723.

\bibitem[Bajwa {\it et~al}(2007)]{Nowak2007}
{\sc Bajwa}, W. U., {\sc Haupt}, J. D., {\sc Raz}, G. M., {\sc
Wright}, S. J. and {\sc Nowak}, R. D. (2007). Toeplitz-structured
compressed sensing matrices. {\it Proc. SSP' 07, Madison, WI, Aug.
2007}, 294--298.


\bibitem[Bickel and Levina(2008)]{Bickel2008a}
{\sc Bickel}, P. and {\sc Levina}, E. (2008).  Covariance regularization by thresholding.  {\it Ann. Statist.},  {\bf 36}(6), 2577--2604.



\bibitem[Cai {\it et~al}(2010)]{CaiLuo2010}
{\sc Cai}, T.,  {\sc Liu}, W. and {\sc Luo}, X. (2010).
A constrained $\ell^1$ minimization approach to sparse precision matrix estimation.
{\it J. Amer. Statist. Assoc.},  {\bf 106}(494), 594-607.

\bibitem[Candes and Plan(2009)]{Candes}
{\sc Candes}, E. and {\sc Plan},  Y.  (2009).  Near-ideal model selection by $\ell^1$-minimization. {\it Ann. Statist.}, {\bf  37}(5),  2145--2177.

\bibitem[Candes and Tao(2007)]{CandesTao}
{\sc Candes}, E. and {\sc Tao}, T. (2007). The Dantzig selector: statistical estimation when $p$ is much larger than $n$ (with discussion). {\it Ann. Statist.},   {\bf 35}(6), 2313--2404.

\bibitem[Chen {\it et~al}(1998)]{Chen}
{\sc Chen}, S., {\sc Donoho}, D.  and {\sc Saunders}, M. (1998). Atomic decomposition by basis pursuit. {\it SIAM J. Sci. Computing},  {\bf 20}(1), 33--61.

\bibitem[Dinur and Nissim(2003)]{Nissim}
{\sc Dinur}, I.  and {\sc Nissim}, K. (2003).  Revealing information while preserving privacy.  {\it Proceedings of the twenty-second ACM SIGMOD-SIGACT-SIGART symposium on principles of database systems}, 202--210,  ACM Press.

\bibitem[Donoho(2006a)]{Donoho2006a}
{\sc Donoho}, D. (2006a). For most large underdetermined systems of linear equations the minimal $\ell^1$-norm solution is also the sparsest solution. {\it Comm. Pure Appl. Math}.,  {\bf 59}(7),  907--934.


\bibitem[Donoho(2006b)]{Donoho2006b}
{\sc Donoho}, D. (2006b). Compressed sensing. {\it IEEE Trans. Inform. Theory},
{\bf 52}(4), 1289--1306.


\bibitem[Donoho and Huo(2001)]{DonohoHuo}
{\sc Donoho}, D. and {\sc Huo}, X. (2001).
Uncertainty principle and ideal atomic decomposition.
{\it      IEEE Trans. Inform. Theory}, {\bf 47}(7), 2845-2862.

\bibitem[Donoho and Jin(2004)]{DonohoJin2004}
{\sc Donoho}, D.  and {\sc Jin}, J. (2004).  Higher criticism for  detecting sparse heterogeneous mixtures. {\it Ann. Statist.}, {\bf 32}(3), 962--994.


\bibitem[Donoho and Jin(2008)]{DonohoJin2008}
 {\sc Donoho}, D. and {\sc Jin}, J. (2008).  Higher criticism thresholding: optimal feature selection when useful features are rare and weak. {\it Proc. Natl. Acad. Sci.}, {\bf 105}(39), 14790--14795.


\bibitem[Donoho and Stark(1989)]{UP2}
{\sc Donoho}, D. and {\sc Stark}, P (1989). Uncertainty
principles and signal recovery. {\it SIAM J. Appl. Math.}, {\bf 49}(3), 906-931.



\bibitem[Efron {\it et~al}(2004)]{LARS}
{\sc Efron}, B.,  {\sc Hastie}, H.,   {\sc Johnstone}, I.  and {\sc Tibshirani}, R. (2004).
Least angle regression (with discussion). {\it Ann. Statist.}, {\bf 32}(2), 407-499.

\bibitem[Fan {\it et~al} (2011)]{Hanetal2011}
{\sc Fan}, J.,  {\sc Han}, X.  and {\sc Gu}, W. (2012).
Estimating false discovery proportion under arbitrary covariance dependence.  {\it J. Amer. Statist. Assoc}.,   {\bf 107}(499),  1019-1035.


\bibitem[Fan and Li(2001)]{FanLi}
{\sc Fan}, J. and {\sc Li}, R. (2001). Variable selection via nonconcave penalized likelihood and its oracle properties. {\it J. Amer. Statist. Assoc.}, {\bf 96}(456), 1349--1360.


\bibitem[Fan and Lv(2008)]{FanLv}
{\sc Fan}, J. and {\sc Lv}, J. (2008).
Sure independence screening for ultrahigh dimensional feature space. {\it J. Roy. Soc. Statist. B},  {\bf  70}(5),  849--911.

\bibitem[Foster and George (1994)]{RIC}
{\sc Foster}, D. P. and  {\sc George}, E. I. (1994). The risk inflation criterion for multiple regression. {\it Ann. Statist.},   {\bf  22}(4),  1947--1975.

\bibitem[Friedman {\it et~al}(2008)]{glasso}
{\sc Friedman},  J.,  {\sc Hastie},  T. and {\sc  Tibshirani},   R. (2008).  Sparse inverse covariance estimation with the graphical lasso.
{\it Biostatistics}, {\bf 9}(3), 432--441.

\bibitem[Friedman {\it et~al}(2010)]{Friedman2010}
{\sc Friedman}, J. H., {\sc Hastie}, T. and {\sc Tibshirani}, R. (2010).  Regularization paths for generalized linear models via coordinate descent. {\it J. Statist. Soft. Soft},
{\bf 33}(1). \\
\verb+http://cran.r-project.org/web/packages/glmnet/index.html+

\bibitem[Frieze and Molloy(1999)]{Frieze}
{\sc Frieze}, A.M. and  {\sc Molloy}, M. (1999).   Splitting an expander graph. {\it J.   Algorithms}, {\bf 33}(1),   166--172.


\bibitem[Genovese {\it et~al}(2012)]{GJW2011}
{\sc Genovese}, C., {\sc Jin}, J.,  {\sc Wasserman}, L., and {\sc Yao}, Z.   (2012).
A comparison of the lasso and marginal regression. {\it J. Mach. Learn. Res}. {\bf 13}, 2107-2143.


\bibitem[Hopcroft and Tarjan(1973)]{Hopcroft}
{\sc Hopcroft}, J. and {\sc Tarjan}, R. (1973).
Algorithm 447: efficient algorithms for graph manipulation. {\it Commun. ACM}, {\bf 16}(6), 372-378.

\bibitem[Horn and Johnson(1990)]{Horn1990}
{\sc Horn}, R. and {\sc Johnson} (1990).
Matrix analysis.  {\it  Cambridge University Press}.

\bibitem[Ingster{\it et~al} (2009)]{Ingster2009}
{\sc Ingster}, Y.,  {\sc Pouet}, C.  and {\sc Tsybakov}, A. (2009).  Classification of sparse high-dimensional vectors.     {\it Phil. Trans. R. Soc. A},   {\bf 367}(1906),  4427-4448.


\bibitem[Ising(1925)]{Ising}
{\sc Ising}, E. (1925).
 A contribution to the theory of ferromagnetism.
{\it Z.Phys.}, {\bf 31}(1), 253-258.


\bibitem[Ji and Jin(2011)]{JiJin}
{\sc Ji}, P. and {\sc Jin}, J. (2011).
UPS delivers optimal phase diagram in high dimensional variable selection.
{\it Ann. Statist.},   {\bf 40}(1),  73-103.

\bibitem[Jin(2009)]{PNAS}
{\sc Jin}, J.  (2009).  Impossibility of successful classification when useful features are rare and weak. {\it Proc. Natl Acad. Sci.},  {\bf 106}(22),  8859-8864.


\bibitem[Li and Li(2011)]{LiLi}
{\sc Li}, C. and {\sc Li}, H. (2011).
Network-constrained regularization and variable selection for analysis of genomic data.
{\it Bioinformatics}, {\bf 24}(9), 1175-1182.


\bibitem[Meinshausen and Buhlmann(2006)]{Meinshausen}
{\sc Meinsausen}, N. and {\sc Buhlmann}, P. (2006). High dimensional graphs and variable selection with the Lasso. {\it Ann. Statist.},  {\bf 34}(3), 1436--1462.

\bibitem[Meinsausen and Buhlmann(2010)]{MeinB2010}
{\sc Meinsausen}, N. and {\sc Buhlmann}, P. (2010). Stability Selection (with discussion).
{\it  J. Roy. Soc. Statist. B},  {\bf 72}(4),   417-473.

\bibitem[Meinsausen and Rice (2006)]{MeinshausenRice}
{\sc Meinsausen}, N. and {\sc Rice}, J. (2006).   Estimating the proportion of false null hypotheses among a large number of independently tested hypotheses.
 {\it Ann. Statist.}, {\bf 34}(1), 323-393.


\bibitem[Pearl(2000)]{Pearl}
{\sc Pearl}, J. (2000). {\it  Causality: models, reasoning, and inference}.  Cambridge University Press.

\bibitem[Peng {\it et~al}(2009)]{Pengetal2009}
{\sc Peng}, J.,   {\sc Wang}, P.,  {\sc Zhou}, N. and {\sc Zhu}, J. (2010).
Partial correlation estimation by joint sparse regression model.  {\it J. Amer. Statist. Assoc.},  {\bf 104}(486),735-746.


\bibitem[Schwarz(1978)]{BIC}
{\sc Schwarz}, G. (1978). Estimating the dimension of a model. {\it Ann. Statist.},   {\bf 6}(2), 461--464.



\bibitem[Tibshirani(1996)]{Tibshirani}
{\sc Tibshirani}, R.  (1996). Regression shrinkage and selection via the lasso. {\it J. Roy. Statist. Soc. B}, {\bf 58}(1),  267--288.



\bibitem[Wainwright(2009)]{Wainwright}
{\sc Wainwright},  M. (2009).  Sharp threshold for high-dimensional and noisy recovery of sparsity using $\ell_1$ constrained quadratic programming (Lasso). {\it IEEE Trans. Inform. Theory}, {\bf 52}(5), 2183--2202.


\bibitem[Wang {\it et al}(2013)]{Wang}
{\sc Wang}, Z., {\sc Liu}, H., and {\sc Zhang}, T. (2013). Optimal computational and statistical rates of convergence for sparse nonconvex learning problems. {\it ArXiv e-print}, arXiv:1306.4960.

\bibitem[Wasserman and Roeder(2009)]{Wasserman}
{\sc Wasserman}, L. and {\sc Roeder}, K. (2009).  High-dimensional variable selection.  {\it Ann. Statist.},  {\bf 37}(5), 2178--2201.


\bibitem[Ye and Zhang(2010)]{Zhang}
{\sc Ye}, F. and {\sc Zhang}, C-H. (2010).  Rate minimaxity of the Lasso and Dantzig selection for the $\ell_q$ loss in $\ell_r$ balls. {\it J. Mach. Learn. Res.},  {\bf 11}, 3519-3540.


\bibitem[Zhang(2010)]{Zhang10}
{\sc Zhang}, C.-H. (2010). Nearly unbiased variable selection under minimax
concave penalty. {\it Ann. Statist.},  {\bf 38}(2), 894-942.

\bibitem[Zhang and Zhang(2010)]{ZhangZ11}
{\sc Zhang}, C.-H. and {\sc Zhang}, T. (2011). A general theory of concave regularization for high
dimensional sparse estimation problems. {\it Statist. Sci.},   {\bf  27}(4),  576-593.

\bibitem[Zhang(2011)]{ZhangT11a}
{\sc Zhang, T.} (2011). Adaptive forward-backward greedy algorithm for learning sparse representations.
{\it IEEE Trans. Inform. Theory}.,   {\bf 57}(7),  4689-4708.

\bibitem[Zhao and Yu(2006)]{YuB}
{\sc Zhao}, P. and {\sc Yu}, B. (2006). On model selection consistency of LASSO. {\it J. Mach. Learn. Res.},   {\bf 7},   2541--2567.

\bibitem[Zou(2006)]{Zou}
{\sc Zou},  H. (2006).  The adaptive lasso and its oracle properties. {\it J. Amer. Statist. Assoc.}, {\bf 101}(476), 1418--1429.


\end{thebibliography}

\medskip\noindent
{\bf Acknowledgments}: We thank   Pengsheng Ji and Tracy Ke for help and pointers. JJ and QZ were partially supported by NSF  CAREER award DMS-0908613.  QZ was also partially supported by the NIH Grant P50-MH080215-02. CHZ was partially supported by the NSF Grants DMS-0906420, DMS-1106753 and NSA Grant H98230-11-1-0205. The research was supported in part by the computational resources on PittGrid.


\end{document}